\documentclass[11pt,fleqn,leqno]{article}
\usepackage{fullpage}
\usepackage{setspace}
\usepackage{titlesec}
\titleformat*{\section}{\large\bfseries}
\titleformat*{\subsection}{\large\bfseries}

\usepackage{amsmath,amsthm,amssymb,latexsym,enumerate,bibentry, nicefrac, graphicx,dutchcal,covington,stmaryrd}
\usepackage{spverbatim}
\usepackage[T1]{fontenc} 

\usepackage[authoryear,sort]{natbib}

\newtheorem{theorem}{Theorem}[section]
\newtheorem{definition}[theorem]{Definition}

\newtheorem{prop}[theorem]{Proposition}
\newtheorem{corollary}[theorem]{Corollary}

\def\phi{\varphi}        
\def\A{\forall}           
\def\E{\exists} 
\def\mand{\, \wedge \, }    
\def\mor{\, \vee \,} 
\def\imp{\rightarrow}

\def\proves{\vdash}

\usepackage{hyperref}
\usepackage[usenames,dvipsnames]{color}
\hypersetup{colorlinks=true, urlcolor=blue,citecolor=Purple}

\begin{document}

\title{On the methodology of informal rigour: \\ set theory, semantics, and intuitionism\thanks{Forthcoming in \textit{Intuitionism, Computation, and Proof: Selected themes from the research of G. Kreisel}, M. Antonutti Marfori and M. Petrolo (editors), Springer.  Please do not cite without permission.}}

\author{Walter Dean\footnote{University of Warwick, \href{mailto:W.H.Dean@warwick.ac.uk}{\texttt{W.H.Dean@warwick.ac.uk}}.} \ \ \& \ Hidenori Kurokawa\footnote{Kanazawa University, \href{mailto:hidenori.kurokawa@gmail.com}{\texttt{hidenori.kurokawa@gmail.com}}.}}

\date{}

\maketitle

\noindent {\footnotesize  \textsl{In his paper Kreisel adumbrates a crucial insight into the nature of mathematics and foundations of mathematics by focusing on the notion of ``informal rigor.'' It seems to me that philosophy of mathematics should pay much more attention to this notion than has been the case.} \hfill \citep[p. 13]{Isaacson2011a}}

\section{Introduction}
\label{intro}

Georg Kreisel is best known amongst philosophers for his defense of the model theoretic analysis of first-order logical validity and his argument that the Continuum Hypothesis (CH) possesses a definite truth value. The former -- which is often referred to as Kreisel's \textit{squeezing argument} -- continues to play a role in debates about logical consequence and validity after its popularization by authors such as \citet{Field1989,Field2008}, \citet{Boolos1985b}, and \citet{Etchemendy1990}.  The latter has been influential in subsequent debates in philosophy of mathematics regarding the role of second-order logic in securing the determinacy of set theoretic truth  -- e.g. \citet{Weston1976}, \citet{Shapiro1985a,Shapiro1991}, \citet{Potter2004}, \citet{Isaacson2011a}.  

The \textit{locus classicus} for both of Kreisel's arguments is his paper `Informal rigour and completeness proofs' \citeyearpar{Kreisel1967b}.  It is fair to say that this paper is well-known.  But it is  also not a substantial exaggeration to say that the sections of \citeyearpar{Kreisel1967b} which contain Kreisel's validity and CH arguments are the \textit{only} part of his otherwise substantial output which has been widely read or appreciated by philosophers.  This is to some extent understandable.   For not only is Kreisel's style sometimes obscure, even the papers which he directed at philosophical audiences are typically informed by mathematical and historical considerations which are not readily appreciated without an understanding of the milieu in which they were composed.  

This is particularly true of \citeyearpar{Kreisel1967b}.   For although this paper is one of several places in which Kreisel juxtaposed his validity and CH arguments, it also provides the most complete exposition of a technique which he referred to as \textsl{informal rigour}.   Kreisel presented these arguments as examples of the technique.   But it is not entirely straightforward to discern what they have in common as applications of a more general method.    Thus although `informal rigour' is often taken to be a hallmark of Kreisel's work, there is little consensus on how this method should be characterized or the sorts of questions it is intended to address.   

The primary goals of this paper will thus be to distill from Kreisel's various expositions a general schema for what we will refer to as an \textsl{informally rigorous argument} and to illustrate in detail how Kreisel's examples conform to it.    In this context we will additionally demarcate what has come to be called a \textsl{squeezing argument} -- as paradigmatically illustrated by Kreisel's validity argument -- and illustrate how such arguments can be understood as a subspecies of Kreisel's more general method.   

We will present our proposed schematization of informal rigour in \S\ref{oninfrig}. In \S\ref{casestudies} we will then present detailed reconstructions of the three primary examples which Kreisel considered in \citeyearpar{Kreisel1967b}.  In addition to his validity and CH arguments, this includes his analysis of the notion of a \textsl{creating subject} which figures in Brouwer's  development of intuitionistic analysis.   Little subsequent attention has been paid to this case as an example of informal rigour.  But it is of particular interest not only because Kreisel's exposition contains a novel mathematical result (whose proof we will reconstruct in \S\ref{gmpsf}) but also because it illustrates the central role which Kreisel assigned to reflection on what he refers to as `new primitive notions' in the operation of informal rigour.  In Appendix \ref{app} we will summarize more briefly Kreisel's proposals about the applicability of informal rigour to the following additional concepts: \textsl{set}, \textsl{nonstandard model}, \textsl{finitist proof}, \textsl{predicative definability}, and \textsl{intuitionistic validity}.  

These examples are illustrative of the many trees in the forest we are about to enter.  But before considering them individually, it will be useful to consider in \S\ref{context} the context in which Kreisel introduced informal rigour, inclusive of his motives for promoting this methodology and his strategy for presenting it to his intended audiences.   In addition to the complexity of Kreisel's relationships with his interlocutors, many aspects of his exposition were also informed by specific mathematical developments to which he either directly contributed or was otherwise an interested party.   

Several of these relate to the other titular concern of \citeyearpar{Kreisel1967b} -- i.e. the significance of \textsl{completeness proofs} in mathematical logic.   About such proofs Kreisel begins by remarking
\begin{quote}
\footnotesize{[Q]uite generally, problems of completeness (of rules) involve informal rigour, at least when one is trying to decide completeness with respect to an intuitive notion of consequence.  \hfill \citeyearpar[p. 139]{Kreisel1967b}}
\end{quote}
By the time he composed \citeyearpar{Kreisel1967b} Kreisel indeed had extensive experience not only with  the technical details of completeness proofs for a variety of logical systems but also with the prior task of distilling notions of consequence and validity from the foundational work of theorists such as Zermelo, Poincar\'e, Hilbert, Brouwer, Heyting, and G\"odel.    It is thus useful to flag the following portions of Kreisel's mathematical work as antecedents to his introduction of the expression `informal rigour': i)  his refinement of Hilbert and Bernays's \citeyearpar{Hilbert1939} arithmetization of G\"odel's \citeyearpar{Godel1929} original completeness proof for classical first order logic \citeyearpar{Kreisel1950,Kreisel1953}; ii) his metamathematical analysis of the completeness of intuitionistic first-order logic \citeyearpar{Kreisel1958,Kreisel1958f,Dyson1961,Kreisel1970}; iii) his investigation on the role of reflection principles in the characterization of finitary mathematics \citeyearpar{Kreisel1958,Kreisel1968a,Kreisel1970a}; iv) his analysis of predicativity in terms of the hyperarithmetical sets \citeyearpar{Kreisel1959,Kreisel1960e,Kreisel1960,Kreisel1962}; and v) his anticipation of the method of forcing in set theory \citeyearpar{Kreisel1958k,Kreisel1961a}.\footnote{The following sources contain discussion of these developments: i) \citep{Smorynski1984}, \citep{Dean2020e}; ii) \citep{Sundholm1983}, \citep{McCarty2008}, \citep{Dean2016c}; iii) \citep{Feferman1995a}, \citep{Dean2015a}; iv) \citep{Feferman2005}, \citep[\S 2.3]{Dean2017b}; v) \citep[\S III]{Kreisel1980} as well as the other sources cited in note \ref{forcingnote} below.}

Such developments may seem far afield to contemporary readers of  `Informal rigour and completeness proofs'.   But perhaps the most interesting aspect of this paper is how Kreisel understood results in mathematical logic to be related to debates in mid-20th century analytic philosophy -- e.g. the possibility of conceptual analysis,  the viability of logical positivism, and the tenability of realism, idealism, and formalism in the foundations of mathematics.    But although in some instances Kreisel's proposals have gone on to have significant downstream effects, it is still difficult to concisely summarize his philosophical contributions.  In fact Kreisel's influence has mostly been felt at one level of remove via his associates and interpreters who have often described his views incompletely or out of mathematical and historical context (often out of expository necessity).

In \S\ref{legacy} we will discuss a source which helps to illustrate the relationship of Kreisel's views with some of his contemporaries -- i.e.  his exchange with Yehoshua Bar-Hillel recorded in the comments and replies published along with \citeyearpar{Kreisel1967b}.    Much of this concerns the comparison of Kreisel's conception of informal rigour relative to the method which Carnap called \textsl{explication}.    Kreisel was clearly dissatisfied that this discussion did not engage more deeply with the mathematical details of his examples.   But it is still a useful waypoint in appreciating the enduring significance of informal rigour in relation to contemporary discussions of philosophical methodology.

\section{Context}  
\label{context}

We will see below that Kreisel took the validity argument to be the most paradigmatic of his case studies in informal rigour.    But it is also clear that he understood the significance of this argument in relation to the case for the determinacy of CH which he which to impress on his audiences in the mid-1960s.  This is attested by the fact that he juxtaposed the two arguments in at least the following four sources:\footnote{The bracketed information gives the bibliographical keys which we will use here for these sources / our best guess as to when the final versions were prepared / and the location in the text of Kreisel's presentation of the validity and CH arguments.}

{\footnotesize
\begin{itemize}
\item[-]  Mathematical logic. In T.L. Saaty, editor, \textit{Lectures on Modern Mathematics}, Vol. III. Wiley, 1965.  [\citeyearpar{Kreisel1965a} / early 1964 / \S 1.8, pp. 111-112, pp. 114-118]\\
\item[-] Mathematical logic: What has it done for the philosophy of mathematics. In R. Schoenman, editor, \textit{Bertrand Russell, Philosopher of the Century}. Allen and Unwin, 1967. [\citeyearpar{Kreisel1967a} / December 1964 /\S 3, pp. 253-262]\\ 
\item[-] Informal Rigour and Completeness Proofs. In I. Lakatos, editor, \textit{Problems in the philosophy of mathematics: Proceedings of the International Colloquium in the Philosophy of Science, London, 1965}. North-Holland, 1967. [\citeyearpar{Kreisel1967b} / July 1966 / \S 1-2, pp. 147-157]\\ 
\item[-] \textit{Elements of Mathematical Logic. Model theory} (with Jean-Louis Krivine). Studies in Logic and the Foundations of Mathematics. North-Holland. [\citeyearpar{Kreisel1967c} / Krivine, chapters 0-5: 1960-1961; Kreisel chapters 6-7 and appendices: 1966 / Appendix A, \S 4, pp. 189-194] \\ 
\end{itemize}}
\noindent The task of reconstructing Kreisel's methodology is made non-trivial by the factors which shaped his presentations in these sources.  It will hence be useful to begin by recording some details about the circumstances of their publication.

Kreisel was 42 years old in 1965 and had already published more than 30 papers in mathematical logic spanning a number of areas.\footnote{It was also in this year that Kreisel took up a permanent position as Professor of Logic and the Foundations of Mathematics at Stanford where he had previously visited in 1958-1959 and 1962-1963 in conjunction with visits hosted by Kurt G\"odel at the Institute for Advanced Study in Princeton.}   This work forms much of the background to \citeyearpar{Kreisel1965a} which corresponds to an entry on mathematical logic which Kreisel was invited to write for Thomas Saaty's three volume series \textit{Lectures in Modern Mathematics}.  

The series was intended to contain up-to-date surveys on topics in contemporary mathematics accessible to a general audience.\footnote{The papers were accompanied by lectures given by the contributors at George Washington University in 1962-1964.  The volumes also contain chapters by Halmos (on Hilbert space), Eilenberg (on algrabic topology), Kaplansky (on Lie algebras), Ahlfors (on quasiconformal mappings), Milnor (on differential topology), Coxeter (on geometry), and Erd\H{o}s (on number theory).}  At 101 pages, Kreisel's paper was by far the longest contribution to the collection.  It contains a particularly detailed section on set theory -- inclusive of a formulation of G\"odel and Cohen's (at that point very recent) independence result for CH as well as a sketch of the method of forcing itself.  This is followed by sections on intuitionistic mathematics (including the constructive interpretation of the logical connectives, Church's Thesis, choice sequences and an axiomatization of intuitionistic analysis),  proof theory (including a description of Hilbert's program, G\"odel's incompleteness theorems, cut elimination, and ordinal analysis), and what Kreisel termed `impredicative analysis' (under which he includes results on second-order arithmetic and hyperarithmetical theory).\footnote{Kreisel explains the omission of model theory and computability theory from his survey (\S1.74, pp. 113-114) by first suggesting that since satisfaction is defined by recursion on formulas and models `recursion theory appears as a branch of model theory'.   But since models are themselves structured sets, he also suggests that model theory -- while a mathematical subject in is own right -- is a `product of set theory'.}   

Although \citeyearpar{Kreisel1965a} is directed at mathematicians, it is in this paper in which Kreisel first introduced a framework which he would repeatedly employ in his later philosophical work.  He begins by distinguishing between what he refers to as \textit{realist} (or \textit{platonist}), \textit{idealist}, and \textit{formalist} conceptions of mathematics.  Although Kreisel's initial characterization of these views was quite broad, he suggested that particular refinements lead respectively to set theoretic realism (as embodied by Zermelo's characterization of the cumulative hierarchy), intuitionism (as embodied by Brouwer and Heyting's constructive understanding of the logical connectives), and finitism (as embodied in Hilbert's finitist consistency program). Kreisel then provides a number of examples which illustrate how specific mathematical developments have grown out of reflection on these foundational conceptions.

Kreisel employed a similar expository strategy in \citeyearpar{Kreisel1967a}, \citeyearpar{Kreisel1967b}, and  \citeyearpar{Kreisel1967c}.  Although the first two of these sources are directed at philosophers rather than mathematicians, they are also the outcomes of other high-profile invitations he received around the same time.  \citeyearpar{Kreisel1967a} appeared in \textit{Bertrand Russell, Philosopher of the Century}.  This volume was edited by Ralph Schoenman and also includes contributions from  A.J. Ayer, Hilary Putnam, W.V.O. Quine, Hans and Maria Reichenbach, and Dana Scott. As its title suggests, the goal of the paper is in some sense the converse of \citeyearpar{Kreisel1965a} -- i.e. rather than attempting to convince mathematicians that technical results in logic grew out of reflection on foundational concerns,  he suggests that philosophers should take interest in such results because they provide insights which bear on broad debates such as the viability of mechanistic views of reasoning or the contrast between realism and idealism. 

It was in this context in which Kreisel first introduced the phrase `informal rigour' (p. 202ff).  We will examine his precise words in \S\ref{words}.  But what he goes on to describe in this initial presentation is a process of first engaging in reflection on how concepts are used -- both informally and in mathematical practice -- to sharpen traditional questions and then citing mathematical results in an attempt to provide decisive answers.  He also suggests here that this technique can be contrasted with the methods preferred by followers of what he refers to as `positivist' and `pragmatist doctrines'. 

Both the details of informal rigour and Kreisel's motives for promoting it come into sharper focus in \citeyearpar{Kreisel1967b}.   This paper is an expanded version of an address which he delivered at the International Colloquium in the Philosophy of Science which took place in July 1965 in London.  The event was organized by Irme Lakatos and featured parallel tracks in philosophy of mathematics, philosophy of science, and inductive logic, with the proceedings of the first track subsequently being published as \citep{Lakatos1967}.  In addition to Kreisel's talk, the philosophy of mathematics track also included talks by Paul Bernays, L\'aszl\'o Kalmar, Stephan K\"orner, Andrzej Mostowski, and Abraham Robinson.\footnote{Kreisel's talk was scheduled at 8:30 PM on 13 July 1965 in a session chaired by Stephen Kleene.}   Most of the participants in this track took Cohen's recent independence results in set theory as one of their central topics.

Much of what Kreisel says in this paper in order to motivate the method of informal rigour in \citeyearpar{Kreisel1967b} is intended to distinguish it from what he takes to be the methods of his envisioned positivist, pragmatist, and formalist interlocutors.  Neither in this paper nor in \citeyearpar{Kreisel1967a} does he attach names to these designations.  However A.J. Ayer, Yehoshua Bar-Hillel, Rudolph Carnap, Michael Dummett, W.V.O. Quine, Patrick Suppes, and Alfred Tarski were all in attendance at the conference.   In \S\ref{legacy} we will also discuss  an exchange between Kreisel and Bar-Hillel which makes particularly clear the extent to which Kreisel understood himself to be championing an unfashionable method akin to `old fashioned' conceptual analysis against a prevailing spirit of relativism.

Although Kreisel frequently alluded to his validity and CH arguments throughout the 1970s and 1980s, the last place he appears to have juxtaposed their presentations in print is in \citeyearpar{Kreisel1967c}.  This is a textbook which Kreisel co-authored with Jean-Louis Krivine.  The first five chapters were written by Krivine based on lectures delivered in Paris in 1960-61.   These provide an exposition of elementary logic and model theory through quantifier elimination.  Kreisel wrote the two final chapters on definability and infinitary logic as well as two lengthy appendices, respectively on the axiomatic method and the foundations of mathematics.   These contain several technical refinements which clarify how Kreisel understood the arguments which we will discuss in \S\ref{casestudies}.

\section{On informal rigour and squeezing arguments} 
\label{oninfrig}

As we have just noted, Kreisel first employed the expression `informal rigour' in \citeyearpar{Kreisel1967a} before making it one of the titular concerns of \citeyearpar{Kreisel1967b}.    But in neither these papers nor elsewhere does he offer what can reasonably be understood as a \textit{definition} of this method.   What he does do, however, is to characterize the origins and aims of informal rigour, compare it to what he takes to be rival methodologies, and provide detailed illustrations of its application.  

Taken together, Kreisel's examples might be understood to characterize informal rigour with sufficient precision to allow us to ``know it when we see it''.   But after recording some of Kreisel's own characterizations in  \S\ref{words}, our main  goal in this section will be to suggest that it is possible to go a bit further, essentially by turning Kreisel's methodology on itself.  In particular, we will suggest in \S\ref{3irs} that by reflecting further on the case studies he presents it is possible to provide a formal schematic characterization of what Kreisel intended by an \textsl{informally rigorous argument}.  In \S\ref{3ss} we will illustrate how so-called \textsl{squeezing arguments} can be understood as a subcase of this schematization.  

By electing to focus on informally rigorous \textsl{arguments} we will thus largely pass over the possibility that informal rigour can also be understood more broadly to  encompass a method by which axioms characterizing fundamental \textsl{mathematical structures} are discovered or otherwise identified.  Such an interpretation is suggested by several passages which we will consider below and has also gone on to inform the influential accounts of \citet{Shapiro1985a,Shapiro1991} and (more explicitly) \citet{Isaacson2011a}.   In \S\ref{set} we will sketch Kreisel's specific proposal that Zermelo's \citeyearpar{Zermelo1930a} axiomatization of set theory should be understood as an attempt to characterize the `intuitive notion of the cumulative type structure' via informal rigour.   The issues flagged there exemplify why we have decided not to develop this case further in parallel to the other examples of  which he developed in greater detail.

\subsection{In Kreisel's words}
\label{words}

In a late paper otherwise outside the scope of the present survey, Kreisel provides the following characterization of what he took to be the lineage of informal rigour:\footnote{In addition to \citeyearpar{Kreisel1987a}, Kreisel also returned to discuss informal rigour in passing in his retrospective accounts \citeyearpar{Kreisel1987b} and \citeyearpar{Kreisel1989}.}
\begin{quote}
{\footnotesize  `[I]nformal rigour' $\ldots$ is a venerable ideal in the broad tradition of analysing precisely common notions or, as one sometime says, notions implicit in common reasoning.  \hfill \citep[p. 499]{Kreisel1987a}}
\end{quote}
This passage makes clear that Kreisel understood informal rigour as descending from -- but not necessarily identical to -- what has traditionally been called \textsl{conceptual analysis}.   Both points are evident in the following passage in which he originally introduced the expression in
\citeyearpar{Kreisel1967a}:
\begin{quote}
\footnotesize{
\textit{Successes of mathematical logic} Time and again it has turned out that traditional notions in philosophy have an essentially unambiguous formulation when one thinks about them  $\dots$ [A]lso, when so formulated by essential use of mathematical logic, they have non trivial consequences for the analysis of mathematical experience. This discovery conflicts with one's naive impression: for, a first examination of the traditional notions almost always reveals some unexpected ambiguities, and the shock leads one to suppose that further examination might produce an \textit{endless} chain of ambiguities: in other words, that there is nothing behind these notions. Instead, in many cases relatively few basic distinctions were enough to get decisive results $\ldots$ Among them are the well known cases
\begin{enumerate}[(i)]
\item the notion of mechanical process, its stability in the sense that apparently different formulations lead to the same results $\ldots$
\item the notion of aggregate which is analysed by means of the hierarchy (theory) of types, and, of course 
\item the notions of logical validity and logical inference which are analysed $\ldots$ by means of (first order) predicate logic $\ldots$
\end{enumerate}

Besides their intrinsic interest, the results are important as \textit{object lessons}: once one has seen the simple considerations in \S3(a) concerning G\"odel's completeness theorem \textit{one cannot doubt the possibility of philosophical proof} or, as one might put it, of \textit{informal rigour} $\ldots$  \hfill \citeyear[pp. 202]{Kreisel1967a}}
\end{quote}

Here Kreisel comes as close as he ever does to providing a definition of informal rigour by likening it to `philosophical proof'. Of course such an expression might well be taken to embroil us in just the sort of ambiguities Kreisel envisions -- e.g. What distinguishes a `philosophical' proof from a `mathematical' one?  What axioms and rules are allowed?  How are they justified? But Kreisel's introduction of this term is directly linked in the text to the section of \citeyearpar{Kreisel1967a} in which he originally formulates his validity argument  (pp. 253-255).   As we will discuss further in \S\ref{validity},  he appears to have regarded this argument as a particularly paradigmatic example of such a `proof'.

Kreisel provides a more general characterization of the goals and methods of informal rigour in the passage which begins \citeyearpar{Kreisel1967b}:
\begin{quote}
\footnotesize{
It is a commonplace that formal rigour consists in setting out formal rules and checking that a given derivation follows these rules; one of the more important achievements of mathematical logic is Turing's analysis of what a formal rule is. Formal rigour does not apply to the discovery or choice of formal rules nor of notions; neither of basic notions such as set in so-called classical mathematics, nor of technical notions such as \textit{group} or \textit{tensor} product $\ldots$

The `old fashioned' idea is that one obtains rules and definitions by analyzing intuitive notions and putting down their properties. This is certainly what mathematicians thought they were doing when defining length or area or, for that matter, logicians when finding rules of inference or axioms (properties) of mathematical structures such as the continuum $\ldots$  What the `old fashioned' idea assumes is quite simply that the intuitive notions are \textit{significant}, be it in the external world or in thought (and a \textsl{precise} formulation of what is significant in a subject is the result, not a starting point of research into that subject).
Informal rigour wants 
\begin{enumerate}[(i)]
\item to make this analysis as precise as possible (with the means available), in particular to eliminate doubtful properties of the intuitive notions when drawing conclusions about them; and
\item to extend this analysis, in particular not to leave undecided questions which can be decided by full use of evident properties of these intuitive notions.
\end{enumerate}
Below the principal emphasis is on intuitive notions which do not occur in ordinary mathematical practice (so-called new \textit{primitive notions}), but lead to new axioms for current notions. \\  \hspace*{1ex}  \hfill \citeyearpar[pp. 138-139]{Kreisel1967b}  
}
\end{quote}

The characterization (i)-(ii) taken together with Kreisel's analogy to `philosophical proof' will serve as the basis of the schematization of informally rigourous arguments we will propose in \S\ref{3irs}.  Before turning to this, however, it will be useful to consider both how Kreisel understood the distinction between informal rigour and formal rigour and also why he took it to be important to stress this distinction in the immediate context of \citeyearpar{Kreisel1967a,Kreisel1967b}.   

Kreisel was careful to distinguish informal rigour as a \textit{method} from what he refers to as `formal rigour' and also the latter from the \textsl{doctrine} he refers to as `formalism'.  A succinct characterization of the former is provided at the beginning of the previously cited passage -- i.e. `formal rigour consists in setting out formal rules and checking that a given derivation follows these rules'.  In order for such a description to provide a definite means of distinguishing formal and informal rigour requires that the notion of `formal rule' can itself be given a precise definition.   But as the following passage makes clear, it is evident that Kreisel understood the classical analyses of the notion of \textsl{mechanical} (or \textsl{effective}) \textsl{procedure} -- and in particular that of Turing -- to have already provided such an analysis:\footnote{Further evidence to this effect is provided by Kreisel's adoption at the beginning of his section on proof theory in \citeyearpar[p. 149]{Kreisel1965a} of Smullyan's \citeyearpar{Smullyan1961} characterization of a `formal system' -- i.e. essentially what we now call a \textsl{recursively axiomatizable theory} together with a precise definition of derivability with such a system.}

\begin{quote}
{\footnotesize
The possible ambiguity concerning the notion of mechanical procedure $\ldots$ has been investigated in the theory of recursive functions. The principles to be used in analysing such a notion as mechanical process are necessarily more delicate than those used in [analyzing validity], but Turing's analysis is quite convincing. The notion is certainly stable for quite a spectrum of alternatives: the mathematical results which can be used to establish this are strong closure properties of the class of recursive functions.   Thus the basic concept of the formalist doctrine, that of formal system, has not only a clear meaning, but a precise extension $\ldots$ \hfill \citeyearpar[p. 227]{Kreisel1967a}}
\end{quote}

 At least at the time of the sources we are considering, Kreisel not only appears to have accepted what we now call \textsl{Church's Thesis} but in fact employed it as a \textsl{presupposition} in his characterization of  the view he calls `formalism'.\footnote{Although Kreisel states that the notion of mechanical procedure is amenable to informal rigour, in the central sources considered here he did not supply a more extended argument to this effect but rather appears to have assumed that his audience would accept that \citet{Turing1936} had already provided such an analysis.   See rescpetively \citep{Odifreddi1996} and \S\ref{intval} below for an account of the complexities involved in Kreisel's later views about Church's Thesis -- in particular in regard to the exigencies of arguments via `equivalent definitions' -- and its role within intuitionistic mathematics.}   This is evident from the following passage:

{\footnotesize
\begin{quote}
Formalism can be considered at all only because of the discovery of formalization [and can be] sharpened $\ldots$ because mathematical notions cannot only be expressed axiomatically but the relation $\ldots$ $A$ is a logical consequence of $B$, can be defined by means of purely mechanical rules $\ldots$ (formalization of predicate logic); and now the axioms are not interpreted as true about abstract objects, but the whole deductive system consisting of axioms and rules of inference is regarded as a compact description $\ldots$ of the outward (syntactic) forms of mathematical language, separated from their meaning and informal uses. The most important point of the formalist doctrine is this: all questions which go beyond such elementary acts of recognition are regarded as outside mathematics. \hfill \citeyearpar[p. 224-225]{Kreisel1967a}
\end{quote}
}

Kreisel's characterizations of formal rigour then build on such a description of formalism -- e.g.
\begin{quote}
\footnotesize{
[P]hilosophical doctrines can shape the whole style of mathematics, as, for instance formalism $\ldots$ shaped modern axiomatic mathematics $\ldots$ Formalism is responsible for the ideal of formal rigour, so much so that a text book like Bourbaki begins with a set of formal rules of inference; this is not very serious because these rules are never mentioned in the later development which shows that the evidence of the proofs in the main text depends on an \textit{understood} notion of logical inference.} \hfill \citeyearpar[p. 210]{Kreisel1967a}
\end{quote}

This passage is also typical of how Kreisel often inveighed against confusing either the subject matter or practice of a  mathematical discipline with the study of the consequences or metamathematical properties of particular axiomatic systems.   As is already evident in \citeyearpar{Kreisel1965a}, Kreisel thought such a tendency had persisted in the wake of the Hilbert program.  But by the time of \citeyearpar{Kreisel1967a,Kreisel1967b} he also believed (with some justification, as we will see below) that such a view had been reinvigorated in the 1960s by Cohen's independence results in set theory.   As the following passage illustrates, it was against this backdrop which Kreisel felt compelled to promote the continued application of informal rigour:
\begin{quote}
\footnotesize{(\S0) \textit{The case against informal rigour} (or: antiphilosophic doctrines). The present conference showed beyond a shadow of doubt that several recent results in logic, particularly the independence results for set theory, have left logicians bewildered about what to do next: in other words, these results do not `speak for themselves' (to these logicians). I believe the reasons underlying their reaction, necessarily also make them suspicious of informal rigour. \hspace*{1ex} \hfill \citeyearpar[p. 140]{Kreisel1967b}}
\end{quote}

The fact that Kreisel spent the first ten pages of \citeyearpar{Kreisel1967b} defending \textsl{just the possibility} of informal rigour against his envisioned formalist, positivist, and pragmatist interlocutors speaks to the intellectual climate in which he understood himself to be working.   This context also appears to have contributed to his rhetorical motives for juxtaposing his CH and validity arguments.  But before turning to such details, our next goal will be to propose a general schema for understanding what Kreisel hope to achieve via informal rigour which is independent of the historical contingencies which led him to originally promote this method.

\subsection{Schematizing informal rigour}
\label{3irs}

The goal of this section it to propose a general schema for understanding what Kreisel meant by an \textsl{informally rigorous argument}.   The possibility of such a characterization naturally suggests itself in light of the case studies he presented in \citeyearpar{Kreisel1967b} and elsewhere.   And while providing such a characterization is itself nothing more than a straightforward exercise in formalization,  it will be of help in framing several observations about Kreisel's method which would be difficult to formulate on the basis of either his general comments or specific examples alone.

On the model we wish to propose, informal rigour is understood as a method for answering questions about one or more \textsl{common} (or `intuitive' or `traditional') \textsl{concepts}  $\mathcal{C}_1,\ldots, \mathcal{C}_k$ which are either left open by our current understanding or in regard to which additional justification is sought for a given answer.   What distinguishes common concepts from what we will refer to as \textsl{precise concepts} is that the latter but not the former have accepted mathematical definitions of the sort which could be formalized in an appropriate axiomatic theory were we to care to do so.\footnote{The use of the terms `common' and `precise'  apply naturally to several of the example we will consider in \S\ref{ch}.  But this application becomes strained in the case of Kreisel's CH argument wherein our scheme suggests that notion expressible in the language of first-order set be regarded as `common' where those expressible only in second-order logic be regarded as `precise'.   Our use of these terms should thus not be endowed with all of their (genuinely) common connotations.}   At the same time, common concepts are not only in everyday use, but they too will often be mathematical in character.  Thus prior to the application of informal rigour, it will still typically be the case that there is an accepted body of propositions which relate the common concepts $\mathcal{C}_1,\ldots, \mathcal{C}_k$ not just to one another but also potentially to some family of precise concepts $\mathcal{P}_1,\ldots,\mathcal{P}_{m}$.  We will refer to these as \textsl{constitutive principles} for $\mathcal{C}_1,\ldots, \mathcal{C}_k$.   In such situations, there are two sorts of questions to which an informally rigorous argument may be addressed:
\begin{enumerate}[1)]
\item Is a currently unsettled proposition $\Phi$ formulated using either common or precise concepts true or false?  
\item Is it possible to define one of the common concepts $\mathcal{C}_i$ by constructing a coextensive property as a complex of the precise concepts $\mathcal{P}_1,\ldots,\mathcal{P}_{m}$?   
\end{enumerate}

In order to describe the structure of an informally rigorous argument, it will be convenient to switch from a description given in terms of concepts and propositions to one given in terms of predicates and sentences.\footnote{Although Kreisel speaks frequently of `intuitive notions' and `informal concepts', the primary sources considered here provide little direct insight into what sort of theory of concepts or propositions he might have preferred.   But at the same time,  a transition from concepts to some sort of linguistic representations is presumably necessary in order to substantiate Kreisel's analogy between informal rigour and `philosophical \textsl{proof}'.   For this reason, many of Kreisel's applications of informal rigour raise the same issues about the interaction between conceptual analysis and formal derivation which were brought into focus by Frege's famous exchange with Hilbert about the nature and content of geometrical axioms (see, e.g., \citealp{Blanchette1996}). Such a parallelism is particularly apparent in the context of Kreisel's CH argument wherein he can be understood as maintaining a maintaining a position similar to Frege's (i.e. that set theory has a fixed subject matter which determines the cardinality of the continuum) against the formalist position he attributed to Cohen and Robinson (which can be likened to Hilbert's perspective on geometry).  But other than calling attention to these issues, we will make no attempt to reconstruct Kreisel's background views about the nature of concepts here beyond the briefly remarking on his (apparent) commitment to conceptual realism which emerges in light of his exchange with Bar-Hillel in \S\ref{legbar} below (see also \citeyear[pp. 204-205]{Kreisel1989}).}   To this end, we will suppose that the common concepts $\mathcal{C}_1,\ldots, \mathcal{C}_k$ are expressed by primitive predicates $C_1(\vec{x}),\ldots,C_k(\vec{x})$ and the precisely defined concepts $\mathcal{P}_1,\ldots,\mathcal{P}_{m}$ involved in the constitutive principles for the former are expressed as either primitive or complex formulas $\pi_1(\vec{x}), \ldots, \pi_n(\vec{x})$ of a mathematical language (where it is allowed that both $C_i$ and $\pi_i$ can vary in both arity and logical type).  We now consider the languages $\mathcal{L}_{C} =  \{C_1,\ldots,C_k\}$ of the common concepts, $\mathcal{L}_P$ of the mathematical theory in which $\pi_1(\vec{x}), \ldots, \pi_m(\vec{x})$ are defined as well as the \textsl{joint} language $\mathcal{L}_{J} = \mathcal{L}_C \cup \mathcal{L}_P$.   We additionally suppose that the constitutive principles for $\mathcal{C}_1,\ldots, \mathcal{C}_k$ have been formulated as a set of $\mathcal{L}_J$-sentences $\Gamma_1 = \Gamma^C_1 \cup \Gamma^{J}_1$ where $\Gamma^C_1$ consists of $\mathcal{L}_C$-sentences expressing relations among the common concepts themselves and $\Gamma^{J}_1$ consists of $\mathcal{L}_{J}$-sentences expressing principles which bridge between the common and precise concepts.  

With this framework in place, it is now straightforward to schematize what we take Kreisel to have intended by an \textsl{informally rigorous argument}:
\begin{example}[(IR)] \textbf{Informal rigour schema}
\label{irs}
\begin{enumerate}[I)]
\item An informal exploration of  $\mathcal{C}_1,\ldots, \mathcal{C}_n$ is undertaken which may lead to the discovery of three additional sets of principles $\Gamma^C_2, \Gamma^J_2$ and $\Gamma^K_2$ of the following sorts:
\begin{enumerate}[a)]
\item $\Gamma^C_2$ is comprised of $\mathcal{L}_C$-sentences stating additional relations between the common concepts  themselves. 
\item $\Gamma^J_2$ is comprised of $\mathcal{L}_J$-sentences stating additional \textsl{bridging principles} between the common and precise concepts.
\item $\Gamma^K_2$  contains principles stated in a \textsl{Kreiselian language} $\mathcal{L}_K = \mathcal{L}_J \cup \mathcal{L}_N$ where the latter consists of predicates $N_1,\ldots, N_n$ intended to denote \textsl{novel} primitive concepts $\mathcal{N}_1,\ldots,\mathcal{N}_n$ not involved in the statement of the original constitutive principles $\Gamma_1$.    The statements in  $\Gamma^K_2$ will thus formalize the relations which the novel concepts bear to themselves as well as to the common and precise concepts.   
\item We form $\Gamma_2 = \Gamma_1 \cup  \Gamma^C_2 \cup \Gamma^J_2 \cup \Gamma^K_2$ (allowing that any of these sets may be empty). 
\item  Since the precise concepts $\mathcal{P}_1,\ldots,\mathcal{P}_m$ will typically be mathematical in character, we adjoin a set of $\mathcal{L}_P$-axioms $\mathsf{T}_P$ appropriate to their subject matter -- e.g. those of arithmetic, analysis, or set theory -- to form a \textsl{Kreiselian theory} $\mathsf{T}_K = \mathsf{T}_P \cup \Gamma_2$.   
\end{enumerate}
\item Reasoning in $\mathsf{T}_K$ we now attempt to answer a question of type 1) or 2) by demonstrating one of the following:
\begin{itemize}
\item[1$'$)] $\mathsf{T}_K \vdash \varphi$ or $\mathsf{T}_K \vdash \neg \varphi$ where $\varphi$ is an $\mathcal{L}_C$- or $\mathcal{L}_P$-sentence expressing the proposition $\Phi$ about common or precise concepts we hoped to answer.
\item[2$'$)] $\mathsf{T}_K \vdash \forall \vec{x}(C_i(\vec{x}) \leftrightarrow \delta(\vec{x}))$ where $\delta(\vec{x})$ is an $\mathcal{L}_P$-formula expressing a property which gives a precise definition for the common concept $\mathcal{C}_i$ we hoped to provide.
\end{itemize}
\end{enumerate}
\end{example}
\noindent We will illustrate in \S\ref{3irs} how squeezing arguments can be understood as a special case of this schema before considering Kreisel's specific application of informal rigour in \S\ref{casestudies}.   First, however, a few words are in order about the general setup we have proposed.

We have assumed above that for each set of common concepts $\mathcal{C}_1,\ldots,\mathcal{C}_n$ to which we wish to apply informal rigour it will be possible to identify a set of constitutive principles $\Gamma_1$ which are \textsl{accepted} before an argument of type IR is attempted.   As we will see in \S\ref{casestudies}, this is plausible for the sort of mathematical concepts in which Kreisel was most interested.   But even in this case, opponents of the `old fashioned' method of conceptual analysis  may be prone to object that inasmuch as such principles correspond to `conceptual' or `analytic' truths about the concepts $\mathcal{C}_i$, any candidate principle will inevitably be contentious.   In fact according to the most ardent of Kreisel's envisioned critics, we may even have $\Gamma_1 = \emptyset$ for many common concepts.

The genuinely `informal' aspects of informal rigour which come into play at step I of our reconstruction and are presumably intended to speak to this concern.   For it is at this stage where Kreisel proposed that we employ `the idea of pushing a bit farther than before the analysis of the intuitive notions' so as to make them `as precise as possible' and to extend them `so as not to leave undecided questions which can be decided by full use of [their] evident properties' \citeyearpar[pp. 138-139]{Kreisel1967b}.  

Of course such adjurations fail to constitute explicit instructions for \textit{how} to carry out the envisioned process of conceptual reflection in a manner which would have satisfied Kreisel's envisioned interlocutors.  But it is in this regard that he calls attention to the importance not just of setting down or deriving connections between common and recognized precise concepts (steps Ia,b) but also of discovering `\textsl{new primitive notions}' \citeyearpar[p. 139]{Kreisel1967b} by which they may also be  related (step Ic).    Such discoveries presumably cannot be made on the basis of formal rigour alone as our previously formalized theories will contain neither expressions denoting the novel concepts nor axioms describing them.  It is thus also in this regard in which Kreisel was particularly critical of his interlocutors for failing to make additional use of what he regarded as the `evident' (but as yet untapped) properties of various common concepts.\footnote{\citet{Kreisel1970a} would later describe the sort of procedure he envisioned as one of reflection on principles `implicit' in common concepts.   We will discuss the relevance of such comments to his characterization of the concept \textsl{finitist proof}  in \S\ref{finproof}.}   To the extent that this process is successful, Kreisel was thus potentially in a position to argue that $\Gamma_2 \neq \emptyset$ even in cases where his interlocutors will refuse to initially recognize the existence of any `accepted' constitutive principles.  

The basis of Kreisel's \citeyearpar{Kreisel1967a} analogy between informal rigour and `philosophical proof' comes into sharper focus at step II.   Of course nothing precludes the informal reflection on concepts called for at step I as itself taking the form of a deduction from the statements $\Gamma^C_1$ which refer to common principles  alone.   But it is striking that each of the case studies which we will consider below is mediated at this stage by non-trivial \textsl{mathematical theorems} which Kreisel either derived himself or sought to repurpose in the service of informal rigour.\footnote{Some specific examples which we will consider below as follows: G\"odel's Completeness Theorem (\S \ref{comp}), Kreisel's refutation of Generalized Markov's Principle (\S \ref{neggmp}),  Zermelo's Quasi-Categoricity Theorem (\S \ref{zqct}), Kreisel's theorem on the definability of nonstandard models (\S \ref{fund2}), Kleene's theorem on the characterization of the hyperarithmetical sets (\S \ref{kleenethm}), Kreisel theorem (\S \ref{incomp}) on the complexity of the set of intuitionistic validities.}   

In the case of questions of type 1) -- i.e. where we hope to settle an open question about common or precise concepts -- the theorem in question will  be formulated and derived in a mathematical theory $\mathsf{T}_P$, possibly together with $\mathcal{L}_N$-axioms for relevant novel concepts discovered at step Ic.   The theorem may then aid in establishing that the Kreiselian theory $\mathsf{T}_K$ extends the either the constitutive principles $\Gamma_1$ for the common concepts or the or the precise theory $\mathsf{T}_P$  \textsl{non-conservatively}.  As we will see, this is prototypically accomplished by allowing  mathematical reasoning conducted in $\mathsf{T}_P$ can be applied to reach conclusions stated in the languages $\mathcal{L}_C$ or $\mathcal{L}_P$ via the principles $\Gamma^J_2$ and $\Gamma^K_2$ which bridge between common, precise, and novel concepts which cannot be derived from $\Gamma_1$ or $\mathsf{T}_P$ alone.   

\subsection{Squeezing as an instance of informal rigour}
\label{3ss}

Kreisel's best known example of an informally rigorous argument is the validity argument which he presented (with minor variations) in each of the four sources cited at the beginning of \S\ref{context}.  This argument -- to whose specifics we will return in \S\ref{validity} -- is aimed at providing a precise definition of the common concept of \textsl{first-order logical validity}. What we will do in this section is to abstract away from its details so as to provide a schematic characterization of what has come to be known as a \textsl{squeezing argument}.\footnote{Kreisel did not explicitly distinguish this form of argument himself.  The expression `squeezing argument' itself appears to have first been used to describe the validity argument by \citet[p. 121]{Field1989} whose engagement with Kreisel we will return in \S\ref{valrec} below.}  Once such a generalization has been undertaken, it then becomes evident that the methodology of squeezing may potentially be applied to other concepts.\footnote{In \S\ref{preddef} we will discuss how Kreisel's analysis of the concept \textsl{predicative definability} may also be understood in this way.   Others (e.g. \citealp{Smith2011a,Dean2016b}) have also proposed that certain arguments for Church's Thesis (in computably theory) and the Cobham-Edmond Thesis (in computational complexity theory) may also be understood in this manner.}

A squeezing argument is a special case of schema IR wherein a precise definition in the sense of II.2$'$ is provided for a common concept $\mathcal{C}$ by demonstrating that its extension is contained between two precise concepts $\mathcal{P}_n$ and  $\mathcal{P}_w$.   For reasons which will become clear, we will respectively refer to these as the  \textsl{narrow concept} and the \textsl{wide concept}.  Suppose we introduce the primitive predicate symbol $C(x)$ to denote $\mathcal{C}$ and the (primitive or complex) formulas $\pi_n(x)$ and $\pi_w(x)$ to denote $\mathcal{P}_n$ and  $\mathcal{P}_w$.  A squeezing argument can now be schematized as follows:
\begin{example}[(S)]\textbf{Squeezing schema}
\label{ss}
\begin{enumerate}[1)]
\item  An informal argument is given to show that $\forall x(\pi_n(x) \rightarrow C(x))$ -- i.e. that falling under the narrow concept $\mathcal{P}_n$ is a \textsl{sufficient condition} for an item of the appropriate sort to fall under $\mathcal{C}$.
\item  An informal argument is given to show that $\forall x(C(x) \rightarrow \pi_w(x))$ -- i.e. that falling under $\mathcal{C}$ is a \textsl{necessary condition} for an item of the appropriate sort to fall under the wide concept $\mathcal{P}_w$.
\item A mathematical proof is then provided which shows that $\forall x(\pi_n(x) \leftrightarrow \pi_w(x))$ -- i.e. that $\mathcal{P}_n$ and $\mathcal{P}_w$ are \textsl{coextensive}.
\item On the basis of 1) - 3), it is finally concluded that $\forall x(C(x) \leftrightarrow \pi_n(x))$ and also $\forall x (C(x)  \leftrightarrow \pi_w(x))$ -- i.e. that $\mathcal{C}$ is coextensive with both $\mathcal{P}_n$ and $\mathcal{P}_w$ and hence precisely definable by either $\pi_n(x)$ or $\pi_w(x)$.
\end{enumerate}
\end{example}
 
It should now be evident why instances of this pattern of reasoning have become known as `squeezing arguments'.  For note that steps 1 and 2 respectively demonstrate the set-theoretic inclusions $\mathcal{P}_n \subseteq \mathcal{C}$ and $\mathcal{C} \subseteq \mathcal{P}_w$.   The demonstration at step 3 that $\mathcal{P}_n$ and $\mathcal{P}_w$ coincide thus has the effect of defining $\mathcal{C}$ precisely by `squeezing' its extension between those of the narrow and wide concepts.
 
It should also be evident how schema S can be understood as a subschema of IR.  For note that the statements $\forall x(\pi_n(x) \rightarrow C(x))$ and $\forall x(C(x) \rightarrow \pi_w(x))$ required at steps 1) and 2) exemplify the sort of bridging principles between common and precise concepts which may potentially either be counted as constitutive principles for $\mathcal{C}$ or can be obtained as the result of the sort of reflection which Kreisel envisions occurring at stage Ib in an informally rigorous argument.   These statements will thus be contained in the set $\Gamma^J_2$.   The statement $\forall x(\pi_n(x) \leftrightarrow \pi_w(x))$ required at step 3) of schema S is stated in $\mathcal{L}_P$ and typically provable in $\mathsf{T}_P$.   As we will see in \S\ref{casestudies}, this is prototypical of the manner in which mathematical theorems figure in Kreisel's examples.     The final derivation of $\forall x(C(x) \leftrightarrow \pi_n(x))$ or $\forall x (C(x)  \leftrightarrow \pi_w(x))$ -- each of which suffices to provide a precise definition of $\mathcal{C}$ -- called for at step II of schema IR is thus reduced to the trivial matter of putting these implications together in the theory $\mathsf{T}_K$ (which by definition contains both $\Gamma^J_2$ and $\mathsf{T}_P$).

In his general descriptions of informal rigour Kreisel repeatedly stressed how reasoning involving what he referred to as \textsl{new primitive notions} may be of help in answering questions about already familiar domains.    We have codified this in step Ic of schema IR whereby connections between novel, common, and precise concepts are recorded which may then be deductively exploited in conjunction with mathematical principles at stage II.   Note, however, that the statements $\forall x(\pi_n(x) \rightarrow C(x))$, $\forall x(C(x) \rightarrow \pi_w(x))$, and $\forall x(\pi_n(x) \leftrightarrow \pi_w(x))$ which form the major premises of a squeezing argument  do not contain novel vocabulary.   This observation does not rule out that principles involving novel concepts might still figure behind the scenes in their individual derivations.  But this suggests that cases in which squeezing is applicable correspond to ones in which a precise definition of $\mathcal{C}$ is already \textsl{implicit} in either its constitutive principles or other bridging principles to our current inventory of precise concepts which are introduced at step Ib.  For of course what a squeezing argument does is precisely to provide an \textsl{explicit} definition for $\mathcal{C}$ by showing that it is coextensive with the precise concepts $\mathcal{P}_n$ and $\mathcal{P}_w$.\footnote{An obvious point of comparison for the method of squeezing is thus Beth's Definability Theorem.   For what a squeezing argument shows is that the constitutive, plus bridging, plus novel principles regarding the concept $\mathcal{C}$ do in fact provide an implicit definition for the predicate $C(x)$ in the manner required to apply the theorem over the mathematical theory $\mathsf{T}_P$.   But at least in the cases considered by Kreisel, the explicit definition is  given directly by the precise predicates $\pi_n(x)$ and $\pi_w(x)$ rather than extracted proof theoretically (e.g. by interpolation) in the manner typically employed in proving Beth's Theorem.}

\section{Case studies in informal rigour}
\label{casestudies}

In this section we will reconstruct and relate the three central arguments which Kreisel presented in \citeyearpar{Kreisel1967b} as examples of informal rigour as well as briefly evaluating them relative to subsequent developments.  While we have already stressed the centrality of Kreisel's validity and CH arguments both to his own thinking about informal rigour and its subsequent reception, his creating subject argument has had proportionally less effect on the discussion of his work within philosophy.     But this latter example is of considerable interest not only in regard to the influence within intuitionistic analysis, but also how it attests to the importance Kreisel invested in exploring the potential of new primitive notions in the operation of informal rigour.  We will discuss several other case studies -- inclusive of a fourth example which Kreisel added as an appendix to \citeyearpar{Kreisel1967b} -- in \S\ref{app} below.\footnote{In the presentation which follows we have reordered the sequence in which Kreisel himself presented his examples in light of various dependencies which will emerge below.   In a postscript to \citeyearpar{Kreisel1967b} Kreisel supplied for its reprinting as \citeyearpar{Kreisel1969c} he also suggested that his arguments should be updated in light of several subsequent publications.}    

\subsection{First-order validity}
\label{validity}

We have noted that Kreisel viewed what we have called his \textit{validity argument} as the most paradigmatic amongst his applications of informal rigour.   It is indeed straightforward to set out the argument in conformity with the schema S so that it may be understood as showing how the common concept of  \textit{first-order logical validity} ($\mathcal{Val}$) is squeezed between the precise concepts of \textit{first-order derivablity} ($\mathcal{D}$) and \textit{first-order model-theoretic validity} ($\mathcal{V}$). But of equal interest is the manner in which Kreisel motivates its individual steps, in particular in regard to the delineation of its formal and informal components and their relation to his about set theory which informed his CH (which we will return to discuss in \ref{ch}).  Despite some changes in emphasis, these features remain largely the same across the presentations discussed in \S\ref{context}.

\subsubsection{Initial schematization}
\label{valis}

In \citeyearpar{Kreisel1967b}, Kreisel introduced the predicates $\mathit{Val}, D$, and $V$ of sentences to denote the concepts $\mathcal{Val}, \mathcal{D}$ and $\mathcal{V}$.  $\mathit{Val}$ should thus be understood as a primitive predicate of the common language $\mathcal{L}_{C}$ whereas $D$ and $V$ should be understood as  defined mathematical predicates in the precise language $\mathcal{L}_P$.  We will return below to the choice of this language, the specific definitions of $D$ and $V$,  and also the background theory $\mathsf{T}_P$ over which the mathematical portion of the validity argument can be carried out.  But for expository purposes, it may safely be assumed that first-order derivability and validity are defined as in the familiar manner of modern logic textbooks -- i.e. that derivability is fixed with respect to some fixed proof system for first-order logic (say natural deduction) and that validity is defined as `truth in all models' in the manner which was ultimately made precise by \citet{Tarski1956c}.   For concreteness, we henceforth will officially adopt the definitions stated in van Dalen's \textsl{Logic and Structure} \citeyearpar{Dalen2008} -- i.e.
\begin{definition}
\label{DVdefns}
\begin{enumerate}[i)]
\item $D(\varphi^1)$ iff $\vdash_1 \varphi^1$ is derivable from no premises in the natural deduction system for first-order logic given by \textnormal{\citet[\S 3.8-10]{Dalen2008}}.
\item $V(\varphi^1)$ iff $\models_1 \varphi^1$ iff for all models $\mathfrak{M} \models_1 \varphi^1$ where $\mathfrak{M} = \langle A,R_1, R_2, \ldots \rangle$ with the domain $A$ a set, $R_1, R_2, \ldots$ are relations (functions, constants) on $A$ of appropriate arities interpreting the non-logical symbols of $\varphi^1,$ and $\models_1$ denotes Tarski's inductive definition of truth in a first-order model as given by \textnormal{\citet[\S 3.4]{Dalen2008}}.
\end{enumerate} 
\end{definition}
\noindent 

We have here also adopted Kreisel's convention of annotating symbols with their order so that, e.g., $\varphi^1$ signifies that $\phi$ is a sentence of a first-order language, $\vdash_1$ denotes first-order derivability, and $\models_1$ denotes first-order logical consequence.   Kreisel is less specific in defining these symbols -- speaking (e.g.) of derivability from `Frege's rules for first-order logic' (without citation).  But such conventional definitions coincide well with the manner in which he employs the predicates $D$ and $V$. The status which Kreisel assigned the common predicate $\mathit{Val}$ is more complicated.   But we will pass over this momentarily so as to first lay out the structure of the argument itself.

The argument proceeds from the premises
\begin{example}
\label{valarg}
\begin{enumerate}[i)]
\item $\forall \varphi^1(D(\varphi^1) \rightarrow \mathit{Val}(\varphi^1))$
\item $\forall \varphi^1(\mathit{Val}(\varphi^1) \rightarrow V(\varphi^1))$
\item $\forall \varphi^1(V(\varphi^1) \rightarrow D(\varphi^1))$
\end{enumerate}
\end{example}
to the conclusion 
\begin{example}[]
\begin{enumerate}[i)]
\setcounter{enumi}{3}
\item $\forall \varphi^1(V(\varphi^1) \leftrightarrow \mathit{Val}(\varphi^1))$
\end{enumerate}
\end{example}
This latter statement is, of course, the desired characterization of the common or `intuitive' notion of validity in terms of the precisely defined notion of first-order model theoretic validity.  And since this statement follows from (\ref{valarg}i-iii), it is straightforward to see that the validity argument conforms to the format of a squeezing argument as given by the schema S.  

As we have suggested, squeezing arguments of this form may be regarded as instances of the format for informally rigorous arguments we have proposed in \S\ref{3irs}.   In the case of argument (\ref{valarg}), however, the significance of several of the questions we have just glossed over comes into sharper focus when we seek to align it with the more general schema IR.  Amongst these are the following:
\begin{enumerate}[i)]
\item How should we understand the relationship between the common concept $\mathcal{Val}$ of first-order validity and the precise concept $\mathcal{V}$ of first-order model-theoretic  validity?
\item At what stage in the schema IR are premises (\ref{valarg}i-ii) justified -- i.e. should they be taken to be constitutive principles of $\mathcal{Val}$ accepted before the argument begins or must they be argued for at stages Ib or Ic?
\item In what formal language $\mathcal{L}_P$ should we understand Definitions \ref{DVdefns}i, ii to be stated  and in what theory $\mathsf{T}_P$ is the derivation of the mathematical premise (\ref{valarg}iii) carried out?  
\item What is the significance of Kreisel's restriction to \textit{first-order} logic in the formulation of the argument?
\end{enumerate} 

\subsubsection{On `the meaning of $\mathit{Val}$'}
\label{valmeaning}

As Kreisel's thinking about iv) is related to the goals he wished to achieve via his CH argument, we will postpone a proper discussion of this question until the end of \S\ref{ch}.  However question iv) is also related to question i), which is in turn difficult to discuss without taking into account some of the textual complexities we have mentioned above.   For although latter-day commentators of the validity argument have focused on a single characterization of $\mathit{Val}$ given in \citeyearpar{Kreisel1967b}, Kreisel himself discusses several possible interpretations in the sources surveyed in \S\ref{context}.  But of course the significance which can be assigned to the argument turns largely on how plausible it is to regard Kreisel's characterization of $\mathit{Val}$ as describing a genuinely `common' or `intuitive'  concept.   

There is in fact a high degree of conformity across Kreisel's presentations of the validity argument on four basic points.\footnote{Kreisel took these points mostly for granted in the central sources considered here. But in evaluating them it should be kept in mind that much of his work in the 1950s-1970s related to completeness phenomena with respect to different validity notions (classical, intuitionistic, finitist, and predicative) some of which we will discuss in \S\ref{standpoints}. In this category, the papers \citeyearpar{Kreisel1950}, \citeyearpar{Kreisel1952c}, \citeyearpar{Kreisel1954}, \citeyearpar{Kreisel1955a}, \citeyearpar{Kreisel1958b}, \citeyearpar{Kreisel1958k}, \citeyearpar{Kreisel1975} provide technical background while \citeyearpar{Kreisel1976},  \citeyearpar{Kreisel1980}, and (especially) \citeyearpar{Kreisel1956} contain relevant historical discussions.    Some traces of Kreisel's historical orientation also remain in the main sources described in \S\ref{context} as follows:  \citeyearpar[pp. 149-151, 178]{Kreisel1965a},  \citeyearpar[pp. 203, 234-235, 253, 255]{Kreisel1967a}, \citeyearpar[pp. 140-141, 153, 154]{Kreisel1967b},  \citeyearpar[pp. 167-168, 186, 196-197]{Kreisel1967c}. \label{kb}}  First, Kreisel was interested in the concept of logical validity (and cognates like logical consequence)  which he understood to figure in mathematical practice from at least the mid-nineteenth century onwards.   Second, while Kreisel took it as given that there was agreement amongst practitioners about instances of validity and consequence \textsl{within mathematics}, he says little which would indicate that he was concerned with the analysis of validity in everyday discourse.  Third, the broad conception which he took to stand behind the notion of validity employed in mathematical practice validates classical first-order logic.   Fourth, it is such an understanding of validity which Kreisel took to inform the work in metamathematics initiated by Frege and Hilbert which led to the formulation and proof of G\"odel's Completeness Theorem for first-order logic.

It still notable that Kreisel refrained from giving a truly general explanation of what he refers to as `intuitive validity'.  But what he does say is also compatible with the following familiar glosses:
\begin{center} $\phi$ is logically valid iff \\ it is impossible for $\phi$ to be false iff \\ there is no state of affairs in which  $\phi$ is false
\end{center}
In conformity with his interest in mathematical statements, Kreisel preferred the term `structure' over expressions like `state of affairs', `situation', or `possible world' which have subsequently been favored by philosophers while also leaving open exactly how such entities should be understood.    But suppose that we now assume -- as apparently did Kreisel -- that it is unproblematic to distinguish the logical constants appearing in $\varphi$ (of any order) and also that we let $\varphi^{\mathfrak{S}}$ denote the result of interpreting the non-logical expressions appearing in a sentence $\varphi$  as an appropriate sort of semantic entity relative to the structure $\mathfrak{S}$.   We then arrive at the following \textsl{interpretational} proto-analysis of validity:\footnote{Given his focus on G\"odel's Completeness Theorem, it in seems likely that Kreisel's specific waypoint for such a characterization was the notion of \textsl{universal validity} (\textsl{Allgemeing\"ultig}) which figured in the early development of metamathematics within the Hilbert school -- i.e. `a formula of the predicate calculus is called \textit{logically true} or, as we also say, \textit{universally valid} only if, independently of the choice of the domain of individuals, the formula always becomes a true sentence for any substitution of definite sentences, of names of individuals belonging to the domain of individuals, and of predicates defined over the domain of individuals, for the sentential variables, the free individual variables, and the predicate variables respectively' \citep[p. 55-56]{Hilbert1938} (see also  \citealp[p. 223]{Hilbert1929} and \citealp[p. 8]{Hilbert1934}).\label{uvnote}}
\begin{example} $\varphi$ is logically valid just in case $\varphi^{\mathfrak{S}}$ is true for all structures $\mathfrak{S}$. \label{interp}
\end{example}

This is in fact quite close to what Kreisel says under the heading `\textit{Meaning of $\mathit{Val}$}':
\begin{quote}
\footnotesize{The intuitive meaning of $\mathit{Val}$ differs from that of $V$ in one particular: $V(\alpha)$ (merely) asserts that $\alpha$ is true in all structures in the cumulative hierarchy, i.e., in all sets in the precise sense of set above, while $\mathit{Val}(\alpha)$ asserts that $\alpha$ is true in all structures $\ldots$   A current view is that the notion of arbitrary structure and hence of intuitive logical validity is so vague that it is absurd to ask for a proof relating it to a precise notion such as $V$ or $D$, and that the most one can do is to give a kind of plausibility argument. \hfill \citeyearpar[p. 153]{Kreisel1967b}}
\end{quote}

While everything said here is in apparent conformity with (\ref{interp}), two additional points stand out.   First, it is clear that Kreisel not only adopted the conventional model-theoretic definition of $V$ we have recorded as Definition \ref{DVdefns}ii, but he  assumed that models can be understood as \textit{sets} which appear as structures within the cumulative hierarchy.\footnote{This passage also points towards Kreisel's engagement with what might be called `the entanglement of set theory and semantics' -- i.e. the apparent inter-independence of the conventional semantic definitions of (e.g.) well-formedness, truth, and satisfaction in a model with a theory of sets which may itself require an axiomatic formulation in a formal language.   But although Kreisel in \citeyearpar{Kreisel1967b} seems to suggest that the latter must come before the former, some of his other treatments (e.g. \citeyear[\S 1.8]{Kreisel1965a}) suggest that this order should either be inverted or that set theory and semantics must be accorded even footing.} Second, the passage suggests that the only feature which he took to distinguish $V$ from $\mathit{Val}$ is that while the former requires truth with respect to all structures which are sets, the latter requires truth in \textit{all structures whatsoever}.

Taken on its own, this passage does not tell us how Kreisel understood the difference between models and structures.   But the immediate context in which it is embedded in \citeyearpar{Kreisel1967b} provides one initially plausible interpretation -- i.e. while $V(\varphi^1)$ requires only the truth of $\varphi^1$ in all \textsl{models} -- i.e. a tuple consisting of a \textsl{set-sized} domain $A$ together with relations, functions and constants on $A^k$ -- $\mathit{Val}(\varphi^1)$ requires the truth of $\varphi^1$ in structures whose domain is potentially a \textit{proper class} (but are otherwise like models in how they assign denotations to non-logical terms). 

Such an understanding of the distinction between $\mathit{Val}$ and $V$ is indeed consistent with much of what Kreisel says in \citeyearpar[\S 2]{Kreisel1967b}.   It has also been adopted by the majority of subsequent commentators on the validity argument.\footnote{E.g. \citep{Shapiro1987}, \citep{Etchemendy1990,Etchemendy2008}, \citep{Hanson1997}, \citep{Field1989,Field2008}, \citep{Smith2011a}, \citep{Halbach2020a,Halbach2020b}.}  There are, however,  several reasons to be dissatisfied with it.  Primary amongst these is that it seems out of keeping with Kreisel's historical orientation to suggest that the set/class distinction could play a role in characterizing the `intuitive' notion of validity.  For whereas the distinction originated with the set theoretic paradoxes around 1900 -- and was then codified later within axiomatic theories such as G\"odel-Bernays or Morse-Kelley set theory starting in the 1920s -- Kreisel is clearly interested in a `traditional' notion of validity which he assumed to have been employs mathematics well before this time.  Adopting this distinction as our preferred means of understanding $\mathit{Val}$ would thus prejudicially constrain how the validity argument may be viewed by suggesting that the common concept which it seeks to make precise is one which can only be appreciated via a technical distinction originating within set theory itself.\footnote{In fact Kreisel explicitly suggests in the first section of \citeyearpar{Kreisel1967b} that the notion of class which figured in our practices prior to the discovery of the paradoxes was ambiguous between `sets \textit{of} something' (i.e. defined relative to a bounding set) and `properties or intensions where one has no \textit{a priori} bound on the extension' (p. 143).  He then claims that Zermelo's reflection on the cumulative hierarchy has led to axioms which do not -- as a matter of practical fact -- require us to countenance a `multifurcation' in the concept of class or set.  (We will discuss the more extensive argument which  \citet[pp. 99-101]{Kreisel1965a} provides for this point in \S \ref{set} below.) But what Kreisel  goes on to say about `structures' is compatible with regarding them as being comprised of entities from any of these types.   This includes the formal theory of `explicitly definable properties' described in Appendix A of \citeyearpar{Kreisel1967b} to formalize the proof of the reflection principle described  in \S\ref{fieldnote2} below.  In fact Kreisel closes this section by remarking `Though the class formalism was originally introduced to deal with purely formal questions of finite axiomatisability, even for technical purposes it is good to ask oneself what classes \textsl{are} -- informal rigour!' (p. 165).}

Not only does Kreisel appear sensitive to this concern, but much of what he says (and does not say) about the meaning of $\mathit{Val}$ seems designed to exclude such a simplistic interpretation.\footnote{In fact  after the foregoing passage he continues `Let us go back to the fact (which is not in doubt) that one reasons in mathematical practice, using the notion of consequence or of logical consequence, freely and surely, (and, recall $\ldots$ the  `crises' in the past in classical mathematics $\ldots$ were not due to lack of precision in the notion of consequence).'  \citeyearpar[p. 153]{Kreisel1967b}.}  For instance, in \citeyearpar{Kreisel1965a} and \citeyearpar{Kreisel1967c} $\mathit{Val}$ is introduced as a primitive predicate glossed simply as `logical validity' or `intuitive validity'.  These sources provide less auxiliary explanation than does \citeyearpar{Kreisel1967b}.   But Kreisel asserts in each that we have a sufficient grasp of this notion to characterize its properties on the basis of `the notion of validity implicit in mathematical practice' \citeyearpar[p. 190]{Kreisel1967c}.

These point speak to the significance which Kreisel assigned to the validity argument -- i.e. 
\begin{quote}
\footnotesize{\textit{Nobody will deny that one knows more about $\mathit{Val}$ after one has established its relations with $V$ and $D$; but that doesn't mean that $\mathit{Val}$ was vague before.} \hspace*{1ex} \hfill \citeyearpar[p. 154]{Kreisel1967b}}
\end{quote}
Of course even if we accept the latter point about the `meaning of $\mathit{Val}$', the question of relating it to $V$ and $D$ remains.  This in turn brings us to the question of the justification of premises (\ref{valarg}i) and (\ref{valarg}ii).   Although there is again a reasonable degree of conformity among Kreisel's presentations, a few words are in order about each.

Principle (\ref{valarg}i) expresses that the axioms and rules we have adopted for first-order logic are \textsl{sound} with respect to the intended interpretation of $\mathit{Val}$ -- i.e. that if a statement is derivable from these axioms and rules then it is `true in all structures'.   Kreisel's most extended comment on this is as follows:
\begin{quote}
\footnotesize{[I]t is generally agreed that at the time of Frege who formulated rules for first order logic, Bolzano's set-theoretic definition of consequence had been forgotten (and had to be rediscovered by Tarski); yet one recognised the validity of Frege's rules $D_F$. This means that implicitly
$$\forall i \forall \alpha(D_F(\alpha^i) \rightarrow \mathit{Val}(\alpha^i))$$
was accepted, and therefore certainly $\mathit{Val}$ was accepted as meaningful. \hfill \citeyearpar[p. 153]{Kreisel1967b}}
\end{quote}
Several oddities may  be noted here.\footnote{It may seem strange both that Kreisel choose to speak of \textsl{Frege's} rules (as opposed to, e.g., those of Hilbert and Ackermann or of Gentzen) and also that he describes them as an axiomatization of \textsl{first-order} logic (as the system described in \citep{Frege1879} was higher-order).  One explanation of the first point is simply that by choosing as early an axiomatization as possible, he is able to stress that the `traditional' notion of validity was already in place when Frege formulated his system (which was 50 years before the proof systems with which we are now more familiar).  Note also that although Kreisel speaks of first-order logic, the statement $\forall i \forall \alpha(D_F(\alpha^i) \rightarrow \mathit{Val}(\alpha^i))$ in fact expresses the intuitive soundness of formally derivable statements of arbitrary order in his notation (wherein in the variable $i$ ranges over orders).}   But  this passage still makes clear that Kreisel regarded (\ref{valarg}i) as what we have called a \textsl{constituitive principle} about $\mathit{Val}$ -- i.e. one which he assumed would be accepted by his interlocutors without further argument.    This is further affirmed by the fact that Kreisel explicitly writes that (\ref{valarg}i) `is assumed to be recognized (as an \textsl{axiom})' in \citeyearpar[p. 116]{Kreisel1965a}.  

Understood from the perspective of contemporary debates, this may seem strange in two respects.   First, one might think that even if (\ref{valarg}i) is to be regarded as a constitutive principle of a genuinely common concept of validity, some further argument is required to show that the rules of \textsl{classical} logic -- as opposed, e.g., to those of intuitionistic or some other non-classical logic -- are genuinely \textsl{common} to our practices.   Second, one might also think that even if it is granted that we accept that the rules of classical first-order logic individually preserve intuitive validity, it might still be objected that an additional argument is required to show that all of their derivable consequences -- i.e. the \textsl{theorems} of first-order logic -- are intuitively valid.  

There is little evidence that Kreisel would have been concessive on the first point.  But this is perhaps understandable both in virtue of his focus on the status of validity in `traditional' mathematics and also because he acknowledged that  \textsl{finitist}, \textsl{predicative}, and \textsl{intutionisitc mathematics} may give rise to distinct concepts of validity amenable to informally rigorous analysis.\footnote{See \citeyearpar[p. 157]{Kreisel1967b} and  \S\ref{standpoints} below.  \citet[p. 246]{Read1994} can be understood as making a similar point in regard to the potential applicability of an argument with the structure of (\ref{valarg}) to non-classical notions of validity -- i.e. `[What Kreisel's] point shows is that Tarski-validity is extensionally safe for any provably complete deductive system which one believes is intuitively sound. That is, whatever your scruples, let $D(\alpha)$ represent `$\alpha$ is provable in my preferred first order logic' -- whether classical, intutionistic, relevant or whatever. Then, if you have a completeness proof for this logic relative to its Tarski semantics, its intuitive soundness (for you) will carry over to its Tarski soundness, that is, Tarski-validity, suitably defined $\ldots$'}    On the other hand, Kreisel illustrates his sensitivity to the second point in \citeyearpar[pp. 253-254]{Kreisel1967a} wherein he suggests taking (\ref{valarg}i) not as an axiom but rather as a consequence of another principle stating that the extension of $\mathit{Val}$ is closed under the relation of `immediate consequence' from premises. This suggests that even if it is not granted that the principle (\ref{valarg}i) is among the constitutive principles for $\mathit{Val}$, then it is still possible to present an additional informally rigorous argument for it on the basis of more basic constitutive principles needed to formulate the familiar argument `if the axioms are valid and the rules preserve validity, then the theorems are valid'.\footnote{\textsl{Pace} \citet{Field2008,Field2015a} as we will consider in \S\ref{valrec}.}   Thus regardless of whether $\mathit{Val}$ is regarded as a primitive or a defined concept, it will follow that (\ref{valarg}i) will be included in the set $\Gamma_2$ as defined in the schema IR.\footnote{Despite Kreisel's (apparently official) policy of treating $\mathit{Val}$ as a primitive concept, the foregoing considerations suggest that in the background he regarded the concept of validity \textsl{intensionally} -- i.e. as defined parametrically both in terms of a given set of \textsl{logical rules}, but also a given notion of \textsl{class} (and hence \textsl{structure}) relative to which their individual soundness is evaluated.   On the other hand, the validity argument itself can be taken to show that the principles (\ref{valarg}i,ii) are all that must be assumed of this concept in order to show it coincides \textsl{extensionally} with $V$.  (In light of the results discussed in note \ref{delta02note} this in turn only requires we accept the existence of a limited range of arithmetically definable classes which -- if one wished -- could themselves be understood either intensionally or extensionally).  As we will see in \S\ref{gmprec} and \S\ref{intval} this illustrates a distinction in how Kreisel came to regard his validity and creating subject arguments. \label{intnote}}

Turning now to principle (\ref{valarg}ii), Kreisel's explanations are even briefer.   For instance in \citeyearpar{Kreisel1965a} he simply observes that this principle `is \textsl{evident} since logical validity implies set theoretic validity' (p. 116) whereas in  \citeyearpar{Kreisel1967b} he writes that `one $\ldots$ accepts [\ref{valarg}ii] the moment one takes it for granted that logic \textsl{applies} to mathematical structures' (p. 154).  Both remarks reflect the observation that it will be recognized that (\ref{valarg}ii) follows from the meaning of $\mathit{Val}$ -- understood as truth in \textsl{all} structures -- as soon as it is realized that models are themselves structures of a certain sort.  It is again perhaps unclear whether this will be immediately accepted by all parties as a constitutive principle for $\mathit{Val}$.  But even if it is not, it is evident that Kreisel saw this as a example of how `pushing a bit farther than before the analysis of the intuitive notions' can lead to useful insights of the relationship between precise and common concepts.   In any case, it seems  Kreisel regarded it as uncontentious that (\ref{valarg}ii) will again end up in $\Gamma_2$.

Having addressed the status of the first two premises of the validity argument, we turn finally to (\ref{valarg}iii).  Once Definitions \ref{DVdefns}i,ii are accepted this is, of course, a statement of G\"odel's Completeness Theorem for first-order logic a special case of which can be formulated using the symbolism introduced above in the familiar way as follows:
\begin{theorem}
\label{comp}
For all $\varphi^1$, if $\models_1 \varphi^1$ then $\vdash_1 \varphi^1$.
\end{theorem}
\noindent In other words, if a first-order sentence is true in all models -- i.e. $V(\phi)$ -- then it is derivable in the first-order proof system we have adopted -- i.e. $D(\phi)$.   This closes the circle of inclusions in the validity argument by which $\mathit{Val}$ is `squeezed' between the narrow concept $D$ of first-order derivability and the wide concept $V$ of truth in all first-order models.   

\subsubsection{From schematization to formalization}
\label{valsf}

Kreisel often called attention to the significance of G\"odel's result not only in light of its role in his validity argument fact but for several other reasons related to informal rigour.   The first of these pertains to the historical setting in which G\"odel obtained the result in his dissertation \citeyearpar{Godel1929}.    As Kreisel notes, at this time the definition of model theoretic validity had not yet been been given formally in the manner of Definition \ref{DVdefns}ii but was rather understood in a semi-formal manner similar to that of `universal validity' as described by \citet{Hilbert1928}.\footnote{In fact, G\"odel's source in his dissertation was the first edition of \textsl{Grundz\"uge der Theoretischen Logik} \citeyearpar{Hilbert1928} wherein universal validity is described even less explicitly than in the second edition (see note \ref{uvnote}) as (essentially) `correctness for each substitution of predicates'  (pp. 72-73).   On the other hand, it appears that Definition \ref{DVdefns}ii may not to have been stated precisely in print until as late as \citep{Tarski1956c}.} The question thus arises how G\"odel can be credited as having \textsl{proven} a result which is now taken to have  the significance of (\ref{comp}) before the relevant definitions were in place.

According to Kreisel \citeyearpar[pp. 254-255]{Kreisel1967a}, \citeyearpar[p. 257]{Kreisel1967b}, G\"odel can be understood to have circumvented this problem in two stages.   First, relying on the intuitive concept of validity, he implicitly used the `obvious' observation 
\begin{equation}
\label{valN}
\forall \varphi^1(\mathit{Val}(\varphi^1) \rightarrow V^{\mathbb{N}}(\varphi^1))
\end{equation} 
This expresses in a variant of Kreisel's notation that if $\varphi^1$ is intuitively valid then $\varphi^1$ is true in all models having the natural numbers $\mathbb{N}$ as domain and relations on $\mathbb{N}^k$.   Second, he proved mathematically that 
\begin{equation}
\label{nummod}
\forall \varphi^1(V^{\mathbb{N}}(\varphi^1) \rightarrow D(\varphi^1))
\end{equation} 
-- i.e. that validity with respect to all such countable models is sufficient to imply first-order derivability.   Taken together with (\ref{valarg}i), this shows that $\mathit{Val}$ is squeezed not only between $D$ and $V$ but also between $D$ and $V^{\mathbb{N}}$.   If a common understanding of $\mathit{Val}$ can indeed be taken for granted, this supports Kreisel's contention that G\"odel can be regarded as demonstrating a result  `intuitively equivalent' to (\ref{comp}) by showing that first-order derivability is already guaranteed by truth with respect to a much narrower class of models than is required by the definition of $V$ given by \ref{DVdefns}ii.

Although this strengthening of G\"odel's Completeness Theorem is itself well-known, a sequence of  results which would arise from its metamathematical analysis also bears directly on the third questions left open above -- i.e. in what mathematical theory $\mathsf{T}_P$ should we understand (\ref{valarg}iii) to be demonstrated.   A first step was Bernays's observation that  (\ref{nummod}) can be strengthened to show that truth with respect to all \textsl{arithmetical models} -- i.e. those which have not only domain $\mathbb{N}$ but are also such that each of their non-logical symbols is interpreted by a formula of first-order arithmetic defining a predicate of appropriate arity -- is sufficient to guarantee first-order derivability.  In the early 1950s, it was shown that Bernays's result can be strengthened to show that completeness obtains with respect to truth in all $\Delta^0_2$-arithmetical models.\footnote{Bernays's original result has come to be known as the \textsl{Arithmetized Completeness Theorem} and is demonstrated in the second volume of \textsl{Grundlagen der Mathematik} \citeyearpar[pp. 234 - 253]{Hilbert1939}.  The strengthening of this result just mentioned amounts to the following: i) there exist consistent first-order formulas which are true in all $\Sigma^0_1 \cup \Pi^0_1$-models (and perforce in all $\Delta^0_1$ -- or \textsl{computable} -- models) but are not derivable in first-order logic; ii) on the other hand, truth with respect to all models wherein all non-logical symbols are defined by $\Delta^0_2$-arithmetical formulas \textsl{is} sufficient to ensure first-order derivability (see, e.g., \citealp{Kreisel1950}, \citealp[\S 72]{Kleene1952}, \citealp[\S 13.2]{Kaye1991}, \citealp{Dean2020e}).    As we will see in \S \ref{intval}, results about the definability of the classes which must be assumed to exist in order for the Completeness Theorem to hold also informed Kreisel's work on intuitionistic validity. \label{delta02note}}

These results initiated a metamathematical analysis of G\"odel's Completeness Theorem which was originally carried out in computability theory but later subsumed into the subject now known as Reverse Mathematics.   The latter is often described as directed at addressing the following question:  `\textsl{Which set existence axioms are needed to prove the theorems of ordinary, non-set-theoretic mathematics?}' \citep[p. 2]{Simpson2009}   The results just summarized can all be understood as pointing to the fact that the Completeness Theorem is indeed a piece of `ordinary' mathematics in the sense that its proof requires only minimal set theoretic assumptions.\footnote{This is in contrast to Kreisel's characterization of $\mathit{Val}$ (which quantifies over all structures and thus all sets).   He goes on to say that this contributes to the `remarkable conclusiveness' of the Completeness Theorem in regard to coincidence of $\mathit{Val}$ and $D$ which `one might have been hard put [to see] if the proof of [\ref{valarg}iii] had involved the assumption that there are non-denumarble measurable sets!' \citeyearpar[p. 254]{Kreisel1967a}.}

This observation is also embodied in another well-known result of \citet{Friedman1975a} -- i.e.  when the definitions (\ref{DVdefns}i,ii) are formalized in the language of second-order arithmetic, the statement (\ref{comp}) is provable in the axiomatic system known as $\mathsf{WKL}_0$.   This theory is formulated in the conventional language $\mathcal{L}^2_a$ of second-order arithmetic and consists of $\mathrm{I}\Sigma_1$ (i.e. the fragment of first-order Peano arithmetic with induction restricted to $\Sigma^0_1$-formulas), the second-order comprehension schema restricted to computable (i.e. $\Delta^0_1$-definable) sets, and the $\mathcal{L}^2_a$-formalization of the principle known as \textsl{Weak K\"onig's Lemma} (i.e. every infinite binary tree has an infinite path).  But  it can also be shown that an appropriately formalized version of (\ref{comp}) is in fact \textsl{equivalent} to Weak K\"onig's Lemma over $\mathsf{RCA}_0$ (i.e. the sub-theory consisting of only $\mathrm{I}\Sigma_1$ and comprehension for computable sets).  This imposes both an upper and a lower bound on both the axioms of the theory $\mathsf{T}_P$ over which we are assuming that the mathematical portion of the validity argument is conducted, as well as the expressivity of the language which is required to formulate it.\footnote{A complete statement and proof of these results can be found in \citep[IV.3]{Simpson2009}.\label{rmnote}}

Such a presentation might seem to add little to the version of the validity argument given at the beginning of this section.  But in the context of the sources we have been considering, it may also be taken to illustrate what Kreisel appears to have meant in \citeyearpar{Kreisel1967a} when he likened informal rigour to `philosophical proof'.  For as we have seen, Kreisel regarded the premises  (\ref{valarg}i,ii) as something akin to \textsl{axioms} involving a genuinely common notion of validity.  While he allowed that these premises could be justified by `pushing further intuitive notions', it seems he would have regarded as futile any attempt to provide further mathematical justification for them as such justification would inevitably rely on the very notion of intuitive consequence at issue.   On the other hand, we can now see that the result on which the argument relies -- i.e. G\"odel's Completeness Theorem -- can be regarded as a mathematical statement in its own right derivable from precisely delimited assumptions.

Putting these observations together, we may finally note that the premises of the validity argument can be formulated in the joint language $\mathcal{L}_J$ consisting of $\mathcal{L}_P = \mathcal{L}^2_a$ together with the single primitive predicate letter $\mathit{Val}'$ (where the latter is now understood as a predicate of G\"odel number of sentences in virtue of the need to arithmetize syntax when working over $\mathcal{L}^2_a$).   For suppose we let $\pi_1(x)$ and $\pi_2(x)$ be the formulas given in \citep[\S II.8]{Simpson2009} which respectively formalize first-order derivability and model-theoretic validity in $\mathcal{L}^2_a$.  Then (\ref{valarg}i-iii) respectively go over to
\begin{example}
\label{formvalarg}
\begin{enumerate}[i)]
\item $\forall x(\pi_1(x) \rightarrow \mathit{Val}'(x))$
\item $\forall x(\mathit{Val}'(x) \rightarrow \pi_2(x))$
\item $\forall x(\pi_1(x) \leftrightarrow \pi_2(x))$
\end{enumerate}
\end{example}

We can now see that it will be the $\mathcal{L}_J$ principles (\ref{formvalarg}i,ii) which should ultimately be included in $\Gamma_2$ at step Ib) in the scheme (IR).  Suppose we additionally take $\mathsf{T}_P = \mathsf{WKL}_0$ and $\mathsf{T}_K =  \mathsf{WKL}_0 \cup \Gamma_2$ at the end of step Id.   It will then follow that  (\ref{formvalarg}iii) -- i.e. the formalized version of the Completeness Theorem -- is derivable in $\mathsf{T}_P$.  And from this it follows that 
\begin{equation}
\label{concl}
\forall x(\mathit{Val}'(x) \leftrightarrow \pi_2(x))
\end{equation} 
will be derivable in $\mathsf{T}_K$ as desired.   The foregoing steps make explicit the derivation of the theorem embodying the validity argument to which Kreisel \citeyearpar[p. 190]{Kreisel1967c} alludes.\footnote{In fact in this version Kreisel labeled (\ref{valarg}iv) as a \textsl{theorem} and (\ref{valarg}) as its \textsl{proof}.}   As we will see below, it also possible to conform the other examples of informally rigorous arguments which Kreisel discusses in \citeyearpar{Kreisel1967b} to this model of `philosophical proof'.   But we also propose that it is this standard of rigor -- inclusive of a clear delineation of common and precise language, identification and defense of constitutive principles, and the demarcation of a precise theory capable of deriving the mathematical theorems employed in the argument -- against which these other arguments should be evaluated.  

Before leaving the validity argument, it will finally be useful to contrast the theory $\mathsf{T}_K$ with the sorts of theories introduced by \citet{Montague1963} to formalize reasoning about a notion of logical necessity \textsl{prima facie} similar to Kreisel's $\mathit{Val}$.   Although many of the theories Montague considered were inconsistent, this is typically so because they include a `reflection' axiom which in the present notation would take the form $\mathit{Val}'(\ulcorner \varphi \urcorner) \rightarrow \varphi$.  But as we will discuss further in \S\ref{valrec}, not only does Kreisel fail to adopt such a schema as an axiomatic principle for $\mathit{Val}$, he would have little motivation to do so.   And in the case at hand, it is additionally easy to see that $\mathsf{T}_K$ is \textsl{conservative} over $\mathsf{T}_P = \mathsf{WKL}_0$ and thus is consistent as long as the latter theory is.\footnote{For note that the validity argument itself shows that $\mathsf{T}_K$ possesses a model in which the extension of $\mathit{Val}'$ can be taken to coincide with the predicate $\pi_1(x)$ formalizing first-order derivability.  As the extension of $\pi_1(x)$ (and thus of $\mathit{Val}$) is not $\Delta^0_1$-definable in virtue of Church's Theorem, this construction cannot be carried out in $\mathsf{WKL}_0$ itself.   But as long as $\Delta^0_1$-comprehension is not extended to formulas containing $\mathit{Val}'$ -- which is not required by the argument -- then the provable coincidence of $\pi_1(x)$ and $V$ in $\mathsf{T}_K$ can be used to replace occurrences of the latter in $\mathsf{T}_K$- proofs of $\mathcal{L}_P$-statements.  And this in turn provides an alternative proof-theoretic demonstration that $\mathsf{T}_K$ is a conservative extension of $\mathsf{T}_P$.}  As we will see in \S \ref{gmp}, this is a feature which distinguished the validity from Kreisel's creating subject argument.

\subsubsection{The reception of the validity argument}
\label{valrec}

Amongst his applications of informal rigour, Kreisel's validity argument has had the widest philosophical impact.  This is so largely in virtue of the influence it has exerted on the contemporary literature on logical consequence.   A systematic examination of this influence -- which has often been indirect for the reasons discussed in \S\ref{intro} and \S\ref{context} -- is beyond the scope of the current paper.     What we will do instead is to briefly remark on how the details of the reconstruction we have just presented bear more directly on the reception of the argument by two of its original popularizers as well as a more recent proposal for formalizing the argument in a manner similar to \S\ref{valsf}.

\paragraph{Field}

The term `squeezing argument' appears to have been introduced by Harty Field  to describe Kreisel's validity argument in a postscript to \citep{Field1984} published in \citep{Field1989}.\footnote{Prior to this \citet[pp. 53-56]{Quine1970a} -- who was both  a contributor to the volume in which \citeyearpar{Kreisel1967a} appeared and an attendee at the event where \citeyearpar{Kreisel1967b} was delivered -- made use of a sequence of observations involving the Arithmetized Completeness Theorem in a section entitled `Saving on sets' to argue for a substitutional treatment of validity similar to that which we have employed in the reconstruction given in \S\ref{valsf}.   But he neither names neither the argument nor attributes it Kreisel.}  In his initial engagement with the argument, Field approvingly cites Kreisel's treatment of $\mathit{Val}$ as a primitive notion suggesting that that the validity argument then provides `a much more satisfying understanding of the significance of the completeness theorem $\ldots$ than one gets by taking the notions of implication and consistency to be defined in the Tarskian manner' \citeyearpar[p. 32]{Field1989}.    In \citeyearpar[pp. 119-123]{Field1989} and \citeyearpar{Field1991} he then goes on to develop a modalized version of the argument in support of the claim that his prior use of metalogical results in his program for developing a nomalisitic reconstruction of mathematics can themselves be given a suitably nominalistic interpretation.\footnote{\citet[pp. 16-20, p. 111]{Field1980} had originally made use of the Completeness Theorem for first-order logic in the course of his argument that the adjunction of a mathematical theory to an appropriate nominalized physical one yields a conservative extension.  Given Kreisel's strong tendency towards mathematical realism (see \S\ref{set}), it seems likely that he would have been hostile to Field's background fictionalist project.  Although we will pass over this difference in orientation here, the following points of contact may still be noted: 1) Field's argument requires that his modal replacement COMP for the Completeness Theorem \citeyearpar[p. 13]{Field1991} be formally provable in some mathematical theory $\mathsf{S}$; 2) as Field notes since it is possible to choose $\mathsf{S}$ so that it is finitely axiomatizable (e.g. $\mathsf{WKL}_0$ or in his case G\"odel-Bernays set theory), this allows COMP to be formulated as a single conditional with $\mathsf{S}$ as antecedent; 3) although Field goes on to philosophically advocate for the \textsl{logical necessity} of this conditional, the fact that it is validated \textsl{non-vacuously} will still depend on whether the background semantics for second-order (or multi-sorted first-order) logic which are assumed are sufficient to guarantee the existence of structures satisfying $\mathsf{S}$.}

Field has more recently returned to engage with the validity argument in a different manner by calling attention to the fact that although the principle (\ref{valarg}i) -- i.e. $D(\phi^1) \rightarrow \mathit{Val}(\phi^1)$ -- plays the role of expressing the `intuitive soundness' of our adopted proof system for first-order logic, a curious feature of the argument is that it does not rely on a formal proof of soundness to deliver the conclusion $D(\phi^1) \rightarrow V(\phi^1)$.  In this regard  \citet[\S 2.3]{Field2008}, \citeyearpar[\S 2.6]{Field2015a} stresses two points germane to our prior discussion: i)  while the argument may provide an adequate justification for identifying intuitive validity with (formal) first-order model-theoretic validity, it does not provide a non-circular argument  that (formal) first-order derivability entails truth; ii) there is a related difficulty in regarding the conventional soundness proof for first-order logic as demonstrating that all formally derivable sentences are true.  

These points draw attention to an aspect of what Kreisel originally wished to achieve via the validity argument which may be at odds with the aims of Feild and other subsequent commentators.   For note that the notion of \textsl{truth} which is at issue here is presumably that of `truth simpliciter' -- i.e. truth in `the actual world'.   On the other hand, we have seen that for Kreisel, the truth of a mathematical statement is  understood as \textsl{truth relative to a structure}.   But although $\mathit{Val}(\varphi^1)$ is intended to express truth in \textsl{all} structures, if $\varphi^1$ is a mathematical sentence it will be the exception rather than the rule that there will be a clear sense of what it would mean for $\varphi^1$ to be true `in the actual world'.   This is canonically illustrated by several of the examples which Kreisel himself repeatedly employs -- e.g. if $\phi^1$ is a statement of group theory such as $\forall x \forall y(x \cdot y = y \cdot x)$ or a statement of geometry such as the Parallel Postulate.

This point notwithstanding, in \citeyearpar{Kreisel1967b} Kreisel considers in detail a particular case which can be understood as anticipating part of the argument which Field \citeyearpar[\S 2.4]{Field2008} ultimately gives against investing the conventional soundness proof for first-order logic with its  customary significance.   For consider a sentence $\varphi^1$ in the language of first-order set theory which contains the sole non-logical symbol $\in$.   In this case it might be thought that there is a unique `actual' structure in which mathematical sentences are intended to be interpreted -- i.e.  what is conventionally known as the \textsl{cumulative hierarchy of sets} $\mathbb{V}$ (see \S\ref{set}).  Kreisel considers in particular the statement $\mathit{Val}(\varphi) \rightarrow \varphi_{\in}$ wherein $\varphi_{\in}$ is intended to express that the set-theoretic sentence $\varphi$ is true when its quantifiers range over  $\mathbb{V}$ and its membership symbol is interpreted as the `actual' membership relation of this structure (i.e. $\in$).   But as he points out, it is by no means obvious that all instances of the $V(\varphi) \rightarrow \varphi_{\in}$ will hold.  For whereas $V(\varphi)$ requires only that $\varphi$ is true in all set-size models, $\varphi_{\in}$ requires that $\varphi$ is true in $\mathbb{V}$ (whose domain is not a set).  There is thus indeed a \textsl{prima facie} reason to be concerned about the status of the related proof-theoretic reflection principle $D(\varphi) \rightarrow \varphi_{\in}$ for first-order logic whose justification exemplifies Field's basic concern.\footnote{The same sequence of observations also serves as a motivation for Boolos's \citeyearpar{Boolos1985b} proposal that first-order validity should be defined in a schematic manner (which he refers to as `supervalidity') similar to Kreisel's definition  of second-order consequence in set theory given in \citeyearpar[pp. 256-257]{Kreisel1967a} (see also \S \ref{chis} below).  Boolos acknowledges (p. 343) that the link between this notion and first-order derivability is still mediated by the Completeness Theorem via argument (\ref{valarg}) despite his previous remark (p. 340) that it is `strange' that an appeal to a `non-trivial' result must be made in order to demonstrate that a sentence is true if it is valid. He also suggest that his account can be extended to second-order validity in a manner which avoids commitments to an ontology of classes, essentially by employing a truth definition similar to that described by \citet[pp. 162-165]{Kreisel1967b}.  But in this case the coincidence of validity with second-order derivability can no longer be guaranteed for the reasons we discuss in \S \ref{chsol} below.  See \citep{Halimi2017} for a recent reappraisal of this situation.}

But but only was Kreisel clearly aware of this issue, he also goes on to present a mathematical argument in \citeyearpar[pp. 155-156]{Kreisel1967b} to show that each instance of $V(\varphi) \rightarrow \varphi_{\in}$ \textsl{is} derivable in first-order Zermelo-Fraenkel set theory $[\mathsf{ZF}]$.   One of the central steps is as an application of the Montague-L\'evy model-theoretic reflection schema -- i.e. if for all $\phi \in \mathcal{L}^1_{\mathsf{Z}}$, if $\mathbb{V} \models \phi$ then there exists $\alpha$ such that $R(\alpha) \models \phi^1$.\footnote{See, e.g., \citep[\S 12]{Jech2003}.}  He then uses this to show that all instances of $D(\varphi) \rightarrow \varphi_{\in}$ are in fact derivable in $\mathsf{ZF}$.  But he also observes that if we considered instead a \textsl{finitely axiomatized} set theory $\mathsf{S}$ which conservatively extends $\mathsf{ZF}$ -- e.g. a theory like G\"odel-Bernays set theory $\mathsf{GB}$ which has sorts for both sets and classes -- then there will be $\mathcal{L}_{\mathsf{S}}$-sentences $\psi$ in the new two-sorted language for which $\mathsf{S} \not\vdash D(\psi) \rightarrow \psi$ (although of course all instances of this schema for $\mathcal{L}_{\mathsf{ZF}}$-sentences will still be derivable).\footnote{Kreisel played an important role in disseminating several of the  results on which this observation relies  -- e.g. the formulation of a conservative theory of `explicitly definable classes' $\mathsf{P}$ which he presents in Appendix A of \citeyearpar{Kreisel1967b}, definability of satisfaction for first-order $\mathcal{L}_{\in}$ in $\mathcal{L}_{\mathsf{P}}$, the fact that the cut elimination theorem for pure first-order logic can be proven in $\mathsf{ZF}$ (and in fact in far weaker systems).    Although references to the original sources for these results are mostly suppressed in \citeyearpar{Kreisel1967b}, a more complete account can be patched together from (e.g.) \citep{Kreisel1968a}, \citep[\S 4]{Hajek1998}, and \citep[II.8.8]{Simpson2009}.}

The sequence of mathematical results which Kreisel adduces in his proof of this result parallel the initial steps in the philosophical argument which \citet{Field2008} gives in a section entitled `The unprovability of soundness' (pp. pp. 48-55).   Field suggests that the result in question should inspire us to engage in a sort of iterated process of semantic ascent.\footnote{At the beginning of this process we start out by realizing that our schematically axiomatized theories (like $\mathsf{ZF}$) may be finitely axiomatized by moving to an appropriate multi-sorted conservative extensions (like $\mathsf{GB}$).  Field then suggests that we should come to regard such a finitely axiomatized theories as inadequate in virtue of failing to prove some instance of $D(\psi) \rightarrow \psi$ in its new multi-sorted language.   On this basis, he then proposes that we should (repeatedly) extend our schema to the new language.   At the first step, this (repeatedly) yields a non-finitely axiomatizable theory in which all instances of  $D(\psi) \rightarrow \psi$ are provable.   But at the second step, this (repeatedly) causes us to adopt a new multi-sorted conservative extension which fails to prove the full soundness of first-order logic in its new language (etc.).}   He then employs the conclusion of this section -- i.e. that the conventional interpretation of the soundness proof for first-order logic is `a bit of a hoax' \citeyearpar[p. 48]{Field2008} - as part of the motivation for the theory of truth which he develops later in the book.\footnote{It should be noted, however, that the theory \citet{Field2008} ultimately proposes is not based on a \textsl{typed} notion of truth of the sort which the ascent described procedure in the prior note might seem to suggest.} For his part, however, Kreisel appears to have drawn exactly the opposite moral:  `Looking at the intuitive relation $\mathit{Val}$ leads not only to formal proofs as in [the validity argument] but also incompleteness theorems' \citeyearpar[p. 156]{Kreisel1967b}. \label{fieldnote2}

\paragraph{Etchemendy}

The second contemporary theorist to engage with Kreisel's validity argument appears to have been John Etchemendy in his well-known monograph \textsl{The Concept of Logical Consequence} \citeyearpar{Etchemendy1990}.   Therein he offers an extended critique of Tarski's \citeyearpar{Tarski1936} account of logical validity and consequence.   We have seen that the validity argument is often presented as an elaboration of the Tarskian analyses.   But not only does Etchemendy again attempt to adapt the argument to his own purposes, there are also several points of technical contact between his approach and Kreisel's.  

Etchemendy's full proposal is often described as difficult to describe in generalities and controversial in its details.  But one of his central contentions is that the Tarskian analysis does not provide an adequate  account of what he (also) refers to as the `intuitive' notions of logical validity or consequence.  This is because Etchemendy takes the basic `reductive' principle on which they are based -- i.e. that of analyzing the validity of an argument in terms of the (`ordinary') truth of an appropriate universal generalization derived by quantifying over interpretations -- as failing to provide an `independent guarantee' that the truth of the premises will ensure the truth of the conclusion.\footnote{One of his characteristic formulations of this point is as follows: `Suppose we have an argument form all of whose instances preserve truth, just as the reductive account requires, but suppose that the only way to recognize this is, so to speak, serially -- by individually ascertaining the truth values of the premises and conclusions of its instances. In other words, suppose there is no independently recognizable guarantee of truth preservation, as there is with modus ponens, only the brute fact that the instances preserve truth.   Would an instance of this argument form be logically valid? Clearly not. For example, we could never come to know the conclusion of such an argument in virtue of our knowledge of its premises.'  \citeyearpar[p. 268]{Etchemendy2008} \label{etchnote}}  In order to illustrate this he introduces a sequence of examples which are intended to illustrate how different background assumptions can lead to an account on which the Tarskian analysis both \textsl{overgenerates} and \textsl{undergenerates},  respectively by assessing too many or two few sentences to be valid.   He thus also concludes that there is no reason to suspect that the Tarskian account should in general provide an extensionally adequate analysis of `intuitive validity'.

This general point notwithstanding,  Etchemendy also acknowledges in \citeyearpar[\S 11]{Etchemendy1990} that the validity argument can be repurposed to show that the Tarskian analysis does not overgenerate in the particular case of first-order logic.  Using Etchemendy's notation $\mathit{Ltr}(\phi)$ for `$\phi$ is an intuitive logical truth' (in his sense) his argument can be reconstructed as follows:  i) Kreisel's argument yields $\mathit{Val}(\phi) \rightarrow V(\phi)$ (relative to the interpretation of $\mathit{Val}$ discussed above); ii) $V(\phi) \rightarrow D(\phi)$ follows via the Completeness Theorem; iii) $D(\phi) \rightarrow \mathit{LTr}(\phi)$ can also be accepted as standard first-order proof systems are sound with respect to intuitive validity; thus iv) $\mathit{Val}(\phi) \rightarrow \mathit{LTr}(\phi)$ -- i.e. since all sentences which are true in all structures are also intuitively valid the Tarskian account does not overgenerate.\footnote{See also \citeyearpar[p. 275]{Etchemendy2008} for a similar account.}  

Etchemendy would later make clear that the notion of logical consequence in which he is ultimately interested subsumes everyday reasoning involving modal, epistemic, indexical, and diagrammatic notions which he takes to go far beyond first-order logic  (e.g. \citeyear[p. 282-295]{Etchemendy2008}).  These notions are also not typically expressed in the mathematical languages which we have seen were Kreisel's primary concern.   Nonetheless, Etchemendy also makes substantial use of examples similar to  Kreisel's to illustrate his other contention that the assurances supplied by the validity argument do not extend to show that the Tarskian account fails to undergenerate or that it can be extended to non-first-order languages.    This includes the observation that  if only finite domains are considered the Tarskian analysis would overgenerate by calling valid the negation of the statement that a transitive and irreflexive relation must have a minimal element  \citeyearpar[p. 118]{Etchemendy1990}, that it would undergenerate if `intuitive validity' subsumed the $\omega$-rule but nonstandard models of arithmetic are considered \citeyearpar[p. 148]{Etchemendy1990}, and that when extended to second-order logic it suggests that either the Continuum Hypothesis or its negation is an intuitive validity  \citeyearpar[p. 276-277]{Etchemendy2008}.

It is evident that Etchemendy's ultimate interests lie much closer to natural language than do Kreisel's.   But his occupation with mathematical examples still suggests that the critical part of his program is motivated by the same sort of foundational which inspired Kreisel's formulation of the validity argument.\footnote{This is most evident in his discussion of the position he calls `finitism' -- i.e. the view that there are only finitely many objects.  Although Etchemendy does not express sympathy for this position himself, he observes that in order for the Completeness Theorem to hold there must exist at least countably many objects (i.e. that finitism is false).   This is illustrative of the sort of ontological considerations he takes to bear on whether formal results form metalogic can be used to illuminate `intuitive validity'.   But in this case the results described in note \ref{delta02note} can in fact be understood as further delimiting the commitments of the Completeness Theorem by showing that the $\Delta^0_2$-definable arithmetical models already constitute what Etchemendy calls a  `rich class' -- i.e. if a sentence is true in all models in this class, then it is derivable in first-order logic.}  This in turn suggests that while their motivations are different, Kreisel may have been in a position to offer rejoinders to at least some aspects of Etchemendy's critiques of the Tarskian accounts of  validity and consequence.   For suppose that he were in a position to rule out  Etchemendy's examples of undergeneration -- e.g. by arguing that $\omega$-consequence does not fall under the relevant `intuitive' conception.  Then when combined with the direction of the validity argument which Etchemendy accepts, it would indeed entail that the Tarskian account was extensionally adequate in at least the case of first-order logic.    But if the argument is understood as an example of `philosophical proof' in the manner which Kreisel intended, then it might also be understood as taking some steps towards addressing Etchemendy's underlying concern about its \textsl{intensional adequacy} as well.\footnote{Note in particular that the Completeness Theorem itself can be understood as providing a means of bridging the epistemological gap which Etchemendy suggests (cf. note \ref{etchnote}) stands between our knowledge of the Tarskian validity of an argument (which is akin to a $\Pi^1_1$-statement quantifying over models) and the fact that its premises entail its conclusion in a recognizable manner (e.g. which is akin to a $\Sigma^0_1$-statement asserting the existence of a proof of its conclusion from its premises).  For while the Completeness Theorem may not itself be obvious, a mathematical demonstration that it holds might itself be taken to provide the sort of epistemic connection between the preservation of truth value and derivability which Etchemendy appears to call for.}

\paragraph{Halbach}

The formalization of the validity argument we have given in \S\ref{valsf} can also be compared to a recent treatment by Volker Halbach \citeyearpar{Halbach2020a,Halbach2020b}.  Halbach suggests that Kreisel's  policy of treating $\mathit{Val}$ as a primitive predicate can be improved upon by developing an axiomatic theory of substitutional quantification within which it is possible to formalize the interpretation `$\varphi^1$ is true in all structures'  by quantifying over an appropriately broad class of substitution instances for the non-logical symbols in $\varphi^1$.   One potential advantage of this approach is that it allows us to define $\mathit{Val}(\varphi)$ in the object language which appear to admit proper class-sized interpretations -- e.g. the so-called `homophonic' interpretation wherein the non-logical symbols in $\varphi$ (such as $\in$ -- as in our discussion of Field) are simply replaced by themselves.   

There are, however, several reasons to suspect that Kreisel would have been dissatisfied with such a reconstruction.   First, there is the evident risk that it conflates formal and informal rigour -- e.g. by replacing his informal arguments for (\ref{valarg}i,ii) with mathematical results (e.g. Lemma 4.3 and Theorem 4.4 in \citealp{Halbach2020b}) which must themselves be recognized as consequences of correct mathematical reasoning relative to a prior understanding of validity.   Second, as we have discussed in \S\ref{valmeaning}, it is unclear that by simply providing a linguistic surrogate for class-sized interpretations, the envisioned approach does justice to the understanding of $\mathit{Val}$ which Kreisel appears to have had in mind.   And third, it seems likely that Kreisel would have dissented from Halbach's remark that `One would expect from an adequate conceptual analysis of logical validity that it is obviously adequate and that establishing the adequacy of the analysis does not require an ingenious proof' \citeyearpar[p. 318]{Halbach2020b}.   For on the one hand, Kreisel repeatedly stresses that the coincidence of $\mathit{Val}$ and $V$ requires `philosophical proof' precisely because it is \textsl{non-obvious}.  And on the other, it is also not clear that Halbach's formalization of substitutional validity -- which is carried out over an axiomatic theory of satisfaction adjoined non-conservatively to Zermelo-Fraenkel set theory -- is any less `ingenious' than the formalization of the Completeness Theorem in $\mathsf{WKL}_0$.

\subsection{Intuitionistic analysis and the creating subject}
\label{gmp}

In this section we will discuss Kreisel's engagement with the notion of the \textsl{creating subject} as it figured in Brouwer's development of intituitionistic analysis.   In \citeyearpar{Kreisel1967b} he formulated a version of an argument in which Brouwer had exploited this concept as his third example of informal rigour.\footnote{Kreisel referred to `the thinking subject' and also occasionally `the creative subject'.   Although there is also some variation between these terms in subsequent sources, we have followed \citet{Van-Atten2018} -- who provides a detailed reconstruction of Brouwer's original argument as well as  providing much additional information which is useful for understanding the context of Kreisel's presentation -- in standardizing on the original expression `creating subject'.}     Prior to this, Kreisel had worked extensively on intuitionistic mathematics and formal systems.   As we will discuss further in \S\ref{intval}, this included his investigation of the completeness of  intuitionistic first-order logic, his proposed formalization of the so-called \textsl{proof interpretation} of the intuitionistic connectives in the form of his \textsl{Theory of Constructions} as well the formulation of one of the earliest axiomatizations of intuitionistic analysis.   

But what is most germane to the current context is that Kreisel's immediate audience was likely to have been familiar with all of these developments.  This notably included Heyting and Myhill -- to whose responses Kreisel replied in detail in the comments which are  published at the end of \citeyearpar[pp. 178-186]{Kreisel1967b} -- and also Kleene -- who chaired the session in which Kreisel's paper was delivered and whose textbook on intuitionistic analysis \citep{Kleene1965} had recently been published.   Perhaps in light of this, Kreisel failed to provide a detailed description of the original construction of Brouwer which his own argument was intended to formalize.   We will thus begin by briefly reviewing the relevant background in intuitionistic analysis.

\subsubsection{Background} 
\label{gmpback}

In the course of his critique of classical reasoning in mathematics, Brouwer employed a device which has come to be known as a \textsl{weak counterexample}.  These take the form of implausibility arguments which illustrate how the adoption of certain logical or mathematical principles gives rise to constructively dubious consequences.   It was in the context of providing such a counterexample in which Brouwer first employed the notion of the creating subject, conceived as an idealized mathematical agent engaged with verifying (or refuting) a class of mathematical statements by generating constructive proofs (or refutations) whose existence (or non-existence) is then identified with the truth (or falsity) of the statements in question.  

The generation of such proofs is understood as a sequence of cognitive acts occurring at stages which  potentially be referenced can themselves in a mathematical construction.  This is canonically illustrated by Brouwer's original application of the creating subject in a mathematical argument:

{\footnotesize
\begin{quote}
Let $[A]$ be a mathematical assertion that \textsl{cannot be tested}, i.e. for which no method is known to prove either its absurdity or the absurdity of its absurdity.  Then the creating subject can, in connection with the assertion $A$, create an infinitely proceeding sequence of rational numbers $a_1 a_2 a_3,\dots$ according to the following direction: As long as, in the course of choosing the $a_n$, the creating subject has experienced neither the truth, nor the absurdity of $\alpha$, every $a_n$ is chosen equal to $0$ $\dots$ However, as soon as between the choice of $a_{r-1}$ and that of $a_r$, the creating subject has obtained a proof of the truth of $A$, $a_r$, as well as $a_{n+v}$, for every natural number $v$ is chosen equal to $2^{-r}$. And as soon as between the choice of $a_{s-1}$ and that of $a_s$ the creating subject has experienced the absurdity of $\alpha$, $a_s$ as well as $a_{s+v}$ for every natural number $v$ is chosen equal to $-2^{-s}$. \hfill \citep[p. 478]{Brouwer1948a}
\end{quote}}

This passage occurs within the context of a discussion intended to illustrate the difference between two notions of the inequality of real numbers as understood within intuitionistic analysis.  In this context real numbers are identified with convergent choice sequences of rationals $\langle r_n \rangle_n$.     It is possible to distinguish in this setting between the \textsl{non-identity} of real numbers $x,y$ (denoted $x \neq y$) and their \textsl{apartness} (denoted $x \# y$).   Brouwer characterizes the first of these relations as \textsl{negative}  (in the sense that it expresses `the absurdity of a constructive property') and the second as \textsl{positive} (in the sense that it is intended to express the existence of a difference between $x$ and $y$).  

The distinction between these notions has traditionally been analyzed by first adopting the following definitions of the equality and less than relations on real numbers $x$ and $y$ given by sequences $\langle r_n \rangle_n$ and $\langle s_n \rangle_n$:
\begin{align}
(x = y) &=_{\mathrm{df}} \A n \E m \A k \geq m (|r_k - s_k| < 2^{-n}) \label{eqdefns} \\ 
(x < y)  &=_{\mathrm{df}}  \exists k \exists n \forall m (|s_{n+m} - r_{n+m}| > 2^{-k})
\end{align}
The prior characterization of $x \neq y$ can be made precise as the constructive negation of (\ref{eqdefns}) -- i.e. $(x \neq y) =_{\mathrm{df}} \neg (x = y)$ -- while that of $x \# y$ can be made precise via the additional definition 
\begin{example}
$(x \, \# \, y) =_{\mathrm{df}} (x < y) \ \vee \ (y < x)$
\end{example}

It is not difficult to see from these definitions that the statements $x \neq y$ and $x \# y$ are classically equivalent.    But note that whereas $x \# y$ asserts the existence of a rational number separating the absolute difference of $x$ and $y$ from $0$, $x \neq y$ merely asserts that the assumption of the equality of their leads to a contradiction.  Thus relative to the intuitionistic understanding of the logical connectives, there appears to be no reason to expect that $x \neq y$ should constructively imply $x \# y$.  Brouwer's original creating subject argument can thus be understood as a weak counterexample illustrating that we should not expect to be able to constructively prove that
\begin{example}
$\A x \A y(x \neq y \rightarrow x \# y)$
\label{ineqapart}
\end{example}
Modifying Brouwer's presentation slightly, this argument be presented as follows:
\begin{example}
\begin{enumerate}[i)]
\item Let $A$ be a proposition that is not \textsl{testable} --  i.e. no method is known for proving $\neg A \lor \neg\neg A$.
\item The creating subject now constructs a choice sequence $\langle r_n \rangle_n$ as follows:
\begin{enumerate}[a)]
\item If by the point at which $r_n$ is chosen, the creating subject has proven neither $A$ nor $\neg A$, then $r_n =0$.
\item If between the choice of $r_{m-1}$ and $r_m$ the creating subject has proven $A$, then $r_n$ for all $n > m$ is chosen to be $2^{-m}$.
\item If between the choice of $r_{m-1}$ and $r_m$ the creating subject has proven $\neg A$, then $r_n$ for all $n > m$ is chosen to be $-2^{-m}$.
\end{enumerate}
\item It follows from this definition that the sequence $\langle r_n \rangle_n$ converges to a real number $r$ and also that  $r=0 \leftrightarrow (\neg A \land \neg\neg A)$.  As the righthand side of this biconditional is absurd, the creating subject can conclude that $r \neq 0$.  
\item But now suppose for reductio that we were also able to conclude that $r < 0 \ \vee 0 < r$.  In this case the creating subject could reason constructively by cases as follows:
\begin{enumerate}[a)]
\item Suppose $r < 0$.  It then follows by definition of $<$ that $\neg (r > 0)$ and thus also from the definition of $r$ that the creating subject never proves $A$.   But from this (and the constructive meaning of negation), it follows that $\neg A$ and thus also $\neg A \ \vee \neg \neg A$. But then $A$ has been tested, contrary to our assumption. 
\item Suppose $0 < r$.  In this case a similar argument yields  $\neg \neg A$ and thus again $\neg A \ \vee \neg \neg A$.  But then $A$ has again been tested, contrary to our assumption. 
\end{enumerate}
\end{enumerate}
\label{csarg}
\end{example}

The foregoing illustrates how reasoning about the creating subject may be invoked in an argument conducted within intuitionistic mathematics.   But a notable feature of the argument is its dependence on the existence of untestable propositions as exemplified by currently open problems such as the Goldbach Conjecture or Riemann Hypothesis.  As the existence of such propositions is a contingent feature of our current state of knowledge, Kreisel went so far as to characterize Brouwer's argument for the implausibility of (\ref{ineqapart}) as `empirical' \citeyearpar[p. 139]{Kreisel1967b}.   One of his aims in applying informal rigour to the practice of intuitionistic analysis was thus that of showing how such reasoning could be sharpened by providing a precise account of the properties which are assumed to hold of the creating subject in arguments like (\ref{csarg}).   This is to say that Kreisel hoped to show that Brouwer's weak counterexamples could be turned into \textsl{strong} ones by presenting axioms which allow for the formal refutation of principles such as (\ref{ineqapart}).\footnote{Or as he puts it: `For truly foundational research it is of special interest to derive a purely mathematical assertion from axioms concerning a specifically intuitionistic notion, here: the thinking subject $\ldots$ [A]part from these mathematical consequences, one wants to formulate as fully as possible properties of these basic notions: one learns more about them by getting contradictions $\ldots$ than by trying to avoid the notions!' \citeyearpar[p. 159]{Kreisel1967b}.}

The strong counterexample which Kreisel hoped to provide involved the refutability not just of (\ref{ineqapart}) itself but rather of a form of \textsl{Markov's Principle} from which it can be derived.   The latter is often taken to correspond to the following first-order schema:
\begin{example}[(MP)]  
$\forall x (A(x) \lor \neg A(x)) \to (\neg\neg \exists x A(x) \to \exists x A(x))$
\end{example}
But in fact the statement which figured in Kreisel's creating subject argument is not (MP) itself but rather the following second-order generalization:
\begin{example}[(GMP)]  
$\forall \alpha (\neg\neg \exists x \alpha (x) = 0 \to \exists x \alpha (x) = 0)$
\end{example}
Here the variable $\alpha(x)$ is intended to range over choice sequences, which in turn suggests that the  interpretation of GMP will depend on what principles they are assumed to satisfy.  But if real numbers are defined as choice sequences of rationals as described above, it may be shown relative to the formalization of intuitionistic analysis which we have adopted below that GMP implies (\ref{ineqapart}).\footnote{See, e.g., \citep[p. 205]{Troelstra1988}.}  This in turn explains why Kreisel was interested in GMP rather than (\ref{ineqapart}) directly.  

It should finally be noted that already at the point of their introduction, MP and GMP were considered controversial.   On the interpretation which \citet{Markov1956} had himself proposed, $\alpha(x)$ can be understood as ranging over constructive functions.   In this context it might be thought that GMP can be justified on the basis of the following argument:  i) suppose we know that $\alpha(x)$ is computable by an algorithm $\mathfrak{A}$ and also we can prove that the assumption $\neg \E x \alpha(x) = 0$ leads to a contradiction; ii) then it is justifiable to assert $\E x \alpha(x) = 0$ since we know if we used $\mathfrak{A}$ to successively compute the values $\alpha(0), \alpha(1), \ldots$ we would eventually find an $n$ such that $\alpha(n) = 0$.   But on the other hand, \citet{Kreisel1958c,Kreisel1959c} had also shown that there are instances of MP which are not derivable in Heyting arithmetic or in one of Kleene's early axiomatizations of intuitionsitic analysis.   At the time Kreisel delivered the address on which \citeyearpar{Kreisel1967b} is based, both the interpretation of GMP and its compatibility with various formal systems were thus very much live issues.\footnote{Kreisel's own summary of the situation was as follows: `Kleene calls $\forall \alpha (\neg\neg \exists x \alpha (x) = 0 \to \exists x \alpha (x) = 0)$ a generalization of Markoff's principle, and Heyting said in the discussion that Markoff would formulate the principle only for constructive functions in place of free choice sequences $\alpha$. While this distinction is certainly valid, it seems too technical: Markoff's (implicit) interpretation of logical connectives is so mechanistic that any similarity to the intended interpretation is purely coincidental. [Theorem \ref{neggmp}] no more contradicts Markoff's principle as he understands it than it contradicts the \textsl{classical} reading of $\forall \alpha (\neg\neg \exists x \alpha (x) = 0 \rightarrow \exists x \alpha (x) = 0)$. -- This much, I believe, is clear; it is not so clear that the rules of intuitionistic mathematics are valid at all for Markoff's interpretation (if the latter is made explicit).' \citeyearpar[pp. 161-162]{Kreisel1967b}}

\subsubsection{Initial schematization}
\label{gmpis}
  
Kreisel begins the presentation of his creating subject argument by remarking that `The present section considers a striking use of a new primitive notion $\ldots$ to derive a purely mathematical assertion: $\forall \alpha (\neg\neg \exists x \alpha (x) = 0 \to \exists x \alpha (x) = 0)$' \citeyearpar[p. 158]{Kreisel1967b}.   As we will see, the method by which he refutes GMP can indeed be understood as a paradigmatic example of how   `intuitive notions which do not occur in ordinary mathematical practice' can lead to `to new axioms for current notions'  so as to make such an analysis of such notions `as precise as possible' while eliminating  `doubtful properties of the intuitive notions when drawing conclusions about them'.\footnote{In fact Kreisel later remarked that  `[O]ne of the main purposes of the analysis is to restrict the notion of thinking subject so as to eliminate \emph{accidental} psychological elements, yet to exploit essential ones' \citeyearpar[p. 159]{Kreisel1967b}.}   Kreisel's overall goal can thus be understood as that of replacing the contingent aspects of argument (\ref{csarg}) with reasoning which fulfills the goal of informal rigour `not to leave undecided questions which can be decided by full use of evident properties of these intuitive notions'.   
 
Kreisel's first step was to introduce the notation $\Sigma \vdash_n A$ with the intended interpretation `the creating subject $\Sigma$ has evidence for asserting $A$ at stage $n$' (p. 159).  He states the rationale for this most clearly in his reply to Heyting:

{\footnotesize
\begin{quote}
[V]ery little of the `thinking subject' is used in the derivation. Instead of writing $\Sigma \vdash_n A$ I could write $\Sigma_n \proves A$ and read it as: the $n$th proof establishes $A$. In other words, the essential point would not be the individual subject, but the idea of proofs arranged in an $\omega$-order (each proof of course being a mental, not necessarily finite, object on the intuitionistic conception).  The idea is that one would not make use of any empirical information about the order  in which people come to think of proofs. Also, the sequence $\Sigma_n$ is not itself considered to be given by a rule. \hfill \citeyearpar[p. 179]{Kreisel1967b}
\end{quote}
}

Here Kreisel begins by observing that his derivation of the negation of GMP in fact requires taking into account only a few properties of how the creating subject operates.   These include the fact that they may be understood as constructing proofs at discrete stages which can be ordered $0,1,2,\ldots$ as an $\omega$-sequence.   As Kreisel's second proposed notation suggests, this makes it possible to understand the operation of the creating subject as giving rise to a sequence of monotonically increasing sets of propositions 
$\Sigma_0, \Sigma_1, \ldots$ such that $\Sigma_i$ contains the statement which they have proven (or otherwise have obtained evidence for) at stage $i$ -- i.e. so that for purposes of the argument, the subject is identified with their current state of knowledge.\footnote{Many interpretative details arise at this stage on which Brouwer himself is largely silent -- e.g. Is there one creating subject or many?  Must the stages in the operations of the (or a) creating subject be ordered as a sequence of order-type $\omega$?   Is it possible for the creating subject to acquire evidence for more than one proposition at a given stage?  Should the sets $\Sigma_n$ be understood to be closed under deductive consequence?   Such issues have been discussed in detail subsequent authors such as  \citet{Troelstra1969}, \citet{Dummett2000}, \citet{Van-Atten2018}.   But as the attendant axiomatic choices do not bear on the reconstruction of Kreisel's argument, we will pass over them here.}  

Although these considerations illuminate the sense which Kreisel wished to assign to $\Sigma \proves_n A$ (or $\Sigma_n \vdash A$), we will employ the notation $\Box_n A$ which has become conventional after its introduction by \citet{Troelstra1988}.   The axioms of what has come to be known as \textsl{Kreisel's Theory of the Creating Subject} $[\mathsf{CS}]$ can now be set out as follows:
\begin{examples}
\item[(CS1)] $\Box_n A \lor \neg \Box_n A$
\item[(CS2)] $A \to \neg\neg \exists n (\Box_n A )$ 
\item[(CS3)] $\exists n (\Box_n A) \to A$ 
\end{examples}
In \citeyearpar{Kreisel1967b} Kreisel provides little individual discussion of these principles outside of their use in his formal refutation of GMP.   But as they have been extensively discussed by subsequent commenators, it will suffice to give the following brief indications of their intended justification:
\begin{enumerate}[i)]
\item CS1 formalizes the decidability of statements of the form `the creating subject has evidence for $A$ at stage $n$'.  This finds justification in the traditional view that the proof of a statement should be intuitively recognizable -- or as Kreisel put it, `we can recognize a proof when we see one'  \citeyearpar[p. 202]{Kreisel1962a}.
\item CS2 is reported by \citet[p. 296]{Myhill1967} to have been called by Kreisel  the axiom of \textsl{Christian charity}  `because it says the only grounds  we could have for asserting that a proposition would never be proved are that we already know it to be absurd -- and not e.g. that people are too stupid'.   
\item CS3 can be understood as expressing the \textsl{soundness} of the notion of constructive proof or evidence at issue as it takes the form `if the creating subject has evidence for $A$ at stage $n$, then $A$ is true'.  Kreisel adopts an analogous \textsl{reflection principle} in his theory of constructions \citeyearpar[p. 126]{Kreisel1965a}.
\end{enumerate}

The process of justifying these principles can be understood to correspond to stage Ic in the schema \ref{irs} for informally rigorous arguments we have proposed in \S\ref{3irs}.   For as we have seen, Kreisel regarded the creating subject itself as a \textsl{novel} concept which (as he puts it) `does not occur in ordinary mathematical practice'.   His reflection on Brouwer's use of this notion can thus be understood as leading to the identification of the axioms just stated as constitutive principles  which relate the creating subject to statements in the precise language of intuitionistic analysis.  The vocabulary Kreisel introduced involving the creating subject (i.e. $\Box_n$) can thus be taken as comprising the novel language $\mathcal{L}_N$ while the axioms $\mathsf{CS} = \{\mathrm{CS1}, \mathrm{CS2}, \mathrm{CS3}\}$ comprise the set $\Gamma^K_2$ in the Kreselian language $\mathcal{L}_K$ formed by adjoining the creating subject vocabulary to that of intuitionistic analysis.   Finally since the common language $\mathcal{L}_C$ -- and thus also the sets $\Gamma_1, \Gamma^C_2$ and $\Gamma^J_2$ -- will be empty in this case, this means that $\Gamma_2 = \mathsf{CS}$.\footnote{This partitioning of languages and principles flags an important difference between Kreisel's validity and creating subject arguments.  As we have seen in \S\ref{validity}, the former is intended to provide a mathematical analysis of a concept -- i.e. logical validity -- which Kreisel takes to be a genuinely \textsl{common} (but initially imprecise) component of our mathematical practice.  On the other hand, he viewed Brouwer's use of the creating subject as a genuinely \textsl{novel} incursion into mathematical reasoning aimed at clarifying the status of a statement regarding the intuitionistic continuum -- i.e. (\ref{ineqapart}) -- which is itself stated in a purely mathematical language.   Thus relative to the way we are employing the terminology, the analysis of \textsl{common} concepts does not figure in the creating subject argument.}

\subsubsection{From schematization to formalization}  
\label{gmpsf}

Carrying out phase II in the application of the informal rigour schema IR to Kreisel's creating subject argument poses two challenges both of which involve the formulation of the appropriate Kreiselian theory $\mathsf{T}_K = \mathsf{T}_P \cup \Gamma_2$.    The first pertains to the status of the notation which is used to reason about the creating subject within the language $\mathcal{L}_K$ of this theory.   As we have just seen, Kreisel himself appears to have vacillated between treating $\Sigma_n$ as an object-language variable over a class of creating subjects and as a meta-variable over sets of sentences for which the subject has obtained evidence.   On the other hand, the notation $\Box_n A$ which has been adopted by subsequent commentators has the appearance of a family of propositional operators.

These obserations suggest that it is not entirely straightforward to formalize Kreisel's argument using the conventional notation of first- or higher-order logic.  But there is also no reason to expect that these challenges cannot be overcome in an essentially routine manner -- e.g.  by either employing a form of temporal logic with variables for stages or a so-called modal predicate of the form $S(n,\ulcorner A \urcorner)$.  But rather than pausing to develop the details of such possibilities here, we will again follow the treatment of the Theory of the Creating Subject given in \citet[\S 4.9.2]{Troelstra1988} which directly employs the notation $\Box_n A$ in its object language as well as allowing for numerical quantifiers over stages to bind the index of this operator.\footnote{This choice also allows us to sidestep two concerns which might arise about the consistency of a theory which combines mathematical axioms with principles formalizing the operation of the creating subject.  For on the one hand, the sorts of concerns raised by Montague's \citeyearpar{Montague1963} well-known inconsistency results involving so-called \textsl{modal} (or \textsl{epistemic}) \textsl{predicates} are avoided by treating $\Box_n$ as a propositional operator when taken in conjunction with van Dalen's \citeyearpar{Van-Dalen1982a} proof that Kreisel's axioms $\mathsf{CS}$ are conservative over Heyting arithmetic.   And on the other, by taking these axioms as the \textsl{only} principles which are assumed to hold of the creating subject and also restricting function comprehension in the manner of the system $\mathsf{FIM}^+_0$ below, we can also avoid the problems caused by the existence of `self-referential' choice sequences illustrated by what is known as \textsl{Troelstra's paradox} (see, e.g., \citealp[p. 845]{Troelstra1988a}, \citealp{Van-Atten2016}).}

Kreisel states that the refutation of GMP from the creating subject axioms he gives on \citeyearpar[p. 160-161]{Kreisel1967b} is carried out using `current intuitionistic axioms'.   But as he is not explicit about this, a second obstacle which must be overcome is that of choosing an appropriate theory of intuitionistic analysis.     At the time of his original address, two axiomatizations had been proposed, respectively by Kreisel himself in \citeyearpar[\S 2.5]{Kreisel1965a} and by \citet{Kleene1965}.  As Kreisel repeatedly cites the latter, we will employ a fragment $\mathsf{FIM}_0$ of Dragalin's \citeyearpar[\S 4.1]{Dragalin1988} concise reformulation of Kleene and Vesley's theory which he introduced under the name $\mathsf{FIM}$.\footnote{This name abbreviates the title of Kleene and Vesley's book \textsl{Foundations of Intuitionistic Mathematics} \citeyearpar{Kleene1965} whose first chapter contains their own more extensive presentation of this system.}

The language of $\mathcal{L}_{\mathsf{FIM}_0}$ of $\mathsf{FIM}_0$ consists of that of first-order arithmetic with $x,y,z, \ldots$ as numerical variables together with names $s, t, \ldots$ for primitive recursive functions and additionally variables $\alpha, \beta, \gamma, \ldots$ and quantifiers ranging over choice sequences of type $\mathbb{N} \rightarrow \mathbb{N}$.  Building on the axioms of intuitionistic second-order logic, the mathematical axioms of $\mathsf{FIM}_0$ can be divided into two groups, the first of which  Dragalin calls \textsl{Primitive Recursive Analysis}:
\begin{examples}
\item[($\mathrm{PrAn}_1$)] The axioms of first-order Heyting Arithmetic, inclusive of induction in the full language of $\mathsf{FIM}_0$, the identity axiom $\forall x \forall y(x=y \rightarrow \alpha(x) = \beta (y))$, and the defining equations of all primitive recursive terms.
\item[($\mathrm{PrAn}_2$)] The \textsl{primitive recursive closure} (or comprehension) scheme $\exists \alpha \forall x (\alpha (x) = t(x))$, where $t(x)$ is any term of $\mathcal{L}_{\mathsf{FIM}_0}$ which does not contain $\alpha$ free.
\end{examples}

In order to state the final axiom of $\mathsf{FIM}_0$, it is useful to introduce the following conventional definitions involving choice sequences:  
\begin{example}
\begin{enumerate}[i)]
\item $\bar{\beta}(x) =_{\mathrm{df}} \langle \beta (0), \beta (1), \ldots, \beta (x-1)\rangle$
\item $(x)_i$ denotes the $i$th component of $x$  when viewed as a code for a finite sequence in some conventional manner and $x \star y$ denotes the code of the sequence formed by concatenating the sequences coded by $x$ and $y$  
\item $\alpha \in K_0 =_{\mathrm{df}} \forall x \forall y (\alpha (x) \neq 0 \to \alpha (x \star y) = \alpha (x)) \wedge \forall \beta \exists x (\alpha (\bar{\beta}(x)) \neq 0)$  
\item $(y = \alpha (\beta)) =_{\mathrm{df}} \exists z (y + 1 = \alpha (\bar{\beta}(z)))$ 
\end{enumerate}
\end{example}
$K_0$ is traditionally described as the class of \textsl{continuous functionals}.  This is to say that if $\alpha \in K_0$ then $\alpha$ is either equal to $0$ for all $x$ or is such that not only does it stabilize to a value $y > 0$, but that when viewed as an operation on initial segments of another choice function $\beta$ there exists an $x$ such that the initial segment $\overline{\beta}(x)$ fixes these values.\footnote{See, e.g., \citep[\S 4.1]{Dragalin1988} or \citep[\S 4]{Troelstra1988} for further discussion of the history and motivation of continuity principles in intuitionistic analysis.}

The final axiom of $\mathsf{FIM}_0$ is what is known as \textsl{Brouwer's continuity principle for numbers} (or $\forall \alpha \E n$-continuity):
\begin{example}[(BC-N)]
$\forall \alpha \exists x A(\alpha, x) \rightarrow \exists \gamma \in K_0 \forall \alpha A(\alpha, \gamma (\alpha))$
\end{example}
This principle can be understood as expressing that if $\forall \alpha \exists x A(\alpha, x)$ holds, the value of the number $x$ must depend continuously on an initial segment of the choice function $\alpha(x)$.   Suppose we now let  $\mathsf{FIM}_0 = \mathrm{PrAn}_1 + \mathrm{PrAn}_2 +$ BC-N.   This theory can already be seen to be in conflict with classical logic as is illustrated by the following result which figures implicitly in Kreisel's refutation of GMP:\footnote{See, e.g.,  \citep[p. 84]{Kleene1965} or \citep[p. 209]{Troelstra1988} where it is shown that the so-called \textsl{Generalized Law of the Excluded Middle} is refutable from a weaker form of continuity derivable from BC-N.}
\begin{prop}
$\mathsf{FIM}_0 \proves \neg \forall \alpha(\A x (\alpha(x) = 0) \mor \A x(\alpha(x) \neq 0))$.
\label{negglem}
\end{prop}

In our reconstruction of Kreisel's argument, $\mathcal{L}_{\mathsf{FIM}_0}$ will correspond to the precise language $\mathcal{L}_P$ and $\mathsf{FIM}_0$ will correspond to the precise theory $\mathsf{T}_P$.  In order to extend $\mathcal{L}_{\mathsf{FIM}_0}$ and $\mathsf{FIM}_0$ to a language and theory appropriate for formalizing Kreisel's argument several additional steps must be taken.   First, we extend the definition of the class of formulas and sentences of  $\mathcal{L}_{\mathsf{FIM}_0}$ in the standard manner to include formulas of the form $\Box_x A$ and also $\E x \Box_x A$ and $\A x \Box_x A$  (here $x$ is a numerical variable which binds the index of $\Box_x$ which we will denote using the range $n,m,p,\ldots$).   Second, we define a formula in this language to be \textsl{extended}-$\Delta^0_0$ just in case it is composed using propositional connectives and bounded first-order quantifiers from atomic arithmetical statements and statements of the form $\Box_n B$ where $n$ is a free variable and $B$ is a sentence of $\mathcal{L}_{\mathsf{FIM}_0}$.\footnote{The intention is that extended $\Delta^0_0$-formulas will be \textsl{decidable} in the theory $\mathsf{FIM}^+_0$ defined below either in the conventional sense which holds for $\Delta^0_0$ arithmetical formulas in Heyting Arithmetic or in virtue of the decidability of formulas of the form $\Box_n B$ which follows as a consequence of CS1 in intuitionistic logic.}   We also define an extended language   $\mathcal{L}^+_{\mathsf{FIM}_0}$ consisting of $\mathcal{L}_{\mathsf{FIM}_0}$ together with a class of new atomic function symbols $\chi_{A}(x)$ for all extended-$\Delta^0_0$ formulas $A(x)$.

We next extend $\mathrm{PrAn}_1$ and $\mathrm{PrAn}_2$ respectively as follows:
\begin{examples}
\item[($\mathrm{PrAn}^+_1$)] In addition to the axioms included in $\mathrm{PrAn}_1$ we also include all statements of the form 
$$\A x(\chi_{A}(x) = 1 \leftrightarrow A(x) \wedge \chi_{A}(x) = 0 \leftrightarrow \neg A(x))$$
for $A(x)$ an extended-$\Delta^0_0$ formula expressing that $\chi_{A}(x)$ defines the characteristic function expressed by this formula.
\item[($\mathrm{PrAn}^+_2$)] The \textsl{primitive recursive closure} scheme $\exists \alpha \forall x (\alpha (x) = t(x))$ is extended to allow that $t(x)$ may be any term of the language $\mathcal{L}_{\mathsf{FIM}_0}^+$ which does not contain $\alpha$ free.
\end{examples}
Finally let $\mathsf{FIM}^+_0 =  \mathrm{PrAn}^+_1 + \mathrm{PrAn}^+_2 +$ BC-N.   

With these definitions in place, we can now define the Kreiselian language to be $\mathcal{L}_K$ to be $\mathcal{L}^+_{\mathsf{FIM}_0}$ and the Kreiselian theory to be $\mathsf{T}_K = \mathsf{FIM}^+_0 + \mathsf{CS}$ -- i.e. a fragment of Dragalin's theory $\mathsf{FIM}$ extended in the indicated manner to the language containing $\Box_n$ together with Kreisel's creating subject axioms $\mathsf{CS}$.   Our next goal is to demonstrate that this is indeed sufficient to refute GMP by showing the following:\footnote{In the following proof we have adhered to the structure of the argument which Kreisel sketches on \citeyearpar[pp. 160-161]{Kreisel1967b} by retaining the sequences of his main claims and supporting statements (we have respectively labeled these as Claims i) - iii) and  a) - e).  The task of providing an exact reconstruction is complicated by Kreisel's omission of many steps (some of which we have attempted to fill in as 1 - 12) and several apparent typographical errors.   As we will discuss in \S\ref{gmprec}, several results in the vicinity of Theorem \ref{neggmp} are now well-known.  But we are unaware of any prior attempt to verify Kreisel's original proof.}

 \newcounter{csc}
\begin{theorem}
 $\mathsf{FIM}^+_0 + \mathsf{CS} \proves \neg \forall \alpha (\neg\neg \exists x \alpha (x) = 0 \to \exists x \alpha (x) = 0)$.
 \label{neggmp}
\end{theorem}

{\small
\begin{proof}

We first define a predicate $P(\beta, \alpha, m)$ as follows:
$$P(\beta, \alpha, m) =_{\mathrm{df}} (\alpha(m) = 0 \leftrightarrow [(\exists x < m) (\beta(x) \neq 0) \vee \Box_m \forall x (\beta(x) = 0)])$$

\noindent We now reason in  $\mathsf{FIM}^+_0 + \mathsf{CS}$ via the following question of claims.\\

\noindent \underline{Claim i)} $\forall \beta \forall \alpha [ \forall m P(\beta, \alpha, m) \rightarrow \neg\neg \exists  n (\alpha(n) = 0)]$\\

\noindent \textsl{Proof of Claim} i): We show this via the following two subclaims:\\

\noindent  Subclaim i1) : If $\forall m P(\beta, \alpha, m)$ and $n \in \mathbb{N}$ is such that $\beta (n) \neq 0$, then $\alpha (n+1) = 0$.\\

\noindent \textsl{Proof of subclaim} i1).  Suppose $\forall m P(\beta, \alpha, m)$ and $\beta(n) \neq 0$.  Since $n < n+1$, $\exists x < n+1 (\beta (x) \neq 0)$ and thus also $\exists x < (n+1) (\beta(x) \neq 0) \vee \Box_{n+1} \forall x (\beta(x) = 0)$.  Thus by the assumption that $\forall m P(\beta, \alpha, m)$ holds, it follows that $P(\beta, \alpha, n+1)$ and thus also $\alpha (n+1) = 0$.  \\

\noindent Subclaim i2) : If $\forall m P(\beta, \alpha, m)$ and $\forall m (\alpha(m) \neq 0)$, then $\forall x (\beta(x) = 0)$.\\

\noindent \textsl{Proof of subclaim} i2).    Suppose $\forall m P(\beta, \alpha, m)$ and $\forall m (\alpha (m) \neq 0)$.  Then in particular, the latter implies $\alpha (n+1) \neq 0$.  Hence, by Subclaim i1) it follows that $\forall m P(\beta, \alpha, m) \mand \forall m (\alpha(m) \neq 0) \rightarrow \neg \beta(n) \neq 0$ -- i.e. $\forall m P(\beta, \alpha, m) \mand \forall m (\alpha(m) \neq 0) \rightarrow \neg \neg \beta(n) = 0$.   But since $\beta(n)$ is a fixed natural number, it follows from the decidability of numerical equality in $\mathsf{FIM}_0$ that $\neg \neg \beta(n) = 0$ implies $\beta(n) = 0$.   But since $n$ was arbitrary, it follows that $\forall m P(\beta, \alpha, m)$ and $\forall m (\alpha(m) \neq 0)$ jointly imply $\forall x (\beta (x) = 0)$ as desired.\\

\noindent Returning to the proof of Claim i), note that as an instance of axioms CS2 we have
\begin{example}[(a)]
$\forall x \beta(x) = 0 \rightarrow \neg\neg \exists n (\Box_n \forall x \beta(x) = 0)$
\end{example}

\noindent By the definition of $P(\beta, \alpha, m)$ and predicate logic we have 

\begin{example}[(1)] $\exists p (\alpha(p) = 0) \leftrightarrow [\exists p (\exists x < p)(\beta(x) \neq 0) \lor \exists p \Box_p \forall x (\beta(x) = 0))] $ 
\end{example}

\noindent From this we can derive 

\begin{example}[(2)]  $ \neg\neg [\exists p (\exists x < p)(\beta(x) \neq 0) \lor \exists p \Box_p \forall x (\beta(x) = 0))] \to  \neg\neg \exists p (\alpha(p) = 0) $ 
\end{example}

\noindent On the other hand, by  propositional logic 

\begin{example}[(3)]  $ \neg\neg \exists p \Box_p \forall x (\beta(x) = 0) \to \neg\neg [\exists p (\exists x < p)(\beta(x) \neq 0) \lor \exists p \Box_p \forall x (\beta(x) = 0))] $
\end{example}

\noindent  Thus (4)  $\neg\neg \exists p \Box_p \forall x (\beta(x) = 0) \to  \neg\neg \exists p (\alpha(p) = 0)$.\\

\noindent By Subclaim i2), together with statements a) and (4)

\begin{example}[(b)]
$ \forall m P(\beta, \alpha, m) \wedge \forall m (\alpha(m) \neq 0) \to  \neg\neg \exists p (\alpha(p) = 0) $
\end{example}

\noindent Then
(5) $\forall m P(\beta, \alpha, m) \wedge \neg \exists p (\alpha(p) = 0) \rightarrow  \neg\neg \exists p (\alpha (p) = 0)$. \\

\noindent Finally by propositional logic
 
\begin{example}[(c)]
$\forall m P(\beta, \alpha, m) \rightarrow \neg\neg \exists p (\alpha (p) = 0) $
 \end{example}

\noindent as desired to establish Claim i). \hfill $\boxtimes$ \\

\noindent \underline{Claim ii)} GMP implies that  $\forall \beta \exists m' [(\exists x <m') (\beta(x) \neq 0) \vee \Box_{m'} \forall x (\beta(x) = 0)]$.\\

\noindent \textsl{Proof of Claim} ii).  By Claim i) and GMP

\begin{example}[(d)]
$\forall \beta \forall \alpha [\forall m P(\beta, \alpha, m) \to \exists m' (\alpha (m') = 0)] $
\end{example}

\noindent By the definition of $P(\beta, \alpha, m)$

\begin{example}[(e)]
$\forall \beta \forall \alpha [ \forall m P(\beta, \alpha, m) \rightarrow \exists m' ((\exists x < m') (\beta(x) \neq 0) \vee \Box_{m'} \forall x (\beta(x) = 0))] $
\end{example}
and thus
(6) $(\forall \beta) [ \exists \alpha \forall m P(\beta, \alpha, m) \to \exists m' ((\exists x < m') (\beta x \neq 0) \lor \Box_{m'} \forall x (\beta x = 0))]$.
\vspace{1em}

\noindent Subclaim ii1): $\exists \alpha \forall m P(\beta, \alpha, m)$ --  i.e.

\begin{example}[(7)] $\exists \alpha \forall m [\alpha(m) = 0 \leftrightarrow  ((\exists x < m) (\beta(x) \neq 0) \lor \Box_m \forall x (\beta(x) = 0))]$
\end{example}

\noindent \textsl{Proof of Subclaim} ii1): Note that formula $A(m) =_{\mathrm{df}} (\exists x < m) (\beta(x) \neq 0) \lor \Box_m \forall x (\beta(x) = 0)$ on the right-hand side of the matrix of (7) is extended-$\Delta^0_0$.    It thus follows from $\mathrm{PrAn}^+_1$ that there is a term $\chi_A(m)$ which defines the characteristic function of the predicate $A(m)$ over $\mathsf{FIM}^+_0$.  But then the existence of a function $\alpha(m)$ witnessing (7) follows from $\mathrm{PrAn}^+_2$.\footnote{Kreisel is not explicit about the principles which he takes to justify the existence of $\alpha(m)$.   Statement (7) could also be obtained directly by appeal to a principle such as the Kripke scheme (as discussed in \S\ref{gmprec}).  However Kreisel also stresses \citeyearpar[p. 161]{Kreisel1967b} that  what justifies this statement is the decidability of the relation $\Box_n$ (as codified by CS1) in the definition of the formula $A(m)$.}\\

\noindent Next note that Subclaim ii1) allows us to obtain 

\begin{example}[(8)] $\forall \beta \exists m' [(\exists x <m') (\beta(x) \neq 0) \vee \Box_{m'} \forall x (\beta(x) = 0)]  
\label{asd}$
\end{example}
under the assumption of GMP.   This establishes Claim ii). \hfill $\boxtimes$\\  

\noindent \underline{Claim iii)} \ $\mathsf{FIM}^+_0 + \mathsf{CS} \vdash \forall \beta \exists m' ((\exists x <m') (\beta(x) \neq 0) \lor \forall x (\beta(x) = 0))$.\\

\noindent This follows immediately from (8) the axiom CS3 and logic. \hfill $\boxtimes$\\

\noindent We finally aim to derive a contradiction from Claim iii) via the axiom BC-N of $\mathsf{FIM}_0$.\footnote{In his original proof, Kreisel simply states that (8) `contradicts continuity'.  But as we will see below,   one way to reconstruct what he seems to have meant is that BC-N can be used to obtain Proposition \ref{negglem}.}  To this end, define $A(\beta, m^\prime) =_{\mathrm{df}} (\exists x <m') (\beta(x) \neq 0) \vee \forall x (\beta(x) = 0)$.   Note that still under the hypothesis of GMP, we have from (8) that $\A \beta \E x A(\beta,x)$.   It thus follows from the relevant instance of BC-N, that $\exists \gamma \in K_0 \A \beta A(\beta,\gamma(\beta))$ -- i.e. $\exists \gamma \in K_0 \forall \beta  ((\exists x < \gamma (\beta)) (\beta(x) \neq 0) \vee \forall x (\beta(x) = 0))$.  This is in turn equivalent to 

\begin{example}[(9)] $\exists \gamma \in K_0 \forall \beta (\exists x(x  < \gamma (\beta) \land \beta x \neq 0) \lor \forall x (\beta x = 0))$
\end{example}

\noindent Letting $\gamma_0 \in K_0$ be an appropriate witness, we then have

\begin{example}[(10)]  $\forall \beta(\exists x(x  < \gamma_0(\beta) \land \beta(x) \neq 0) \lor \forall x (\beta(x) = 0)) $
\end{example}

\noindent But this implies (11)  $\forall \beta(\exists x \neg(\beta(x) = 0) \vee \forall x (\beta(x) = 0))$ which is in turn implies

\begin{example}[(12)]  $\forall \beta (\neg \A x(\beta(x) = 0) \vee \forall x (\beta(x) = 0))
\label{glem}$
\end{example}

\noindent Note finally that (12) corresponds to the \textsl{Generalized Principle of the Excluded Middle} which is refutable in $\mathsf{FIM}^+_0$ by Proposition \ref{negglem}.   Formalizing the foregoing reasoning as a \textsl{reductio}, we are thus able to derive the negation of GMP in $\mathsf{FIM}^+_0 + \mathsf{CS}$ as desired.   

\end{proof}

}
 
Since we have taken $\mathsf{FIM}^+_0 + \mathsf{CS}$ as the Kreiselian theory $\mathsf{T}_K$, Theorem \ref{neggmp} constitutes the relevant informally rigourous refutation of GMP.    It appears that Kreisel initially took this result as settling the question about the tenability of Markov's  Principle from Brouwer's perspective which he had set out to answer.    This is evident for instance from one of his remarks in reply to Myhill:   

{\footnotesize
\begin{quote}
 The argument in the text gives a correct deduction of $$\forall \beta [\neg\neg \exists x (\beta(x) = 0) \to \exists x (\beta(x) = 0)] \rightarrow \forall \beta [\exists x (\beta(x) \neq 0) \lor \forall x (\beta (x) = 0)]$$
by essential use of the axioms for $\vdash_n$ (essential, since the hypothesis is applied to an empirically defined sequence $\alpha$).  \hfill \citeyearpar[p. 185]{Kreisel1967b}
\end{quote}
}

As we will see in \S\ref{gmprec}, the situation quickly became more complicated in light of subsequent developments.    But  it may finally be noted that some additional principles beyond $\mathsf{FIM}_0$ are indeed needed to obtain a formal refutation of GMP.   For as we have already noted, Kreisel had previously shown that MP is independent of first-order Heyting arithmetic via a modified realizability interpretation.   Using their related method of \textsl{special realizability} \citet[p. 131]{Kleene1965} extended this to show that MP is additionally consistent with $\mathsf{FIM}$. \citet[p. 43]{Vesley1972} extended this further to show that GMP is additionally consistent with $\mathsf{FIM}$.    Putting this together with Kreisel's result, this immediately yields the following:
\begin{prop}  $\mathsf{FIM}^+_0 + \mathsf{CS}$ is a non-conservative extension of $\mathsf{FIM}_0$ $($and in fact of $\mathsf{FIM})$. 
\end{prop}

This observation confirms Kreisel's remark that the sort of extension of  $\mathsf{FIM}$ which is embodied by his creating subject axioms is `essential' for the refutation of GMP.   It also illustrates our prior claim at the end of \S\ref{3irs} that in the case of an informally rigorous argument carried out via the scheme \ref{irs}, the Kreiselian theory $\mathsf{T}_K$ should extend the precise theory $\mathsf{T}_P$ \textsl{non-conservatively}, thereby settling the status of a principle which cannot be resolved by `precise' or `common' reasoning alone.  We will see in \S\ref{ch} that this is a feature which distinguishes Kreisel's creating subject and CH arguments.

\subsubsection{The reception of the creating subject argument}
\label{gmprec} 
 
The foregoing reconstruction illustrates how Kreisel's refutation of GMP from his creating subject axioms conforms to the model of informal rigour we have described in \S\ref{3irs}.   Kreisel's argument also had a demonstrable effect on the development of intuitionistic analysis at the time it was originally presented.   Nonetheless, the role of Theorem \ref{neggmp} itself has subsequently been eclipsed by results obtained by Kripke and Myhill either in close parallel to or in light of Kreisel's presentation.  But on the other hand, the bearing of informal rigour on the interpretation of these results is also evident from the  exchanges between Kreisel, Heyting, and Myhill which are recorded at the end of  \citeyearpar{Kreisel1967b}.   

Kripke's contribution to these developments took the form of his introduction of what has come to be called the \textsl{Kripke Schema}.   This may be understood informally as a comprehension scheme for choice sequences whose values may be defined relative to an arbitrary proposition $A$.   The schema is typically formulated as having both a \textsl{weak} and \textsl{strong} form as follows:
\begin{examples}
\item[(wKS)] $\exists \alpha ((\neg A \leftrightarrow \forall x \alpha(x) = 0 \land (\exists x \alpha (x) \neq 0 \to A))$
\item[(sKS)] \ $\exists \alpha ( A \leftrightarrow \exists x (\alpha (x) \neq 0))$
\end{examples}
The priority for this principle with respect to Kreisel's creating subject axioms is itself unclear.\footnote{\citet[p. 497]{Kripke2019} reports that he formulated wKS in a letter to Kreisel written in `around 1965'.  But the first references to the Kripke scheme in print appear to occur in Myhill's reply to Kreisel \citeyearpar[p. 174]{Kreisel1967b} and in \citep[pp. 294-295]{Myhill1967}.}  But in terms of motivation, it is generally agreed that the justification of wKS (or sKS) are derived from the same considerations which lead to Kreisel's axioms $\mathsf{CS}$ -- i.e. if we imagine the creating subject carrying out their work at stages $n = 0,1,2,\ldots$, then we can imagine defining a choice sequence $\alpha(n)$ relative to statements involving whether the creating subject has proven a certain proposition at a given stage in the manner illustrated by the definition in Kreisel's proof of Theorem \ref{neggmp}.\footnote{Indeed \citet[p. 128]{Kreisel1970} would later remark that `The [Kripke] schema is \emph{justified} by reference to the `thinking subject' or, more objectively,  to the analysis of mathematics into $\omega$ stages'.  \citet[pp. 295-296]{Myhill1967} additionally remarks that Kreisel's axioms give `a deeper analysis' and also that `As long as we have a  sufficiently clear idea of the meaning of ``$A$ has been proved by the $n$th stage''  $\ldots$ to justify [the axioms CS1-CS3, the] derivation of [wKS] is a simple exercise'.}   

In fact it is straightforward to show over $\mathsf{FIM}_0$ that wKS is equivalent to the conjunction of the axioms $\mathsf{CS}$ and also that sKS is equivalent to the theory obtained by replacing CS2 with the following strengthened form:
\begin{example}[(sCS2)]
$A \to \exists n (\Box_n A )$ 
\end{example}
From this it immediately follows that GMP is also refutable in $\mathsf{FIM}_0 + \mathrm{wKS}$.   But as \citet[p. 128]{Kreisel1970} also notes, the latter principle is already stated in the language of intuitionistic analysis itself.   As a consequence, it has become more common to reconstruct Brouwer's original creating subject arguments in terms of wKS (or sKS) rather than engaging with the details of the `intensional' dimensions of Kreisel's formalization we will discuss below.\footnote{The discussions in  \citep[\S 16.3]{Troelstra1988a} and \citeyearpar[\S4.9.3]{Troelstra1988} are typical in this regard.} 

A result in which Kreisel took more immediate interest at the time derives from an observation of Myhill that a strengthening of the continuity assumption of the background theory over which the $\mathsf{CS}$ axioms are added leads to an outright inconsistency rather than a non-conservative extension of $\mathsf{FIM}_0$ in the manner of Theorem \ref{neggmp}.    In order to formulate Myhill's result, consider the principle known as \textsl{Brouwer's continuity principle for functions} (or $\forall \alpha \exists \beta$-continuity):\footnote{Here $\gamma|\alpha$ is defined by the condition  $\gamma | \alpha = \beta$ if and only if $\forall x(\lambda n.\gamma(\langle x \rangle \star n)(\alpha) = \beta(x))$ where $\langle x \rangle$ denotes the finite sequence coded by $x$.}
\begin{example}[(BC-C)]
$\forall \alpha \exists \beta A (\alpha, \beta) \rightarrow \exists \gamma \in K_0 \forall \alpha A (\alpha, \gamma | \alpha)$
\end{example}
Such a principle can be understood as an attempt to express relative to the intended constructive interpretation of function quantifiers that if $\forall \alpha \exists \beta A (\alpha, \beta)$ holds then the dependence of $\beta$ on $\alpha$ must be given by a continuous functional.  Myhill's result can now be stated follows:\footnote{In his comments on Kreisel's address (which were submitted afterwards) Myhill indicates the following proof: `Let $\beta (n)$ be $0$ until $\alpha$ is known to be rational and $1$ thereafter: take $A(\alpha, \beta)$ be $[(\forall x)(\beta (x) = 0 \leftrightarrow \neg (\alpha$ is rational$)] \wedge [(\exists x) (\beta (x) \neq 0 \to \alpha$ is rational$)]$; this is extensional but $\beta$ cannot be chosen to be a continuous function of $\alpha$' \citeyearpar[pp. 173-174]{Kreisel1967b}.  A more complete proof of the incompatibility of BC-C with wKS (and thus with $\mathsf{CS}$) is given by \citet[p. 135-136]{Dragalin1988}.}\begin{theorem} $\mathrm{PrAn} + \textnormal{BC-C} + \mathsf{CS}$ is inconsistent. \label{myhillthm} \end{theorem}
Both \citet[p. 135]{Kreisel1965a} and \citet[p. 73]{Kleene1965} had included a version of BC-C in their axiomatizations of intuitionistic analysis.   But BC-C was also a subject of debate at the time and in this case has also remained controversial.\footnote{See, e.g., \citep[\S 7]{Van-Atten2018}.}  Indeed much of Kreisel's exchange with Myhill is devoted to the question of how we should react to Theorem \ref{myhillthm}.  Kreisel introduced one of the central issues of their exchange as follows:

{\footnotesize 
\begin{quote}
Myhill's analysis suggests that there may have been a definite error, namely a failure to distinguish between extensional and intensional operations on free choice sequences. If anything, I myself was perhaps a little too faithful to Brouwer! I was indeed struck by the distinction which I lamely characterized as a difference between mathematical constructions on free choice sequences and those involving empirical concepts. But the sharper and better formulation is due to Myhill who also realized that, in view of the distinction, a new proof of $\neg \forall \alpha [\exists x 
\alpha(x) = 0 \vee \neg \exists x \alpha(x) =0]$, used in my derivation, is necessary. \\ \hspace*{1ex} \hfill \citeyearpar[p. 179]{Kreisel1967b}
\end{quote}
}

Relative to the way the terminology is employed here, Kreisel's creating subject axioms are `intensional' in the sense that they refer to the temporal stages in the operation of the creating subject while the Kripke scheme is `extensional' in the sense that reference to such stages is suppressed.   Once such a distinction is introduced, a further question is whether schematic principles like BC-N or BC-C should be restricted to disallow instances of the formula $A(\alpha,x)$ or $A(\alpha,\beta)$ that make use of vocabulary which references stages directly.  Kreisel and Myhill's discussion thus turned on the relative priority of preserving the practice of intuitionistic analysis, the different motivations which might be given for BC-N and BC-C, and whether intensional notions need to be represented in the object language in order to differentiate such principles.   

It is evident that not only were Kreisel and Myhill intrigued by such questions but also that they agreed they were illustrative of the sorts of issues to which informal rigour might be applied.  On the other hand, the mathematical morals which they drew from Theorem \ref{myhillthm} were framed tentatively.\footnote{With respect to the first point, Myhill concludes his comments with the following observation [his emphasis]: `\textsl{The whole dialectic of this chapter in mathematical philosophy is a delightful example of how our formalizations correct our intuitions while our intuitions shape our formalizations}'  \citeyearpar[p. 175]{Kreisel1967b}.   With respect to the latter point,  in \citeyearpar{Myhill1967} Myhill provided not only a proposal of an alternative axiomatization of intuitionistic analysis, but also a detailed argument for why BC-N should be accepted but BC-C rejected once intensional considerations about choices sequence are properly taken to account.  On the other hand, \citet[p. 137]{Troelstra1977} presents a different argument which appears to have largely discouraged further investigation of the role of intensionality in the subsequent development of intuitionistic analysis.}   Rather than providing a further reconstruction here, we will thus conclude this section by highlighting one of Kreisel's more programmatic conclusions:

{\footnotesize
\begin{quote}
[T]he actual practice of intuitionistic \textsl{mathematics} seems most elegantly formulated by restricting oneself to extensional operations.  So there may be some practical conflict here $\dots$ For foundations it is evident that intensional operations are fundamental, by the principle stressed throughout my paper: \textsl{extensional operators can be defined in terms of intensional ones, but not conversely}. \hfill \citeyearpar[p. 184]{Kreisel1967b} 
\end{quote}
}

The other example of the intensional/extensional distinction which Kreisel discusses in \citeyearpar[p. 143]{Kreisel1967b} is in regard to the `crude mixture' of concepts which he took contribute to the informal understanding of the notion of \textsl{set} -- i.e. `(i) finite sets of individuals $\ldots$ or (ii) set \textsl{of} something $\ldots$ but also (iii) properties or \textsl{intensions} where one has no \textsl{a priori} bound on the the extension (which are very common in ordinary thought but not in mathematics)'.   As we have suggested in \S\ref{validity}, Kreisel understood his validity argument to show that the coincidence of the class of logically valid sentences -- i.e. those true in all \textsl{structures} -- with those true in all set-sized \textsl{models} holds regardless of which of these informal concepts is used to precisify the notion of `structure'.   On the other hand, he also stresses that classical mathematics has evolved to prefer the second component of the mixture -- i.e. `sets of something' -- and also that Zermelo's axiomatization of set theory on this basis is `marvellously clear and comprehensive' (see \S\ref{set}).   This presumably includes the adoption of the Axiom of Extensionality for sets.

On the other hand, the considerations exposed by Kreisel's exchange with Myhill suggest that he came to view the development of intuitionistic mathematics differently.    For as Brouwer's creating subject arguments illustrate, this is a context in which intensional notions are employed, albeit in a manner which Kreisel acknowledges is \textsl{novel}.   Thus despite the reluctance expressed by Heyting \citeyearpar[p. 173]{Kreisel1967b} to embrace such notions as a proper part of intuitionism or Myhill's tentative proposal \citeyearpar[p. 174]{Kreisel1967b} to replace the creating subject axioms with the Kripke scheme, Kreisel appears to have maintained that there is in fact a rationale for continuing to reason explicitly with intensional notions.   For if suppressing the vocabulary needed to formalize the relevant distinctions leaves us unable to decide the status of questions such as GMP, then we are failing to discharge the second duty of informal rigour `not to leave undecided questions which can be decided by full use of evident properties of these intuitive notions'.

\subsection{Mathematically definite problems and the Continuum Hypothesis}
\label{ch}

Kreisel presented what we have referred to above as his \textsl{CH argument} not only in the sources reviewed in \S\ref{context}, but on several other occasions.\footnote{The account given in Kreisel's biographic memoire of G\"odel \citeyearpar[\S III]{Kreisel1980} is one of the clearest and most extensive.}  But in addition to the usual textual complexities, the task of providing a clear account of the argument -- and why Kreisel chose to present it in the manner he did -- is also complicated by the contextual factors described in \S\ref{context}.   It will again be useful to expand briefly on the setting of his original presentations before considering the argument itself.

\subsubsection{Background}

Although Kreisel had worked in descriptive set theory during the 1950s, prior to the mid-1960s he had not engaged extensively with general set theory in the tradition of Cantor, Zeremelo, Fraenkel, Skolem, and G\"odel.  But as we have noted in \S\ref{context}, the first section of his survey of mathematical logic \citeyearpar{Kreisel1965a} provides an overview which was state-of-the-art for its time, inclusive of summaries of G\"odel's consistency proof for the axioms of Constructibility and Choice, the L\'evy-Shoenfield Absoluteness Lemma, and statements of some early results and open problems about measurable cardinals.   

In \citeyearpar[pp. 144-145]{Kreisel1967b} Kreisel also calls attention to the fact that the process by which Zermelo provided his original axiomatization of set theory by reflecting on the cumulative hierarchy of sets -- or as Kreisel calls it the \textsl{cumulative type structure} -- might itself be understood as an instance of informal rigour.\footnote{This proposal is developed further in \citep{Isaacson2011a}.}  Although he does not argue extensively for this point, his claims can be understood as descending both from his discussion of mathematical realism in \citeyearpar[\S 1]{Kreisel1967a} and also his earlier presentation of an axiomatic theory of sets and types in \citeyearpar[\S 1.1]{Kreisel1965a} which he had previously described as a reconstruction of `Zermelo's informal derivation of his axioms'.  We will postpone further consideration of this proposal until \S\ref{set} where we will describe an obstacle to assimilating his claims about Zermelo's formulation of axioms to the model of an informally rigorous \textsl{argument} we have proposed in \S\ref{oninfrig}.   

In addition to this, \citeyearpar[\S 1.6]{Kreisel1965a} also contained a detailed summary of the method of forcing by which Paul Cohen had proven the independence of the Continuum Hypothesis from the axioms of $\mathsf{ZF}$ in the spring of 1963.   A number of subsequent commentators have called attention to Kreisel's role (together with G\"odel, Feferman, and Scott) not only in anticipating similar techniques but also more directly in the discovery of Cohen's proof.\footnote{See \citep[\S 2]{Kripke1965}, Scott's introduction to \citep{Bell1977a}, \citep[p. 192-201]{Kreisel1980}, \citep{Moore1987}, \citep[\S XII.1]{Odifreddi1999a}, \citep{Cohen2002}, and \citep{Kanamori2008}.   Although a precise estimation of Kreisel's role remains elusive, the following two strands in his prior work on predicativity and intuitionism appear relevant.   The first strand originates with Kreisel's introduction in \citeyearpar[\S 6]{Kreisel1961a} of a relation similar to Cohen's original forcing definition as part of his attempt to analyze reasoning with what he refers to as \textsl{extensionally definite} terms.   This is in turn related to his proposed analysis in \citeyearpar{Kreisel1960} of the notion \textsl{predicative definability} which we will discuss further in \S \ref{preddef} below.  Feferman's paper \citeyearpar{Feferman1964a} -- which extends an abstract he delivered at the same 1963 conference where Cohen presented his results -- contains the related definition of \textsl{arithmetical forcing} -- i.e. the adaptation of Cohen's \citeyearpar{Cohen1963a} definition to the language of first- and second-order arithmetic.  Feferman states (p. 334) that the first theorem (2.12) he obtains using this definition -- i.e. that there exist hyperarithmetical sets which are not \textsl{implicitly definable} as the unique set satisfying a formula $\phi(X)$ in the language of first-order arithmetic with the single free second-order variable $X$  -- answers a question which was posed by Kreisel (see \citealp[p. 447, p. 452]{Rogers1987}, \citealp[Proposition 1.7, p. 294]{Odifreddi1983}).  In his proof, Feferman's proof employs a result of \citet[p. 307]{Kreisel1962c} which that states for all  $\phi(x,y) \in \Pi^1_1$, if $\A x \E y \phi(x,y)$, then there exists a hyperarithmetical function $f(x)$ such that $\A x \phi(x,f(x))$ (see \citealp[Lemma  2.6 p. 31]{Sacks1990}).   \citet[Theorem 2.10, p. 333]{Feferman1964a} used this to show the existence of  hyperarithmetic sets which generic relative to his forcing definition -- i.e. what are now called $\omega$-\textsl{generic} sets.   One of the central results of Feferman's paper (Theorem 3.15, p. 339) states that the structure $\mathrm{RA}_{\omega^{ck}_1}(A)$ obtained by adjoining a generic set $A$ to Kleene's ramified analytical hierarchy up to level $\omega^{ck}_1$  (see \S\ref{preddef}) preserves the property of satisfying hyperarithmetical comprehension hence also $\Delta^1_1$-comprehension by Kleene's Theorem (see \citealp[Proposition 6.2, p. 729]{Odifreddi1983b}, \citealp[Theorem 3.6, p. 96]{Sacks1990}).  The second strand in Kreisel's prior work relates to his introduction in \citeyearpar{Kreisel1958f} of \textsl{lawless sequences} (which Kreisel originally referred to as \textsl{absolutely free choice sequences}) to intuitionistic analysis.  Such a sequence is a function $\alpha: \mathbb{N} \rightarrow \mathbb{N}$ which is understood as being generated by a process constrained only by the restriction that at no stage a law-like restriction imposed on subsequent choices (see, e.g., \citealp{Troelstra1977a} for general discussion).  Kreisel proposed several axioms describing such sequences including what has come to be known as the principle of \textsl{Open Data} -- i.e. $\A \alpha \phi(\alpha) \rightarrow \E n(\alpha \in n \mand \A \beta \in n \phi(\beta))$ -- which states that if a formula $\phi$ holds of a lawless sequence $\alpha$, a finite initial segment $\alpha \upharpoonright n$ must already contain sufficient information to determine that this is the case.   In \citeyearpar[p. 109-110]{Kreisel1965a}, Kreisel analogized generic sets to lawless sequences and observed that Cohen's forcing definition can be obtained from Open Data and his other axioms for choice sequences.\label{forcingnote} }   Whatever attributional issues may have remained unresolved at this time, it is also evident that Kreisel felt that many of the figures who were present at the conference where he delivered \citeyearpar{Kreisel1967b} had derived the wrong moral from the result.    This is already evident in the passage we have reproduced at the end of \S\ref{words}.   But the extent of Kreisel's disdain for what he took to be the overly facile conclusions about the significance of the formal independence of CH drawn by his interlocutors is even more evident in his slightly later papers \citeyearpar{Kreisel1969a,Kreisel1971d}.   Here Kreisel directly attacks the formalist views expressed by \citet{Cohen1971} and \citet{Robinson1965a} which are not explicitly attributed in \citeyearpar{Kreisel1967b} -- in part by reiterating aspects of the argument for the definiteness of CH we are about to consider.

\subsubsection{Initial schematization}
\label{chis}

Kreisel's CH argument itself can be understood as proceeding in two basic steps:\footnote{The basic structure of the CH argument also remains consistent across Kreisel's presentations.   In the reconstruction below we have largely followed \citeyearpar[\S 3b]{Kreisel1967a} -- which in this case it is the most detailed of the treatments -- while also relying on  \citeyearpar[\S 1]{Kreisel1965a}, \citeyearpar[\S 1]{Kreisel1967b}, and \citeyearpar[\S III]{Kreisel1980} to clarify notational conventions and other details.}
\begin{example}
\begin{enumerate}[i)]
\item A precise definition is proposed for the concept of a \textsl{mathematically definite statement}.  
\item A mathematical argument is then presented which shows that the Continuum Hypothesis satisfies the proposed definition of  definiteness.   
\end{enumerate}
\label{charg}
\end{example}
A consequence of the second step is that CH possesses a definite truth value in the sense which Kreisel intended to highlight.   But it is also evident that his ultimate hope was to show that an additional informally rigorous argument could be mounted for one of the following:
\begin{example}
\begin{enumerate}[i)]
\item The Continuum Hypothesis is true.
\item The Continuum Hypothesis if false.
\end{enumerate}
\end{example}
In \S\ref{chnovel} we will discuss how Kreisel hinted that this might be possible by making use of \textsl{novel} set theoretic concepts. But he demurred from taking this step himself.\footnote{\citet[p. 196]{Kreisel1971d} would later assess the situation as follow:  `CH \textsl{is} decided by the full (second order) axioms of Zermelo $\ldots$ Our \textsl{present} analysis of Zermelo's axioms, that is the first order schemata in the usual language of set theory, is not sufficient to decide CH. Put succinctly: not the notion of set, but our analysis (present knowledge) of this notion is at fault.'}

Kreisel's point of departure in each of his presentations of the CH argument is that the original axiomatizations of arithmetic, analysis, and set theory provided by Peano, Dedekind, and Zermelo were categorical.  Upon stressing the fact that these axiomatizations are themselves given in \textsl{second-order languages}, he then makes the two following declarations in \citeyearpar{Kreisel1967a}:
\begin{enumerate}[i)]
\item `For usual axiomatic mathematics the notion [of mathematically definite problem] is confined to first-order statements as in [the validity argument]';
\item `To explain the notion of mathematically definite problem one needs the corresponding notion of second (or higher) order consequence.'  \hfill  \citeyearpar[p. 256]{Kreisel1967a}
\end{enumerate}

In accordance with his ultimate interest in CH, the definition of mathematical definiteness which Kreisel gave in \citeyearpar[p. 257]{Kreisel1967a} was restricted to sentences in the language of first-order set theory.   But to highlight the generality of this definition relative to some of the other examples he used to motivate it, it will be useful to state a more general version applicable to arbitrary finitely axiomatizable second-order theories formulated in relational languages.   To this end, let $\mathsf{T}^2$ be such a theory formulated in a second-order language $\mathcal{L}_{\mathsf{T}^2}$ containing the non-logical predicates $P_1,\ldots,P_k$ or arities $a_1,\ldots,a_k$.  Additionally let $\tau^2(X,Y_1,\ldots,Y_k)$ be the sentence resulting from restricting first-order quantifiers to $X$ and replacing $P_i$ with the second-order variable $Y_i$ in the conjunction of the axioms of $\mathsf{T}^2$.   Next let $\mathrm{Sat}^1_k(X,Y_1,\ldots,Y_k,\ulcorner \phi^1 \urcorner)$ denote a second-order definition of satisfaction for first-order sentences which formalizes the fact that the $\mathcal{L}_{\mathsf{T}^2}$-sentence $\phi^1$ holds when its first-order quantifiers are restricted to $X$ and $P_i$ is interpreted as $Y_i$.\footnote{The definability of such satisfaction for first-order formulas in a second-order language can be traced back to \citep{Mostowski1950} in the case of set theory and \citep[\S 5.2e]{Hilbert1939} in the case of arithmetic.} Finally, let $\models_2 \varphi^i$ ($i \in \{1,2\}$) be the conventional definition of second-order logical validity with respect to the so-called \textsl{standard} semantics wherein second-order quantifiers range over the \textsl{full powerset} of an appropriate Cartesian product of the domain -- e.g. \citep[\S 4.2]{Shapiro1991} or \citep[\S 5]{Dalen2008}.

In parallel to \citeyearpar[p. 157]{Kreisel1967a} we can now state

\begin{definition}
Let $\varphi^1$ be a first-order sentence of $\mathcal{L}_{\mathsf{T}^2}$ -- i.e. one containing no second-order quantifiers or free first-order variables.   Then $\varphi^1$ is said to be a \textnormal{consequence} of $\mathsf{T}^2$ just in case 
$$\models_2 \forall X \forall Y_1 \subseteq X^{a_1} \ldots \forall Y_k \subseteq X^{a_k} [\tau^2(X,Y_1,\ldots,Y_k) \rightarrow \\ \mathrm{Sat}^1_k(X,Y_1,\ldots,Y_k,\ulcorner \varphi^1 \urcorner)]$$
We then say that $\varphi^1$ is \textnormal{decided by} $\mathsf{T}^2$ just in case either $\varphi^1$ is a consequence of $\mathsf{T}^2$ or  $\neg \varphi^1$ is a consequence of $\mathsf{T}^2$.
\label{cddefn}
\end{definition}
\noindent Several complications apply to Kreisel's original formulation of this definition.\footnote{A first complication arises from the fact that although Kreisel's definition is otherwise stated formally in second order logic, he uses the word `satisfies' (in English) rather than employing a satisfaction predicate.  However not only does he describe a similar formal predicate in \citeyearpar[p. 155, App. A]{Kreisel1967b}, his exposition in \citeyearpar{Kreisel1967a} occasionally contains object language quantification over sentences which cannot be formalized without such a device.    Another potential confusion for modern readers arises in light of Kreisel's use of the symbol $\vdash_2$ in both \citeyearpar{Kreisel1967a,Kreisel1967b} to denote second-order logical consequence with respect to the standard semantics rather than derivability in a deductive system for second-order logic (as we have employed $\vdash_2$ here).}   But since these will not be relevant until we discuss the prospects for formalizing Kreisel's argument below, it will suffice for the moment to understand Kreisel's definition of \textsl{$\varphi^1$ is a consequence of $\mathsf{T}^2$} to coincide with the conventional definitions of $\mathsf{T}^2 \models_2 \varphi^1$ -- i.e. $\varphi^1$ is a second-order logical consequence of $\mathsf{T}^2$ with respect to the standard semantics -- and similarly for \textsl{$\varphi^1$ is decided by $\mathsf{T}^2$} -- i.e. either $\mathsf{T}^2 \models_2 \varphi^1$ or  $\mathsf{T}^2 \models_2 \neg\varphi^1$.  
 
For present purposes, it will also suffice to take Kreisel's basic axiomatization to be \textsl{second-order Zermelo set theory} $\mathsf{Z}^2$ -- i.e.  the $\mathcal{L}_{\mathsf{Z}} = \{\in\}$ theory consisting of the axioms of the familiar Zermel-Franekel axioms without the Axiom of Replacement but together with what Kreisel calls the \textsl{Axiom of Comprehensiveness} stating that the intersection of a class and a set is a set.\footnote{Kreisel's axiomatization of what he calls \textsl{Zermelo systems} is given on \citeyearpar[p. 256]{Kreisel1967b} and is based on the more detailed presentation in \citeyearpar[\S 1]{Kreisel1965a} several of whose conventions are then employed in \citeyearpar{Kreisel1967a,Kreisel1967b} without explanation.  A modern presentation can be obtained by the axiomatization given by \citet[p. 85]{Shapiro1991} as follows: i) replace $\in$ by Kreisel's non-logical membership predicate $E$; ii) omit Replacement;  iii) adjoin Comprehensiveness -- i.e. $\forall x \forall X \exists y \forall z[E(z,y) \leftrightarrow (X(z) \wedge E(z,x))]$.  (Recall also that the background deductive system for second-order logic is assumed to contain full impredicative comprehension.)}  It is a familiar fact that the smallest rank in the cumulative hierarchy -- or as Kreisel dubs it \textsl{the cumulative type structure} -- satisfying $\mathsf{Z}^2$ is $R(\omega + \omega)$.  Again simplifying slightly, Kreisel states the following:\footnote{In order to present a more streamlined version of the basic steps in Kreisel's CH argument, most commentators have followed \citet[p. 286]{Weston1976} by formulating Definition \ref{md} relative to $\mathsf{ZF}^2$ rather than $\mathsf{Z}^2$.   However this precludes Kreisel's use of the Axiom of Replacement (or more precisely, certain of its first-order consequences) as examples of \textsl{non-definite} statements.   On the other hand, in order to facilitate even more fine-grained scrutiny of such (weak) `axioms of infinity', \citet[p. 257]{Kreisel1967a} originally introduced another parameter in the definition by formulating it terms of the sequence of theories he called $\mathsf{Z}_{\alpha}$ which can be obtained $\mathsf{Z}^2$ by omitting the Power Set axiom and adjoining the formalization `$R(\alpha)$ exists' for $\alpha = \omega, \omega + 1, \omega + 2, \ldots$  This complexity is suppressed in \citeyearpar{Kreisel1967b} and will also not be relevant here.}

\begin{definition} For all $\varphi^1 \in \mathcal{L}^1_{\mathsf{Z}}$ -- i.e. first-order sentences in the language of set theory -- $\varphi^1$ is \textnormal{mathematically definite} just in case $\varphi^1$ is decided by $\mathsf{Z}^2$.  
\label{md}
\end{definition}

Kreisel's formulation of the Continuum Hypothesis was as follows (\citeyear[p. 267]{Kreisel1967a}, \citeyear[p.150]{Kreisel1967b}):
\begin{example}  For every set $x \subseteq R(\omega+1)$, there exists either an injective function from $x$ into $R(\omega)$ -- i.e. $|x| \leq \aleph_0$ -- or a bijective function between $x$ and $R(\omega+1)$ -- i.e. $|x| = 2^{\aleph_0}$.
\end{example}
It is not difficult to see that this statement is formalizable as an $\mathcal{L}^1_{\mathsf{Z}}$-statement $\chi^1$ in which all quantifiers are bounded by a definable term formalizing the definition of $R(\omega + 4)$.\footnote{For note that if $x \subseteq R(\omega + 1)$, then $x \in R(\omega + 2)$.   It thus follows that a function which maps $x$ either injectively into $R(\omega)$ or bijectively onto $R(\omega + 1)$ is a set of ordered pairs from $R(\omega + 2)$ which will themselves be members of $R(\omega + 3)$ (for the standard definition of pairing).   Such a function will thus be a set in $R(\omega + 4)$.} Suppose that in general we call an $\mathcal{L}^1_{\mathsf{Z}}$-sentence $\alpha$ \textsl{rank-bounded} just in case each of its quantifiers is bounded by such a term formalizing the definition of $R(\beta)$ for $\beta \leq \alpha$.   Kreisel then states a somewhat more general form of the following:
\begin{theorem}
All $\omega + \omega$ rank-bounded $\mathcal{L}^1_{\mathsf{Z}}$-statements are decided by $\mathsf{Z}^2$.
\label{kdef}
\end{theorem} 
Since we have just observed that CH can be expressed as the $\omega + 4$-bounded $\mathcal{L}^1_{\mathsf{Z}}$-statement $\chi^1$, we can then immediately conclude 
\begin{corollary} 
$\chi^1$ is decided by $\mathsf{Z}^2$ and thus mathematically definite.
\label{kch}
\end{corollary}

These results are consequences of the following more general fact which is now often referred to as \textsl{Zermelo's Quasi-Categoricity Theorem}:\footnote{Zermelo stated Theorem \ref{zqct} in \citeyearpar{Zermelo1930a} after which it was formalized in G\"odel-Bernays set theory by \citet{Shepherdson1952} who used it to show that CH cannot be refuted by an inner model construction in the manner of G\"odel's \citeyearpar{Godel1938} original consistency proof.   Following Kreisel's presentation of the foregoing argument in the mid-1960s -- and its subsequent popularization by theorists such as \citet{Weston1974,Weston1976}, \citet{Shapiro1991}, \citet{Potter2004}, and \citet{Isaacson2011a} -- Zermelo's result has become a focus of recent work in the philosophy of set theory in regard to the phenomena known as \textsl{internal categoricity} (to which we will return below).   But already in 1930 it was Zermelo himself who originally made the crucial observation on which Kreisel's CH argument depends: ``From this [\ref{zqct}] already follows $\ldots$ that Cantor's (generalized) conjecture $\ldots$ does not depend on the choice of the model, but that it is decided (as true or false) once and for all by means of our axiom system.'' \citeyearpar[p. 437]{Zermelo1930d}  (Strictly speaking, Zermelo's claim is too strong.   For Theorem \ref{zqct} entails that if $\mathfrak{M}_1 \models_2 \mathsf{ZF}$ and $\mathfrak{M}_2 \models_2 \mathsf{ZF}$ then $\mathfrak{M}_1 \cong R(\lambda)$ and $\mathfrak{M}_2 \cong R(\kappa)$ for strongly inaccessible $\lambda$ and $\kappa$.   But as it may be that (e.g.) $\lambda < \kappa$, it is possible that (say) $2^{\alpha} = \aleph_{\alpha+1}$ holds in both $\mathfrak{M}_1$ and $\mathfrak{M}_2$ for all infinite $\alpha < \lambda$ but fails in $\mathfrak{M}_2$ for larger $\alpha$.   Thus while the theorem does show that CH is definite in Kreisel's sense it does not suffice to show the definiteness of GCH.)}
\begin{theorem} Let $\mathsf{ZF}^2$ denote full second-order Zermelo-Fraenkel set theory $($i.e. $\mathsf{Z}^2$ together with the Axiom of Replacement$)$.   Then if $\mathfrak{M} = \langle M,E \rangle$ is such that $\mathfrak{M} \models_2 \mathsf{ZF}^2$, then $\mathfrak{M}$ is isomorphic $R(\kappa)$ where $\kappa$ is a strongly inaccessible cardinal.
\label{zqct}
\end{theorem}
\noindent It is not difficult to adapt the proof of Theorem \ref{zqct} to show that all full models of $\mathsf{Z}^2$ are isomorphic to $R(\lambda)$ where $\lambda$ is a limit ordinal greater than $\omega + \omega$.  Kreisel's Theorem \ref{kdef} then follows by observing that if $\phi^1$ is $\omega+\omega$-rank bounded, then it is absolute between $R(\omega + \omega)$ and $R(\alpha)$ for $\alpha \geq \omega + \omega$ -- i.e. $R(\omega + \omega) \models \phi^1$ iff $R(\alpha) \models \phi^1$.

If we now wish to understand the foregoing argument as an instance of informal rigour, an initial question is how the concepts figuring in the foregoing definitions should be segregated into common and precise vocabularies.  Some insight into how Kreisel thought about these matters can be gleaned by recalling his remark that mathematical definiteness is `confined to first-order statements'.   This may at first seem at odds with the fact that Kreisel repeatedly drew attention to the fact that the now-familiar first-order axiomatizations of number theory, analysis, and set theory have second-order antecedents which are categorical with respect to the standard semantics.  Such an orientation is reflected in the following oft-cited passage:
\begin{quote}
{\footnotesize The familiar classical structures (natural numbers with the successor relation, the continuum with a denumerable dense base etc.) are definable by second order axioms, as shown by Dedekind. Zermelo showed that his cumulative hierarchy up to $\omega$ or $\omega + \omega$ $\ldots$ is equally definable by second order formulae. Whenever we have such a second order definition there is associated a schema in first order form (in the language considered): For instance, in Peano's axiom $$\forall P(P(0) \wedge \forall x(P(x) \rightarrow P(x+1)) \rightarrow \forall x P(x))$$ 
\begin{spacing}{1}
one replaces the second order quantifier $P$ by a list of those $P$ which are explicitly defined in ordinary first order form $\ldots$ A moment's reflection shows that the evidence of the first order axiom schema derives from the second order schema $\ldots$ \hfill \citeyearpar[p. 148]{Kreisel1967b} \end{spacing}}

\end{quote}

Kreisel's own views notwithstanding, it is also evident that he understood the primacy of first-order axiomatizations over their second-order counterparts to reflect the entrenched conventional wisdom of mathematical logicians in the mid-1960s.\footnote{Kreisel discusses the mathematical and historical context of the distinction between first- and second-order axiomatizations at greater length in his reply to Mostowksi in \citep[pp. 97-103]{Lakatos1967} and in \citeyearpar[\S 4]{Kreisel1968}.}  As such, one way in which his CH argument can be approached is by regarding first-order languages, definitions, and axiomatizations as notion which are \textsl{common} in the sense the sense that they were both employed and accepted by the theorists whom Kreisel was intending to address with his CH argument.   On the hand, Kreisel can also be understood as regarding  second-order languages, definitions, and axiomatizations as \textsl{precise} in the sense that he regarded them as necessary for a proper understanding of their first-order counterparts.   This application of our prior terminology also clarifies Kreisel's initial remark that  `to explain the notion of mathematically definite problem one needs the corresponding notion of second order consequence'.  In this way, the role played by the notions  \textsl{mathematical definiteness} in Kreisel's CH argument can be analogized to that occupied by the common notion of \textsl{validity} in the argument of \S\ref{validity}.  For while Kreisel proposes that we should investigate the restriction of these concepts to first-order statements, he suggests that their analysis may (at least \textsl{prima facie}) require higher-order notions.   

Definition \ref{md} can accordingly be understood as providing a precise analysis of such a notion of mathematical definiteness in the same manner that Kreisel argued in the course of his validity argument that  $V$ -- i.e. truth in all set-sized models -- provides a precise analysis of $\mathit{Val}$.   But at this stage, two additional questions arise: 
\begin{enumerate}[i)]
\item Is it reasonable to regard mathematical definiteness as a genuinely common notion?  
\item If so what sort of argument can be given that Definition \ref{md} provides a \textsl{correct} analysis of mathematical definiteness?
\end{enumerate}

Kreisel does not speak directly to i).  But the notion of definiteness in question seems unlike several of the other concepts which he thought were amenable to informal rigour.  In particular, reflection on mathematical definiteness does not appear to play as clear a role in our prior mathematical practice in the manner Kreisel suggests that the notion of validity (or logical consequence) figured systematically in the evaluation of mathematical arguments long before the model theoretic definition was given (although see note \ref{geonote} below).  But in regard to ii), one might still hope that our grasp of this notion determines it sufficiently that an argument for the adequacy of Definition \ref{md} could be mounted via the method of squeezing.   Suppose, for instance, we introduce the $\mathcal{L}_C$-predicate $\mathit{Def}$ to denote the relevant common concept of definiteness and let $\pi_w(x)$ abbreviate the precise notion defined by the formalization of Definition \ref{md} in an appropriate precise language $\mathcal{L}_P$.    We might then hope that it is possible to give an informally rigorous argument which shows that the extension of $\mathit{Def}$ is `squeezed' between $\pi_w(x)$ and some appropriate narrow precise notion $\pi_n(x)$.   

For reasons we will discuss in \S\ref{chsol}, it will not suffice to mimic the validity argument by taking $\pi_n(x)$ to correspond to provability or refutability in a proof system for second-order logic with a computably enumerable derivability relation.   Not was Kreisel likely ware aware of this, subsequent results illuminate why we should expect that the extension of $\mathit{Def}$ is sufficiently complex that no such argument will be forthcoming.  The fact that Kreisel \textsl{fails} to provide one himself thus represents a notable difference between his CH and validity arguments.  What he does instead is to argue in favor of the extensional correctness of Definition \ref{md} on the basis of several examples in addition to CH itself.  Of these, we will now examine two of the most illustrative.

It is a consequence of Definition \ref{md} that various first-order consequences of the Axiom of Replacement are \textsl{not} mathematically definite in Kreisel's sense.   For instance let $\rho_{\alpha}$ abbreviate the first-order statement formalizing that `$R(\alpha)$ exists' (where $\alpha$ is itself an ordinal definable in $\mathcal{L}^1_{\mathsf{Z}}$).  Then $\rho_{\omega + \omega}$ is not mathematically definite in the sense of Definition \ref{md} since -- e.g. $R(\omega + \omega) \models \mathsf{Z}^2 + \neg \rho_{\omega + \omega}$ but $R(\omega + \omega + \omega) \models \mathsf{Z}^2 + \rho_{\omega + \omega}$ (and of course also $\mathsf{ZF}^2 \proves \rho_{\omega + \omega}$ even though $\mathsf{Z}^2 \not\proves \rho_{\omega + \omega}$).  But of course one might also take the truth of $\rho_{\omega + \omega}$ to follow from the so-called \textsl{iterative conception of set} which we will see in \S\ref{set} Kreisel takes to motivate the axioms of $\mathsf{Z}^2$ itself.   For his own part, however, he stressed that the \textsl{non-definiteness} of such statements illustrates that his definition is non-trivial in the sense that it does apply to \textsl{all} $\mathcal{L}^1_{\mathsf{Z}}$-statements.\footnote{Kreisel goes on to remark `Reflection shows that the logical undecidability results which surprise mathematicians concern mathematically definite problems like the continuum hypothesis, not the existence of [$R(\omega + \omega)$]' \citeyearpar[p. 257]{Kreisel1967a} (see also \citealp[pp. 98-99]{Lakatos1967}).   His apparent point is thus that although opinions may vary on the definiteness of `axioms of infinity' even as weak as $\rho_{\omega + \omega}$, CH concerns the structure of $R(\omega + 2)$.   For as this structure is presumably of more concern to mathematical practice, there is reason to be more circumspect about the potential indefiniteness of CH on the basis of its formal independence results than we might be about the consequences of Replacement in regard to the existence of yet larger ordinals.}       

Kreisel also suggests that another famous example of formal independence -- i.e. the Parallel Postulate of geometry [PP] -- is more akin to Replacement (or $\rho_{\omega + \omega}$) than it is to CH.   For as he observes, PP remains independent of Hilbert's \citeyearpar{Hilbert1899} basic axiomatization of geometry even when the various continuity principles he considers -- e.g. the Archimedean Axiom or the existence of Dedekind cuts -- are formalized as second-order axioms rather than first-order schemas.   It thus follows that if the theory $\mathsf{T}^2$ in Definitions \ref{cddefn} and \ref{md} were taken to be an appropriate axiomatization of geometry,  then PP would also be classified as \textsl{non-definite}.    But not only does this differ from the classification of CH,  Kreisel appears to take this to reflect the manner in which the status of PP was ultimately resolved within the practice of geometry.\footnote{Although Kreisel makes this point only in passing in \citeyearpar[p. 151]{Kreisel1967b}, he spells it out in more detail in \citeyearpar[pp. 109-110]{Kreisel1969a}.   One might indeed take the long history of attempts to assess whether PP is determinately true or false relative to various conceptions of space to counter our prior claim that reflection on definiteness has not played a substantial role in prior mathematical practice.   But not only were the original constructions of Beltrami, Klein, and Poincar\'e of non-Euclidean geometries contentious when they were originally proposed, many of the developments which followed -- perhaps most famously Hilbert's \citeyearpar{Hilbert1899} consistency proofs using analytical models and his subsequent debate with Frege about what they showed -- appear to testify to the apparent \textsl{lack} of consensus about the definiteness of specific geometrical statements.  But again for his own part, Kreisel hoped to deploy the distinction between PP and CH relative to his notion of definiteness to illustrate a problematic feature of the analogy between geometry and set theory which \citet{Cohen1967} and \citet{Robinson1968} had recently exploited to promote a formalist understanding of set theoretic independence results. \label{geonote}}

Even if it is agreed that mathematical definiteness is a common notion, one's intuitions about what such case studies show -- and thus also of the aptness of Kreisel's definition -- may still differ.  But if these difficulties set aside, we can at least attempt to evaluate whether the argument for the definiteness of CH rehearsed above can be formalized as a `philosophical theorem' in the manner which we have suggested is possible in the case of Kreisel's validity and creating subject arguments.   

\subsubsection{From schematization to formalization}
\label{chsf}

We have already taken the initial step of proposing that $\mathit{Def}$ should be regarded as a predicate in the common language $\mathcal{L}_C$ and that we should understand Definition \ref{md} as defining a predicate $\pi(x)$ formulated in the precise language $\mathcal{L}_P$.  It is evident that Kreisel understood the considerations just adduced to establish that the principle
\begin{example}
$\A \phi^1(\mathit{Def}(\phi^1) \leftrightarrow \pi(\phi^1))$
\label{defe}
\end{example}
should be taken to be a member of the set $\Gamma_2$ of bridging principles stated in the joint (i.e. common plus precise) language $\mathcal{L}_J$.   The next question is that of specifying the precise language $\mathcal{L}_P$ so that Definitions \ref{cddefn} and \ref{md} can be formulated in $\mathcal{L}_J$.   It is at this stage where Kreisel's decision to limit the domain of application of his definition of definiteness to first-order formulas becomes significant.   For if this restriction is imposed, then there is no in principle obstacle to formalizing the semantic notions which appear in these definitions so that Corollary \ref{kch} can be regarded as a genuine theorem in an appropriate choice for the Kreiselian theory $\mathsf{T}_K$.  

We will take $\mathcal{L}_P$ to correspond to the language $\mathcal{L}^2_{\mathsf{Z}}$ of second-order set theory and the mathematical theory $\mathsf{T}_P$ on which $\mathsf{T}_K$ is based to be $\mathsf{ZF}^2$ itself.  For in this case, it is straightforward to formalize Kreisel's definition of a \textsl{Zermelo system} \citeyearpar[p. 255-256]{Kreisel1967a} to obtain a single $\mathcal{L}^2_{\mathsf{Z}}$ formula $\zeta^2(X,Y)$ which expresses that the axioms of $\mathsf{Z}^2$ hold of the structure $\langle X,Y \rangle$ where $Y \subseteq X \times X$.   Using standard techniques from the arithmetization of syntax, it is simlarly possible to take the predicate $\mathrm{Sat}^k_1(X,Y_1,\ldots,Y_k,\ulcorner \varphi^1 \urcorner)$ to abbreviate an $\mathcal{L}^2_{\mathsf{Z}}$-formula defining satisfaction for $\mathcal{L}^1_{\mathsf{Z}}$-formulas of the sort which Kreisel describes \citeyearpar[p. 155]{Kreisel1967b}.\footnote{We have also followed Kreisel's presentation in  by suppressing complications arising from the need to formulate substitution for free variables in $\phi^1$ which would be required in a proper definition of a formal satisfaction predicate.  See, e.g., \citep[\S 3.5]{Drake1974} for the relevant details.}    Once these steps are undertaken, a precise definition of the definiteness predicate $\pi(x)$ can be stated in $\mathcal{L}^2_{\mathsf{Z}}$ as follows:\\

\noindent $\pi(\ulcorner \phi^1 \urcorner) =  \A X \A Y \subseteq X^2(\zeta^2(X,Y) \rightarrow \mathrm{Sat}^1_1(X,Y,\ulcorner \phi^1 \urcorner)) \mor \A X \A Y \subseteq X^2(\zeta^2(X,Y) \rightarrow \mathrm{Sat}^1_1(X,Y,\ulcorner \neg \phi^1 \urcorner))$\\

If we now take $\mathit{Def}'$ as an arithmetized stand-in for the common predicate $\mathit{Def}$ applicable to codes of sentences, then Kreisel's informal argument for (\ref{defe}) yields
\begin{example}
$\A x (\mathit{Def}'(x) \leftrightarrow \pi(x))$
\label{defp}
\end{example}

Working within $\mathsf{ZF}^2$ it is also possible to show that if an ordinal $\alpha$ exists then so does the structure $\alpha$th level of the cumulative hierarchy -- i.e. the structure $\langle R(\alpha),\in \upharpoonright R(\alpha) \times R(\alpha) \rangle$.  We will employ $R_{\alpha}$ and $E_{\alpha}$ as abbreviations for terms defining the components of this structure.  We will also use $\mathrm{Iso}(X,Y,U,V)$ to abbreviate an $\mathcal{L}_{\mathsf{Z}}^2$-predicate expressing that $X \subseteq Y^2$ and $U \subseteq V^2$ and there exists an isomorphism between the structures $\langle X,Y \rangle$ and  $\langle U,V \rangle$ -- i.e. a bijection  $F:X \rightarrow Y$ such that $U(x,y)$ if and only if $V(F(x),F(y))$.  

With these conventions in place, it is then straightforward to formalize the proof of Theorem \ref{zqct} to obtain:\footnote{Here $\proves_2$ again denotes derivability in the standard deductive system for second order logic -- see, e.g., \citep[\S 3.2]{Shapiro1991} or \citep[\S 5]{Dalen2008}.   Such a formalization can be compared to the so-called \textsl{internal quasi-catgeoricity} results for theories similar to $\mathsf{ZF}^2$ reported by \citet{Vaananen2015} (Theorem 3), \citet{Button2018} (Corrolary 11.3), and \citet{Vaananen2019} (Theorem 1) to which we will return below.} 
\begin{example}
$\mathsf{ZF}^2 \proves_2 \A X \A Y \subseteq X^2[\zeta^2(X,Y) \rightarrow \exists \lambda \geq (\omega + \omega)(\mathrm{Iso}(X,Y,R_{\lambda},E_{\lambda})]$
\label{intcat}
\end{example}
Next note that it is also possible to formalize in $\mathsf{ZF}^2$ the observation that $\omega+\omega$ rank-bounded formulas such as $\chi^1$ are absolute between $R(\omega + \omega)$ and $R(\alpha)$ for $\alpha \geq \omega + \omega$.  Together with our choice of $\chi^1$ this yields
\begin{example}
$\mathsf{ZF}^2 \proves_2 \A \alpha \geq (\omega + \omega)[\zeta^2(X,Y) \rightarrow  (\mathrm{Sat}^1_1(R_{\omega + \omega},E_{\omega + \omega},\ulcorner \chi^1 \urcorner) \leftrightarrow \mathrm{Sat}^1_1(R_{\alpha},E_{\alpha},\ulcorner \chi^1 \urcorner) ]$
\label{abs}
\end{example}
Finally note that within $\mathsf{ZF}^2$ it is also possible to formalize the conventional argument that isomorphism implies elementary equivalence for first-order sentences -- i.e. 
\begin{example}
$\mathsf{ZF}^2 \proves_2 \A X \A Y \A U \A V \A \ulcorner \phi^1 \urcorner [\mathrm{Iso}(X,Y,U,V) \imp \A \ulcorner \phi^1 \urcorner(\mathrm{Sat}^1_1(X,Y,\ulcorner \phi^1 \urcorner) \leftrightarrow \mathrm{Sat}^1_1(U,V,\ulcorner \phi^1 \urcorner))]$
\label{isoelem}
\end{example}

Relative to the definitions we have adopted, the mathematical definiteness of CH is expressed by $\mathrm{Def}'(\ulcorner \chi^1 \urcorner)$ which by (\ref{defp}) is asserted to be equivalent to the (precise) $\mathcal{L}^2_{\mathsf{Z}}$-statement $\pi(\ulcorner \chi^1 \urcorner)$.   This can be derived by the following expected argument:
\begin{example}
\begin{enumerate}[i)]
\item Either a) $\mathrm{Sat}^1_1(R_{\omega + \omega},E_{\omega + \omega},\ulcorner \chi^1 \urcorner)$ or b) $\neg \mathrm{Sat}^1_1(R_{\omega + \omega},E_{\omega + \omega},\ulcorner \chi^1 \urcorner)$.
\item Assuming a), let $X,Y$ be arbitrary and assume $\zeta^2(X,Y)$.
\item Then by (\ref{intcat}), $\exists \lambda \geq (\omega + \omega) \mathrm{Iso}(X,Y,R_{\lambda},E_{\lambda})$.
\item Fixing $\lambda$, it then follows by a) and (\ref{abs}), $\mathrm{Sat}^1_1(R_{\lambda},E_{\lambda},\ulcorner \chi^1 \urcorner)$.
\item But then $\mathrm{Sat}^1_1(X,Y,\ulcorner \chi^1 \urcorner)$ by ii) and (\ref{isoelem}).  
\item Discharging assumption ii) and generalizing, we obtain $\A X \A Y \subseteq X^2(\zeta^2(X,Y) \rightarrow \mathrm{Sat}^1_1(X,Y,\ulcorner \chi^1 \urcorner))$.
\item Assuming b),  and reasoning in parallel to iii) - v) we obtain $\neg \mathrm{Sat}^1_1(X,Y,\ulcorner \chi^1 \urcorner)$. It then follows from the properties of satisfaction predicate that $\mathrm{Sat}^1_1(X,Y,\ulcorner \neg \chi^1 \urcorner)$.
\item Discharging and generalizing in parallel to vi), we obtain $\A X \A Y \subseteq X^2(\zeta^2(X,Y) \rightarrow \mathrm{Sat}^1_1(X,Y,\ulcorner \neg \chi^1 \urcorner)$.
\item Reasoning from i), vi), and viii) by constructive dilemma we can conclude 
$$\A X \A Y \subseteq X^2(\zeta^2(X,Y) \rightarrow \mathrm{Sat}^1_1(X,Y,\ulcorner \chi^1 \urcorner)) \mor \A X \A Y \subseteq X^2(\zeta^2(X,Y) \rightarrow \mathrm{Sat}^1_1(X,Y,\ulcorner \neg \chi^1 \urcorner))$$
\end{enumerate}
\label{formcharg}
\end{example}

It is evident that the forgoing argument can be formalized to obtain the following:\footnote{Kreisel does not state this result explicitly in \citeyearpar{Kreisel1967a} or  \citeyearpar{Kreisel1967b}.  But in his review \citeyearpar{Kreisel1977b} of \citep{Weston1976} he remarks that \ref{kwt} should be regarded as a `formal theorem'.}
\begin{theorem}
$\mathsf{ZF}^2 \proves_2 \pi(\ulcorner \chi^1 \urcorner)$.
\label{kwt}
\end{theorem}
\noindent It thus follows that if we take the relevant Kreiselian theory to be $\mathsf{T}_K = \mathsf{ZF}^2 \cup \Gamma_2$ -- where we agreed that the latter set includes (\ref{defp}) -- then $\mathsf{T}_K \proves \mathit{Def}'(\ulcorner \chi^1 \urcorner)$ as desired. This makes good on the claim that Kreisel's CH argument can be assimilated to the schema IR by which we have argued in \S\ref{oninfrig} that informally rigourous arguments may themselves be rigorized in the manner of `philosophical proofs'.   This methodological point aside, it will also be useful to examine some additional features of the argument both in its original context and in regard to subsequent developments.

\subsubsection{Novel set theoretic notions?}
\label{chnovel}

One outstanding issue pertains to Kreisel's evident hope that \textsl{novel} concepts might be invoked to extend the foregoing argument in a manner which would settle the truth value of CH on the basis of the intended interpretation of the language of set theory in the manner he viewed his creating subject argument as settling the truth value of GMP relative to the intended interpretation of the language of intuitionistic analysis.  He remarks on this at the end of the first section of \citeyearpar{Kreisel1967b}:

\begin{quote}
{\footnotesize Finally, and this is of course the most direct link between the present section and the main theme of this article, [the] second order decidability of CH  suggests this: new primitive notions, e.g. properties of natural numbers, which are \textsl{not} definable in the language of set theory (such as in the footnote on p. 150), may have to be taken seriously to decide CH; for, what is left out when one replaces the second order axiom by the schema, are precisely the properties which are not so definable. \citeyearpar[p. 152]{Kreisel1967b}}
\end{quote}

The specifics of Kreisel's proposal are clarified by the footnote to which he alludes together with the immediately preceding text:
\begin{quote}
{\footnotesize  One expects [that second-order consequence may be expressed by a first-order formula of set theory] simply because it is always claimed that this first order language is adequate for all mathematics; so if it weren't adequate for expressing second order consequence, somebody would have noticed. [Footnote: One cannot be 100 per cent sure: for instance, consider the so called truth definition. We have here a set $T$ of natural numbers, namely G\"odel numbers $\ulcorner \alpha_i\urcorner$  of first order formulae of set theory, such that $n \in T \leftrightarrow \E i(n = \ulcorner \alpha_i\urcorner \ \& \  \alpha_i) \ldots$ As Tarski emphasised, $T$ is not definable by means of a first order formula (in the precise sense above).] \hfill  \citeyearpar[p. 150]{Kreisel1967b}}
\end{quote}

This makes clear that the particular novel concept that Kreisel had in mind for settling CH was thus that of \textsl{set theoretic truth}.   This is programmatically significant in the sense that we have just illustrated how a formal satisfaction predicate may be used to formalize his CH argument.  As we have seen in \S\ref{validity}, \citet[pp. 155-158]{Kreisel1967b} also used such a predicate to illustrate the relationship between his validity predicate and set theoretic truth in a manner which requires that $\mathrm{Sat}^1_1(X,Y,\ulcorner \phi \urcorner)$ provably satisfied Tarski's T-biconditional for all first-order formulas -- i.e. $\mathsf{ZF}^2 \proves \mathrm{Sat}^1_1(X,Y,\ulcorner \phi^1 \urcorner) \leftrightarrow \phi^1$ for all $\phi^1 \in \mathcal{L}^1_{\mathsf{Z}}$.    But in virtue of this, it follows that this predicate cannot be defined by a formula of $\mathcal{L}^1_{\mathsf{Z}}$ in virtue of Tarski's undefinability theorem.  In this sense, Kreisel is correct that the concept of first-order truth (or satisfaction) is a \textsl{novel} concept if all that is meant it is that it is  not first-order definable.    

Although Kreisel does not develop this point further in \citeyearpar{Kreisel1967b}, he later observed

\begin{quote}
\footnotesize{[W]e know that there are lots of sets, even of natural numbers, which cannot be defined in $\mathcal{L}_E$. More specifically, we have Tarski's implicit definition of the satisfaction relation which cannot be explicitly defined in $\mathcal{L}_E$. If we add the defining proposition to set theory and expand the axiom schemata, we get new theorems formulated in $\mathcal{L}_E$ itself; and also, for instance, a very natural consistency proof for certain reflection principles. Note that \textsl{the traditional `reduction' of mathematics $($arithmetic, continuum$)$ to set theory does not cover these cases}. \hfill \citeyearpar[p. 100]{Kreisel1969a}}
\end{quote}

It is not made clear here whether the language at issue here is intended to be first- or second-order.   But this passage is at least suggestive of the proposal that set theoretic truth (or satisfaction) be treated as a new predicate which can be adjoined  to a theory such as $\mathsf{ZF}^1$ (or perhaps $\mathsf{ZF}^2$) along with axioms formalizing the clauses in the familiar inductive definition of truth, now presumably understood as constitutive principles for a novel concept.  But in this regard, it is also significant that immediately following both of the preceding passages Kreisel cites G\"odel's address \citeyearpar{Godel1946} as evidence that not only are `axioms of infinity $\ldots$ more efficient' than truth-theoretic principles for obtaining the sorts of consequences he has in mind \citeyearpar[p. 152]{Kreisel1967b} but also that the former can be used to `replace' the latter \citeyearpar[p. 100]{Kreisel1969a}.  

In neither instance are Kreisel's remarks sufficiently detailed to be sure of a precise reconstruction.   But consideration of both contemporaneous and subsequent results suggest that adjoining `axioms of infinity' -- i.e. \textsl{large cardinal hypotheses} as we would now call them -- does indeed typically lead to stronger extensions than are obtained by adding the sort of truth-theoretic principles Kreisel apparently had in mind.\footnote{As Parsons suggests in his introduction to \citep[p. 146]{Godel1946} the sort of `replacement' of a primitive notion of satisfaction which Kreisel has in mind is likely to already be effected by the Montague-L\'evy first-order reflection scheme for $\mathsf{ZF}$.  In fact \citet[Theorem 9]{Kreisel1968a} would themselves go on to show that $\mathsf{ZF}$ proves all instances of the proof-theoretic reflection scheme for $\mathsf{Z}$ and thus also the novel arithmetical consequence $\mathrm{Con}(\mathsf{Z})$.  In a similar spirit, the second-order set-theoretic reflection scheme previously introduced by \citet{Levy1960} -- which implies the existence of Mahlo cardinals -- can be used to prove $\mathrm{Con}(\mathsf{ZF})$ and even stronger statements like $\mathrm{Con}(\mathsf{ZF}) + \textrm{`there exists a strongly inaccessible cardinal'}$.  But since Kreisel would presumably have regarded such principles as \textsl{non-definite} in his technical sense it is unclear exactly what he should have made of such results.  On the other hand, the more recent exploration of truth-theoretic extensions of set theory has revealed that systems which mimic the familiar extensions of $\mathsf{PA}$ with a primitive truth predicate satisfying the clauses in Tarski's definition of truth (e.g. $\mathsf{CT}, \mathsf{KF}$, etc.) do not even entail the existence of inaccessible cardinals (see, e.g. \citealp{Fujimoto2012}).}   But on the other hand, what came to be known as \textsl{G\"odel's program} for deciding statements formally independent of $\mathsf{ZF}$ via large cardinal hypotheses is now widely understood to be incapable of determining the truth value of CH itself.\footnote{See, e.g., \citep{Steel2014} for an overview.}   Such limitative results do not rule out the possibility that consideration of some other novel (or simply non-first-order expressible) concepts might someday lead to what might come to be regarded as an informally rigorous argument for either CH or $\neg \mathrm{CH}$.   But they do appear to deflate the hope that reflection on the specific concept of \textsl{set theoretic truth} can settle CH in the manner which Kreisel hoped it might in the immediate wake of the discovery of its formal independence.

\subsubsection{The role of internal categoricity}
\label{chic}

As should now be evident, a mathematical result -- i.e. Zermelo's Theorem \ref{zqct} -- functions as the crux of Kreisel's CH argument in much the same way that G\"odel's Completeness Theorem and Kreisel's own refutation of generalized Markov's Principle are respectively central in his validity and creating subject arguments.   As we have seen, this is one of several traditional categoricity theorems for second-order axiomatizations to which Kreisel draws attention.    But although these results were originally stated and proven informally,  we have also noted that there is no technical obstacle to formalizing the relevant case of Theorem \ref{zqct} so that the object language expression of the quasi-categoricity of $\mathsf{Z}^2$ becomes a formal theorem of a suitable background theory.  But once these steps are initiated, some additional questions about the CH argument come into focus in regard to recent discussions of the phenomenon known as \textsl{internal categoricity}.\footnote{Although the expression `internal catgeoricity' appears to have been introduced in \citet{Walmsley2002}, the basic definitions and results were formulated earlier by \citet{Parsons1990a} and \citet{McGee1997}.    A general theory which encompasses the examples considered in these sources has been developed more systematically by V\"a\"an\"anen \citeyearpar{Vaananen2012,Vaananen2015,Vaananen2020}.}

As formulated by \citet{Vaananen2012} internal categoricity should, in the first instance, be understood semantically relative to the so-called \textsl{Henkin semantics} for second-order logic.\footnote{See, e.g.,  \citep[\S 4.3]{Shapiro1991}.  Recall in particular that a \textsl{Henkin model} of second-order logic is a structure $\langle \mathfrak{M}, \mathfrak{S} \rangle$ where $\mathfrak{M}$ is a model for a given signature in the usual first-order sense and $\mathfrak{S}$ is a collection of subsets of the domain of $\mathfrak{M}$ of appropriate arities which is sufficiently rich to satisfy the comprehension scheme when second-order quantifiers are restricted to $\mathfrak{S}$ but may be a proper subclass of the full powerset of its domain.}  \citet[p. 122]{Vaananen2015} then define a second-order theory $\mathsf{T}^2$ to be internally categorical just in case `all models $\mathsf{T}^2$ within a common Henkin model are witnessed to be isomorphic by the model'.  Recall, however, that \citet{Henkin1950} originally showed that the Henkin semantics leads to a definition of second-order validity $\models^h_2$ which is complete with respect to the standard definition of $\proves_2$.   Adapting our notation from above, the internal categoricity of $\mathsf{T}^2$ can thus also be defined proof-theoretically by the condition
\begin{example}
$\proves_2 \A X_1 \A Y_1 \ldots Y_K \A X' \A Y'_1 \ldots \A Y'_k[(\tau^2(X_1,Y_1,\ldots, Y_k) \mand \tau^2(X'_1,Y'_1,\ldots, Y'_k)) \rightarrow \mathrm{Iso}(X_1,Y_1,\ldots, Y_k,X'_1,Y'_1,\ldots, Y'_k)]$
\end{example}
In other words, the internal categoricity of a theory is equivalent to the formal derivability in second-order logic of the existence of an isomorphism between any structures which satisfy its axioms.   

We have reconstructed Kreisel's CH argument using $\mathsf{ZF}^2$ itself as the relevant `precise' mathematical theory $\mathsf{T}_P$.  However as was shown by \citep{Vaananen2015} (and is reproved in more details by \citealp[\S 11.C]{Button2018}) $\mathsf{ZF}^2$ is internally categorical subject to an additional assumption which ensures that all structures satisfying its axioms have the same `height'.\footnote{E.g. that there is an isomorphism between their ordinals or that there are no inaccessible cardinals $> \omega$.}   These treatments also show that as long as the full comprehension axiom in the language $\mathcal{L}^2_{\mathsf{Z}}$ is subsumed under the definition of $\proves_2$, then the derivation of (\ref{intcat}) also goes through in pure second-order logic. In the version of the argument we have presented above, some mathematical principles are still required to handle the reasoning about the satisfaction predicate.   But as these too can be carried out in a fragment of $\mathsf{ZF}^2$ this suggests that it should be possible to formalize the argument relative to a weaker choice for $\mathsf{T}_P$.\footnote{In fact \citet[\S 12]{Button2018} show in a somewhat different setting that a generalization of the argument (\ref{formcharg}) leading to the definiteness of all $n^{\mathrm{th}}$-order set-theoretic statements can be carried out in pure $n+3^{\mathrm{rd}}$-order logic.  However this degree of generality is not needed here due to Kreisel's decision to restrict his definiteness predicate to formulas in the language of first-order set theory.}

But the fact that the statement $\pi(\ulcorner \chi^1 \urcorner)$ which we have suggested embodies Kreisel's analysis of the definiteness of CH is a formal theorem of $\mathsf{ZF}^2$ already raises a more fundamental question about what the CH argument ultimately should be understood to show.  For recall that \citet{Weston1977} (building on  \citealp{Chuaqui1972}) demonstrated that the formal independence of CH from first-order $\mathsf{ZF}$ extends to $\mathsf{ZF}^2$.   It thus follows that there exist Henkin models $\langle \mathfrak{M}_1, \mathfrak{S}_1 \rangle$ and $\langle \mathfrak{M}_2, \mathfrak{S}_2 \rangle$ such that
\begin{example}
\begin{enumerate}[i)]
\item $\langle \mathfrak{M}_1, \mathfrak{S}_1 \rangle \models^h_2 \A X \A Y \subseteq X^2 (\zeta^+(X,Y) \imp \mathrm{Sat}^1_1(X,Y,\ulcorner \chi^1 \urcorner))$
\item $\langle \mathfrak{M}_2, \mathfrak{S}_2 \rangle \models^h_2 \A X \A Y \subseteq X^2 (\zeta^+(X,Y) \imp \mathrm{Sat}^1_1(X,Y,\ulcorner \neg \chi^1 \urcorner))$
\end{enumerate}
\end{example}
where $\zeta^+(X,Y)$ formalizes that $\langle X,Y \rangle$ satisfies $\mathsf{ZF}^2$.    We thus reach the conclusion that there are Henkin models of $\mathsf{ZF}^2$ in which all structures which satisfy $\mathsf{ZF}^2$ validate CH and also Henkin models in which all structures which satisfy $\mathsf{ZF}^2$ falsify CH.   

This may at first seem like an incongruous situation.   But as \citet{Vaananen2015} observe, all the existence of such models shows is that $\mathfrak{S}_1$ and $\mathfrak{S}_2$ cannot be definably brought together to form a joint Henkin model of $\mathsf{ZF}^2$.   On the other hand,  the fact that the argument (\ref{formcharg}) can be formalized in $\mathsf{ZF}^2$ shows that both $\langle \mathfrak{M}_1, \mathfrak{S}_1 \rangle$ and$ \langle \mathfrak{M}_2, \mathfrak{S}_2 \rangle$ must also satisfy $\pi(\ulcorner \chi^1 \urcorner)$.   It is thus a consequence of the derivability of $\pi(\ulcorner \chi^1 \urcorner)$ that its truth in a given Henkin model $\langle \mathfrak{M}, \mathfrak{S} \rangle$  cannot be diagnostic of the fact that CH assumes the same truth value across all such models but only of those which are `internal' to $\langle \mathfrak{M}, \mathfrak{S} \rangle$ itself.

As we have seen, however, Kreisel took  the statement $\pi(\ulcorner \chi^1 \urcorner)$ to express the mathematical definiteness of CH.   Relative to the model of informal rigour we have presented here, it is thus indeed significant that this statement is formally derivable in an appropriate choice for the Kreiselian theory $\mathsf{T}_K$.  But note in this case $\mathsf{T}_K$ extends what we have taken to be the precise theory  $\mathsf{T}_P$ (i.e. $\ \mathsf{ZF}^2$) by adding only the bridging principle (\ref{defp}) whose role is simply that of arithmetizing (\ref{defe}).  As in this case $\mathsf{T}_K$ is an extension-by-definitions of $\mathsf{T}_P$, the former is evidently \textsl{conservative} over the latter.    In this sense, it would thus indeed be unreasonable to expect that this formulation of Kreisel's argument provides a stronger guarantee of the definiteness of CH than is expressed by any statement derivable in $\mathsf{ZF}^2$.\footnote{One could foresee Kreisel objecting that the formalization of Definition (\ref{md}) which we have proposed conflates with informal and informal rigour (although see \citeyear[p. 257-262]{Kreisel1967a} and note \ref{formnote} below). But this point aside, the conservativity of $\mathsf{T}_K$ over $\mathsf{T}_P$ still marks a significant formal contrast between his creating subject and CH arguments.}

The prior sequence of observations also highlight a feature of Kreisel's CH argument which  \citet{Weston1974,Weston1976} suggests detracts from its significance.  For as we have just seen, in order for the truth or falsity of $\pi(\ulcorner \chi^1 \urcorner)$ to be genuinely diagnostic of the definiteness of CH requires that we restrict attention to full models of $\mathsf{Z}^2$.   Zermelo's Theorem entails that all such models of $\mathsf{Z}^2$ will have an initial segment isomorphic to $R(\omega + \omega)$.   Since the truth value of $\chi^1$ in this structure will thus determine its truth value in all other full models, Weston suggests that Kreisel's argument for the definiteness of CH reduces to its truth or falsity in this single structure viewed as a first-order model.   But since (trivially) CH is either determinately true or determinately false in this model\footnote{Note that this is simply the `dilemma' which is expressed as premise (\ref{formcharg}i) in our prior formalization of the validity argument.} -- i.e. $R(\omega + \omega) \models_1 \chi^1 \mor \neg \chi^1$ -- Weston concludes that Kreisel's detour through second-order logic does not provide us further asurance of the definiteness of CH beyond whatever evidence can be mounted for the existence of a unique intended model of $\mathsf{ZF}$ (or even $\mathsf{Z}$).   

Kreisel does not directly react to Weston's point in \citeyearpar{Kreisel1967b}.   But in both \citeyearpar[p. 257-262]{Kreisel1967a} and \citeyearpar[p. 192-194]{Kreisel1967c}, he anticipates many of the observations which lead up to it.   One of the central morals he derived was aimed at exposing what he took to be a rhetorical instability in the position of his formalist interlocutors: 

{\footnotesize 
\begin{quote}
What is suspect is the \textsl{significance} of [model theoretic] formal independence proofs for someone who $\ldots$ in the same breath uses model theoretic methods. For $\ldots$ when doing so he thinks in terms of a notion of set which makes the formally undecided problem mathematically definite. So there is certainly an informal contradiction between the basic \textsl{importance} of formal independence and the \textsl{acceptance} of the semantic interpretation.  \hfill \citeyearpar[p. 260]{Kreisel1967a} 
\end{quote}
}
\noindent This in turn brings us back to a question which we left open in \S\ref{validity}: how should we understand the relationship between Kreisel's CH argument and his validity argument?

\subsubsection{Squeezing, set theory, and second-order logic}
\label{chsol}

We are now in a position to appreciate that there is indeed at least a \textsl{prima facie} tension between Kreisel's validity and CH arguments:  while the former is intended to highlight the significance of the precise definability of the notion of validity for first-order languages, the latter is premised on the superiority of second-order languages for characterizing mathematical structures and concepts.    As such, a natural question is whether the method of the validity argument can be extended to provide a precise analysis of the intuitive notion of second-order validity.

We have already touched on this question in regard to the hope of providing a squeezing-like argument to support Kreisel's definition of mathematical definiteness.   But in fact Kreisel raises the issue directly himself:
\begin{quote}
\footnotesize{All this [the validity argument] was for first order formulae. \textsl{For higher order formulae we do not have a convincing proof} of e.g. $\A \alpha^2(V(\alpha^2) \leftrightarrow Val(\alpha^2))$ although one would expect one.\\ \hspace*{1ex} \hfill \citeyearpar[p. 157]{Kreisel1967b}}
\end{quote}
Kreisel does not speak further to this question in either \citeyearpar{Kreisel1967a} or \citeyearpar{Kreisel1967b}.   On the other hand, there can be no doubt that he was aware of Henkin's \citeyearpar{Henkin1950} proof of the deductive completeness of second-order logic with respect to the Henkin semantics.\footnote{See, e.g., \citet[p. 120]{Kreisel1952c}.}    The fact that he \textsl{fails} to propose an argument by which $\mathit{Val}(\alpha^2)$ -- i.e. the intuitive validity of the second-order statement $\alpha^2$ -- can be squeezed between $V(\alpha^2)$ -- i.e. the validity of $\alpha^2$ with respect to set-sized structures -- and $D(\alpha^2)$ -- i.e. the derivability of $\alpha^2$ relative to the definition of $\proves_2$ -- is thus already significant.   

Given the connection which Kreisel emphasized between formal completeness and informal rigour in the introduction to \citeyearpar{Kreisel1967b}, it seems unlikely that this was simply an oversight.  But at least at the programmatic level, it is easy to see why he is likely to have opposed an interpretation of $\mathit{Val}(\alpha^2)$ which would permit a straighforward generalization of his validity argument.    For building on our discussion of $\mathit{Val}(\alpha^1)$ and $V(\alpha^1)$ in \S\ref{valmeaning}, we can now see that there are in fact two ways in which we might elect to characterize \textsl{both} of these notions in the second-order case:
\begin{example}
\begin{enumerate}[i)]
\item $\mathit{Val}_{h}(\alpha_2)$ iff $\alpha_2$ is true in all (potentially class-sized) structures with respect to the Henkin semantics.
\item $\mathit{Val}_{s}(\alpha_2)$ iff $\alpha_2$ is true in all (potentially class-sized) structures with respect to the standard semantics.
\item $\mathit{V}_{h}(\alpha_2)$ iff $\alpha_2$ is true in all (set-sized) models with respect to the Henkin semantics.
\item $\mathit{V}_{s}(\alpha_2)$ iff $\alpha_2$ is true in all (set-sized) models with respect to the standard semantics.
\end{enumerate}
\label{valopts}
\end{example}

It is presumably uncontentious that the standard definition of $\proves_2$ is intuitively sound with respect to both refinements of $\mathit{Val}(\alpha^2)$ --  i.e.
\begin{example}
\begin{enumerate}[i)]
\item $D(\alpha^2) \rightarrow \mathit{Val}_s(\alpha^2)$
\item $D(\alpha^2) \rightarrow \mathit{Val}_h(\alpha^2)$
\end{enumerate}
\end{example}
If we were to additionally opt to understand informal and formal second-order validity respectively via (\ref{valopts}i,iii), then we would presumably also be in a position to accept
\begin{example}
$Val_h(\alpha^2) \rightarrow V_h(\alpha)$
\end{example}
And of course Henkin's completeness proof also yields
\begin{example}
$V_h(\alpha^2) \rightarrow D(\alpha^2)$
\end{example}
In this case we would indeed be able to conclude that $\forall \alpha^2(\mathit{Val}(\alpha^2) \leftrightarrow V_h(\alpha^2))$.  

But as we are now in a position to appreciate, it would seem that Kreisel is likely to have insisted that informal and formal second-order validity should be respectively understood via (\ref{valopts}ii,iv).  In this case he would have been in a position to accept 
\begin{example}
$Val_s(\alpha^2) \rightarrow V_s(\alpha)$
\end{example}
But as is well-known -- and in fact is explicitly stressed by \citet[p. 81]{Henkin1950} -- the deductive completeness of $\proves_2$ with respect to the standard semantics \textsl{fails} -- i.e. we do \textsl{not} have
\begin{example}
$V_s(\alpha^2) \rightarrow D(\alpha^2)$
\end{example}
This blocks directly completing an squeezing argument for $\mathit{Val}_s(\alpha^2)$ in a manner analogous to the first-order case.  

But of course we might also try to salvage the argument by offering an independent case for 
\begin{example}
$\mathit{Val}_s(\alpha^2) \rightarrow V_h(\alpha)$
\label{vsh}
\end{example}
Relative to the characterization provided in (\ref{valopts}), $\mathit{Val}_s(\alpha^2)$ does indeed take into account \textsl{more} structures than does $V_h(\alpha^2)$ with respect to cardinality (i.e. classed-sized ones in addition to set-sized ones).  But note that  a standard model for the language of $\alpha^2$ may be regarded as a Henkin model $\langle \mathfrak{M},\mathfrak{S} \rangle$ wherein $\mathfrak{S}$ corresponds to the full powerset of the domain of $\mathfrak{M}$.   It thus follows that $\mathit{Val}_s(\alpha^2)$ takes into account \textsl{fewer} structures than $\mathit{Val}_h(\alpha^2)$ with respect to  potential variations in the range of second-order quantification.  It would thus seem to follow that the dictates of informal rigour themselves prevent us from accepting (\ref{vsh}) -- i.e. if we `push further' our intuitive understanding of second-order validity, we can see why (\ref{vsh}) does not hold.   And this in turn appears to cut off the possibility of providing an informally rigorous characterization of $\mathit{Val}(\alpha^2)$ in a manner which takes advantage of a completeness theorem with respect to a computably enumerable definition of second-order derivability.\footnote{\citet{Kennedy2017} consider a similar sequence of alternatives for extending Kreisel's validity argument to second-order logic.  But they ultimately suggest that $\mathit{Val}_s(\alpha^2) \rightarrow V_h(\alpha)$ should be accepted for, as they put it, `On the informal level it is impossible to see a difference between a standard model and a general [i.e. Henkin] model $\ldots$ The position taken here is that it is contrary to the idea of informal validity that one should be able to survey the situation from outside' (pp. 14-15).  We are personally sympathetic to such a view.  But we also take the example in the next paragraph to illustrate why this option was not open to Kreisel.}

To bring Kreisel's validity and CH arguments into direct contact, note finally that the formal independence of CH from $\mathsf{ZF}^2$ can in fact be used to provide a concrete illustration of the incompleteness of second-order derivability with respect to the standard semantics.   For recall that $\mathsf{ZF}^2$ is finitely axiomatizable -- say by a single $\mathcal{L}^2_{\mathsf{Z}}$-sentence $\zeta^*$.   By combining variants of Chuaqui and Weston's variant of Cohen's consistency proof for $\neg \mathrm{CH}$ and G\"odel consistency proof for CH we can hence see that $\not\proves_2 \zeta^* \rightarrow \chi^1$ and $\not\proves_2 \zeta^* \rightarrow \neg \chi^1$.   But of course the crux of Kreisel's CH argument is that one of $\models_2 \zeta^* \rightarrow \chi^1$ or  $\models_2 \zeta^* \rightarrow \neg \chi^1$ must hold (although we do not know which).   Thus a more precise reason why Kreisel may have demurred from attempting to provide an informally rigorous analysis of second-order validity was that he anticipated that its extension was of considerably greater complexity than that of first-order validity.\footnote{Kreisel develops this essentially epistemological point further in \citeyearpar[p. 191-194]{Kreisel1967c}.  Here he observes not only that the L\"owenheim number of second-order logic is likely to be much larger than that of first-order logic, but also that second-order validity cannot be effectively decidable.   He illustrates this further by observing that there is a broad class of statements $\alpha^2$ (akin to CH) for we can determine that the disjunction $V(\alpha^2)$ \textsl{or} $V(\neg \alpha^2)$ holds without coming to know which disjunct makes it true.  These considerations anticipate the later result of \citet{Vaananen2012} that the set of second-order formulas which are valid with respect to the standard semantics is $\Pi_2$-complete in the L\'evy hierarchy -- i.e. it is thus not $\Sigma^n_m$-definable for any level of the extended analytical hierarchy.  And this indeed does validate Kreisel's apparent suspicion that at least in \textsl{extension}, second-order validity is a considerably more complex concept than first-order validity (which is $\Pi^0_1$-complete in the arithmetical hierarchy).  This case can also be compared with Kreisel's apparent conclusions about \textsl{intuitionistic validity} as we discuss in \ref{intval}. \label{formnote}}

\section{The legacy of informal rigour}
\label{legacy}

As we forewarned in \S\ref{intro}, several obstacles stand in the way of coming to a general understanding of what Kreisel meant by `informal rigour'.   After illustrating some of the relevant textual and historical challenges in \S\ref{context}, we have attempted to factor them out in proposing schemas for formulating Kreisel's arguments in \S\ref{oninfrig}.  We then attempted to illustrate how his specific arguments conform to these schemas in \S\ref{casestudies}.   It is our hope that this framework will be use of use in coming both to a better understanding of Kreisel's contributions to a number of developments lying within the intersection of mathematical logic and philosophy of mathematics.   But more the hope also remains that by coming to a clear understanding of Kreisel's methodology will also allow us to extend the class of cases to which informal rigour may fruitfully be applied.

Rather than taking up this call here,  we will finally attempt to exit the forest of details we entered in \S\ref{casestudies} by discharging two more modest tasks.   First, in \ref{legsum} we will  summarize the case studies considered in \S\ref{casestudies} so as to further highlight some of their connections between them which emerge in light to our proposed schematizations.   Second, in \S\ref{legbar} we will further contextualize Kreisel's methodological perspective by recounting his exchange with Yehoshua Bar-Hillel about the relationship between informal rigour and Carnap's \textsl{method of explication}.

\subsection{Programmatic summary}
\label{legsum}

The examples we have examined in \S\ref{casestudies} suggest that it is ultimately straightforward to conform the examples which Kreisel presented as instances of informal rigour to the schemas IR and S we have presented in \S \ref{oninfrig}.  But although this suggests that these cases are \textsl{good examples} of Kreisel's proposed methodology,  the question also remains as to whether they are \textsl{good arguments} in their own right.  We have already taken some initial steps in evaluating the arguments in the course of considering how they have been received.  What we will do now is to record a few more programmatic observations which can be gleaned by comparing the structure of the arguments themselves.

The validity argument as reconstructed in \S\ref{validity} is Kreisel's best known application of informal rigour, eclipsing his other arguments in secondary sources and even serving as a stand-in for his overall methodology in the eyes of subsequent commentators.   We have suggested that this argument is indeed paradigmatic of the method of squeezing embodied by schema S and also that it is plausible to regard the constitutive principles for validity on which the argument relies as derivable by `pushing a bit farther than before the analysis of the intuitive notions'.  In light of this, we have also suggested that the argument provides an \textsl{extensionally adequate} characterization of the concept of first-order validity ($\mathit{Val}$) despite the fact it does so by treating the notion of validity \textsl{intensionally} (see note \ref{intnote}) -- a contrast we will suggest in \S\ref{legbar} bears on its standing as a successful instance of conceptual analysis.   We have also suggested that the validity argument is characteristic of squeezing arguments in the sense that the relevant Kreiselian theory $\mathsf{T}_K$ in which we have proposed that the argument should be conducted is \textsl{conservative} over the relevant precise theory $\mathsf{T}_P$ (i.e. $\mathsf{WKL}_0$).

 Kreisel's creating subject argument as reconstructed in \S\ref{gmp} has received much less attention outside the context of its original presentation.  But we have suggested that it represents an equally paradigmatic instance of the schema IR.   As Kreisel suggested, the argument highlights the potential role of \textsl{novel concepts} not previously accepted as proper parts of mathematical discourse for addressing open mathematical problems. As we have seen, the argument provides the relevant sort of resolution (i.e. a refutation of GMP) precisely because the Kreiselian theory $\mathsf{T}_K$ in which the principles linking the novel concept (i.e. the creating subject) to the background language mathematical is a \textsl{non-conservative} extension of the relevant precise theory $\mathsf{T}_P$ (i.e. $\mathsf{FIM}_0$).

Kreisel's CH argument as reconstructed in \S\ref{ch} is the most complex of the examples of formal rigour he presented in \citeyearpar{Kreisel1967b} and also the most difficult to assess.   At the macroscopic level, this argument can be understood as an instance of the schema IR whereby an answer is provided to the question: \textsl{Is the Continuum Hypothesis a mathematically definite statement}?   But we have also suggested that the concept \textsl{mathematical definiteness} ($\mathit{Def}$) is itself a common notion in need of further analysis before an answer to Kreisel's question can be accepted as definitive.  There is at least the possibility of providing such an analysis by the method of squeezing.  But rather than proceeding in this manner, Kreisel himself settled for assessing his proposed definition relative to a class of test cases.   We have suggested that reflection on these examples illustrates that his definition of $\mathit{Def}$ may not be extensionally adequate.\footnote{We will see in \S\ref{nonstandard} that a similar problem appears to beset his argument for priority of standard models over nonstandard ones.}  But in addition to this, we have also observed that the precise formalizability of Kreisel's definition together with the fact no additional principles are assumed to hold of $\mathit{Def}$ entails that the relevant Kreiselian theory $\mathsf{T}_K$ is \textsl{conservative} over the mathematical theory $\mathsf{T}_P$ (i.e. $\mathsf{ZF}_2$).   The fact that $\mathsf{T}_K$ proves the statement $\mathit{Def'}(\ulcorner \chi_1 \urcorner)$ expressing the definiteness of CH in Kreisel's sense thus requires further scrutiny before it can be accepted as providing an informally rigorous answer to the question which he hoped to resolve.

\subsection{Kreisel, Bar-Hillel, and Carnap}
\label{legbar}

One of the few portions of \citeyearpar{Kreisel1967b} on which we have not yet touched is Kreisel's exchange with Bar-Hillel.  In his brief reply to Kreisel, Bar-Hillel makes two basic points:

{\footnotesize 
\begin{quote}
Even the heuristic value of reflection can be impaired by taking the expression `reflecting about $\ldots$' too seriously. I would certainly object against taking too seriously the picture $\ldots$ that somewhere there are certain mathematical entities around whose exact nature is somehow veiled to a normal mortal and which reveal themselves only to those who know how to make good use of their reflective capacities.\\

It would also be of help to some of us if we could understand your notions of \textsl{informal rigour} vs. \textsl{formal rigour} as being closely similar to Carnap's pair of \textsl{clarification of the explicandum} vs. \textsl{providing the explicatum}. This identification would be of particular importance if I am correct in assuming that you intend your pair of notions to be used not only in the philosophy of mathematics but in the philosophy of science in general.  \hfill \citeyearpar[p. 172]{Kreisel1967b}
\end{quote}
}

Bar-Hillel's first comment draws attention to Kreisel's use of `reflection about' common concepts in the application of informal rigour.   Although this is not an expression which Kreisel himself employed in \citeyearpar{Kreisel1967a} or \citeyearpar{Kreisel1967b}, we have suggested that some form of reflection plays a role at steps Ib,c in the schema IR.  We will also return in \ref{finproof} to describe Kreisel's promotion of such a process in other sources as well as its relation to formal \textsl{reflection principles} as part of his characterization of finitist mathematics.   But what matters most in the present context is that rather than attempting to explain or qualify his views about `reflection', Kreisel issued the following direct rejoinder:

{\footnotesize
\begin{quote}
What also can make Bar-Hillel suggest that it is extraordinary for normal mortals to use their reflective capacities? Maybe some of us don't use them very well, but do we use our other capacities so much better? \P\  If I were really convinced that reflection is extraordinary or illusory I should certainly not choose philosophy as a profession; or, having chosen it, I'd get out fast.  \hfill \citeyearpar[p. 178]{Kreisel1967b}
\end{quote}
}

The more substantive part of Kreisel's reply pertains to the analogy between informal rigour and Carnap's method of explication which Bar-Hillel proposed.  About this Kreisel remarks 

{\footnotesize
\begin{quote}
Concerning the equation 
$$\frac{\textrm{clarification of the explicandum}}{\textrm{providing the explicatum}} = \frac{\textrm{informal rigour}}{\textrm{formal rigour}}$$

two things are to be said. First, strictly speaking, the equation does not hold because Carnap certainly denies the possibility of informal \emph{rigour} or \emph{proof}; he would not 
accept the problem of finding the \emph{correct} explicatum and proving it, but speaks of `replacing' the prescientific explicandum by an `adequate' explicatum. Carnap does not reject 
the possibility of proof outright, but feels convinced of the impossibility or fruitlessness of such a proof as a result of his experience.  The examples of the paper are intended to remind us of fruitful cases.
\hfill \citeyearpar[p. 176]{Kreisel1967b}
\end{quote}
}

In order to appreciate Kreisel's comments, it is useful to recall one of the best-known passages in which Carnap originally introduced the expressions `explicandum' and `explicatum':\footnote{Although Kreisel treated Bar-Hillel as a proxy for Carnap, Carnap himself had spoken on the same day Kreisel delivered the address on which \citeyearpar{Kreisel1967b} is based.  However no direct interaction between them is recorded in \citep{Lakatos1967}.} 

{\footnotesize
\begin{quote}
The task of \textsl{explication} consists in transforming a given more or less inexact concept into an exact one or, rather, in replacing the first by the second. We call the given concept (or the term used for it) the \textsl{explicandum}, and the exact concept proposed to take the place of the first (or the term proposed for it) the \textsl{explicatum}. The explicandum may belong to everyday language or to a previous stage in the development of scientific language.  The explicatum must be given by explicit rules for its use, for example, by a definition which incorporates it into a well-constructed system of scientific either logico-mathematical or empirical concepts.\hspace*{40ex} \citep[p. 3]{Carnap1962}
\end{quote}
}

There is thus indeed a \textsl{prima facie} similarity between informal rigour and explication in the sense that both Kreisel and Carnap describe their methods as seeking to make sense of the use of informal concepts in everyday practice so that they may be better integrated into mathematical or scientific reasoning.  On the other hand, the following characteristics of explication -- which are highlighted across many of Carnap's presentations which descend from \citeyearpar{Carnap1937} -- may be cited in order to distinguish its intended domain of application, methods, and goals from those of informal rigour:
\begin{enumerate}[i)]
\item Carnap characterized explication as a process by which `inexact' everyday
concepts are replaced by `exact' ones.   He suggests that the need for such replacements arises not only due to the fact inexact concepts are `pre-scientific' and thus potentially ambiguous -- as in the case of the concept \textsl{fish} \citeyearpar[\S1.5]{Carnap1962} -- but also that they are sometimes bound up in philosophical controversies -- as in the case of the concept \textsl{probability} (\citeyear{Carnap1945}, \citeyear[\S II]{Carnap1962}).  The process of explication may thus render what was originally regarded as a single everyday concept into multiple candidates amongst which we may choose an exact replacement.   
\item Carnap regarded the process of replacing pre-scientific concepts with exact ones as taking place within what he called a \textsl{linguistic framework} (which he often described as `formal').   He additionally took the adoption of such frameworks as inducing a distinction between questions which can be meaningfully posed `internally' within a framework itself and those which can only be asked `externally' in regard to the choice between frameworks (e.g. \citeyear[p. 20]{Carnap1950}).   
\item One of Carnap's goals in promoting explication was that of showing how philosophical controversies may be dissolved in favor of logical concerns once it is seen how concepts may be formalized within different frameworks (e.g. \citeyear[p. xiii]{Carnap1937}).    In prototypical cases the only remaining non-formal (i.e. `philosophical') questions will be those relating to the practical choice between frameworks.  After the process of explication has been carried out, it may thus no longer make sense to ask whether the exact concept which has been obtained as a surrogate for a given pre-scientific one provides a correct analysis of the original pre-scientific concept.   This is the basis of Carnap's well-known `principle of tolerance' (e.g. \citeyear[\S 17]{Carnap1937}).   
\end{enumerate}

Kreisel acknowledged \citeyearpar[p. 177]{Kreisel1967b} that he had only a partial understanding of Carnap's project.    But even on the basis of the foregoing summary, it is easy to see why he would have regarded both the methodology of informal rigour and the sorts of results he took it to be capable of establishing to differ considerably from the corresponding aspects of explication.   For instance the following contrasts with i)-iii) are apparent:
\begin{enumerate}[i$'$)]
\item Unlike Carnap, Kreisel leaves open the possibility that some noteworthy common concepts are not vague or ambiguous despite the fact that they may be `pre-scientific' in the sense that they initially lack accepted precise definitions.   This is paradigmatically illustrated by Kreisel's remarks about the concept of \textsl{validity} which we discussed in \S\ref{validity} -- i.e. `Nobody will deny that one knows more about $\mathit{Val}$ after one has established its relation with $V$ and $D$; but that doesn't mean that $\mathit{Val}$ was vague before' \citeyearpar[p. 154]{Kreisel1967b}.  And it is also illustrated by Kreisel's comments about the concept \textsl{set} (e.g. \citeyear[p. 144-145]{Kreisel1967b}) which he explicitly observed has \textsl{not} bifurcated in mathematical practice despite the original appearance that it might have been ambiguous between several distinct notions.
\item According to the model we have proposed in \S\ref{oninfrig}, an informally rigorous argument is understood as one carried out jointly by using constitutive principles for common concepts, bridging principles which connect them to novel and mathematical concepts, as well as a precise  mathematical background theory.   It would thus seem that Kreisel regarded informal rigour as a means by which different domains of concepts can be profitably brought into contact via the multi-stage procedure we have attempted to codify via the schema IR rather than one which segregates them concepts into disjoint `frameworks'.   
\item In the passage reproduced above, Kreisel directly rebukes Carnap for failing to hold open the possibility that explication might lead to anything more than `fruitful' analyses of common concepts -- i.e. that he rejected out of hand that a given analysis could be shown to be \textsl{correct} in the manner of what Kreisel called \textsl{philosophical proof}. Instead of the goal of replacing common concepts by more fruitful exact ones, it is thus more accurate to think of informal rigour as seeking to show that they admit to precise analyses which are provably \textsl{correct}.\footnote{In  \citeyearpar[p. 205]{Kreisel1967a} Kreisel provides a more extensive formulation of this point in the course of his critique of  `positivist or pragmatist doctrines' which he characterizes as assuming that  `traditional philosophical questions are so ill-defined that there is no possibility of a precise solution'.  He continues by remarking that such views do not take seriously `experience in traditional mental philosophy, namely insight into such intuitive concepts as logical validity, mechanical process, elementary proof, to name a few $\ldots$ Consequently, [such views do] not accept as meaningful the question whether certain axioms (laws) for concepts are \textsl{correct} as an analysis of the understood concepts and thus rejects the possibility of informal rigour. Instead, it speaks of \textsl{replacing} these concepts by formally introduced concepts which are supposed to be useful or adequate for certain (more or less unspecified) purposes: clarification and explication (\textsl{sic}) are favourites.'}
\end{enumerate}

Despite the superficial similarities between explication and informal rigour to which Bar-Hillel drew attention, the foregoing comparisons illustrate how sharply Kreisel appears to have understood his methods to differ from those of Carnap.  One point which stands out is that Kreisel at least seems open to a form of \textsl{conceptual realism} within the mathematical domains which were his primary interest.   For in holding that particular analyses of concepts such as \textsl{validity}, \textsl{set}, and \textsl{mechanical procedure} are correct, Kreisel appears to commit himself not only to the fact that these concepts \textsl{exist} independently to a given Carnapian framework but also that we possess methods for \textsl{demonstrating} that certain precise analyses are correct relative to some form of independent standard.\footnote{See, for interest, \citeyearpar[pp. 204-205]{Kreisel1989}.   Such a view can be contrast with Reck's  \citeyearpar[p. 198]{Reck2012} recent assessment that `Carnap basically rejects the assumption that there are concepts in some Platonic sense, existing ``out there'', with which an explicatum could be compared'.} 

Making clear what this standard amounts to would itself be a delicate matter.  But in repeatedly likening his validity argument to a `philosophical proof' it seems reasonable to regard Kreisel as proposing that the method which it exemplifies -- i.e. the squeezing scheme S of \S\ref{3ss} -- should be understood as a means for discovering objective relations between concepts which are already implicit in the constitutive principles we accept about them.   And in accepting that reflection on concepts plays a role in informal rigour, Kreisel thus reveals a bit more about the view of concepts he presumably would have had to maintain to underwrite such a procedure.  

Having highlighted these aspects of Kreisel's proposal, it also becomes clear that had he 
articulated his view about the nature of concepts more systematically he would ultimately have come up against another traditional philosophical concern -- i.e. the so-called \textsl{paradox of analysis} originally described by \citet{Moore1903} and \citet{Langford1942}.   This is a problem which Carnap elected to confront directly and which he introduced via Langford's formulation as follows:
 
 {\footnotesize
\begin{quote}
If the verbal expression representing the analysandum has the same meaning as the verbal expression representing the analysans, the analysis states a bare identity and is trivial; but if the two verbal expressions do not have the same meaning, the analysis is incorrect. \hfill \citet[p. 63]{Carnap1956}
\end{quote} }
 
\noindent  Carnap famously illustrated the problem at issue by comparing the statements
\begin{example}
\begin{enumerate}[i)]
\item The concept \textsl{brother} is identical with the concept \textsl{male sibling}.
\item The concept \textsl{brother} is identical with the concept \textsl{brother}.
\end{enumerate}
\label{brother}
\end{example}
About this contrast Carnap remarked that (\ref{brother}i) conveys `fruitful information' while (\ref{brother}ii) is `quite trivial'.  We are thus left with an instance of Moore's original puzzle: `If the first sentence is true, then the second seems to make the same statement as the first $\ldots$ [B]ut it is obvious that these two statements are not the same' \citeyearpar[p. 63]{Carnap1956}.

Here the traditional terms \textsl{analysandum} and \textsl{analysans} can be understood as counterparts for either Carnap's contrast between \textsl{explicandum} and \textsl{explicatum} or Kreisel's between \textsl{common} and \textsl{precise concepts}.   We will not pause here to consider the details of how Carnap proposed to account for the apparent difference in the `informativeness' of (\ref{brother}i,ii) via his own method of intension and extension.\footnote{See, e.g., \citep[\S14-\S15]{Carnap1956}.}  It would appear, however, that the goal of an informally rigorous argument carried out according to schema S is precisely to isolate conceptual identities in a manner which leads to analogous pairs of statements.   For instance in light of Kreisel's validity argument we are lead to consider the relationship between the following pair:
\begin{example}
\begin{enumerate}[i)]
\item The concept \textsl{validity} ($\mathit{Val}$) is identical with the concept \textsl{true in all set-sized models} ($V$).
\item The concept \textsl{validity} is identical with the concept \textsl{validity}.
\end{enumerate}
\label{valid}
\end{example}

Accounting for the apparent contrast in informativeness between such statements was -- needless to say -- a major preoccupation of mid-20th century analytic philosophy.   Although Kreisel saw himself as working at the boundary of these developments with mathematical logic, it seems likely that he would have claimed greater allegiance  to the latter.  But in likening informal rigour to the `old fashioned' method of conceptual analysis, it would still seem to be incumbent on him to provide a response to Moore's puzzle.   To the best of our knowledge, Kreisel never engaged  with this issue directly.   But we will now close by briefly outlining three considerations he might  plausibly have adduced in this regard.

A first consideration derives from the fact that relative to the schematization we have proposed in \S\ref{3ss}, the conclusion of a squeezing argument takes the form  $\A x(C(x) \leftrightarrow \pi(x))$ -- i.e.  a \textsl{material biconditional} relating a common concept expressed by $C(x)$ (e.g. \textsl{validity}) with a potentially complex predicate  $\pi(x)$ expressed in a precise language (e.g. \textsl{truth in all set-sized models}).  Thus even if such arguments are likened to conceptual analyses, their conclusions should -- at least in the first instance -- only be understood as expressing the \textsl{extensional equivalence} of the relevant analysandum and analysans.   

A related clarificatory consideration derives from Kreisel's own remark that   `For foundations it is evident that intensional operations are fundamental $\ldots$ \textsl{extensional operators can be defined in terms of intensional ones, but not conversely}' \citeyearpar[p. 184]{Kreisel1967b}.  We have suggested in \S\ref{validity} Kreisel took his validity argument as a paradigmatic example of informal rigour was in part because it shows how the intensional aspects of the proto-notion of \textsl{set} (and hence \textsl{structure}) on which it depends are shown to be factored out by the relevant squeezing of $\mathit{Val}$ between $D$ and $V$ (see note \ref{intnote}).   On this basis it might reasonably be concluded that informal rigour is only intended to yield extensionally adequate analyses even in cases where the details of the analysis take intensional features of concepts into account.  And if this were the case, Kreisel would be under no obligation to demonstrate that his methods were capable of delivering analysans which have the `same information value' with their corresponding analysanda (or were otherwise `cognitively synonymous' or `uniformly intersubstitutable' with them, etc.).   

Such an extensional understanding of informal rigour may appear modest relative to the goals which Carnap and other theorists may have hoped to achieve via methods such as explication.   But a final consideration derives from the fact although statements of the form $\forall x(C(x) \leftrightarrow \pi(x))$ express only an extensional coincidence of concepts, the informally rigorous arguments by which they are demonstrated correspond to \textsl{proofs} in suitable Kreiselian theories $\mathsf{T}_K$.  For instance the validity argument may be understood as yielding not only the conclusion that $\mathit{Val}$ is coincident with $V$, but also a demonstration explaining \textsl{why} such a coincidence holds.

This is a feature which appears to distinguish informal rigour not only from explication and traditional conceptual analysis, but also latter-day proposals such as `conceptual engineering' (e.g. \citealp{Cappelen2018}).   But of equal significance is the fact that each of the examples we have considered in \S\ref{casestudies} attests that the `philosophical proofs' by which such extensional conclusions are provided may be genuinely non-trivial.   For  as we have seen, these examples rely not only on the interplay between common, novel, and precise concepts but also on substantial mathematical theorems in the background.   Had he chosen to do so, it would thus seem that Kreisel could have responded to the paradox of analysis by explaining the manner in which an identity between concepts delivered by an informally rigorous argument can be demonstrably \textsl{correct} while simultaneously being \textsl{informative} on the basis of the traditional model of how mathematical theorems may be \textsl{true} but \textsl{non-obvious}.\footnote{At this point it is perhaps otiose to observe that the result which forms the crux of Kreisel's validity argument -- i.e. G\"odel's Completeness Theorem -- itself takes the form of an identity statement between two precise concepts (i.e. the derivable statements $D$ and those true in all set-sized models $V$).   But an equally striking example of the same phenomenon is Kreisel's use of Kleene's identification of the hyperarithmatical and $\Delta^1_1$-definable sets (Theorem \ref{kleenethm}) in his analysis of \textsl{predicative definability} which we have reconstructed in \S\ref{preddef}.}

\appendix

\section{Additional examples of informal rigour}
\label{app}

The goal of the present Appendix is to provide a brief overview of other apparent examples of informal rigour in Kreisel's work which have not been treated above.  These divide into two categories.    First, there are the analyses of the concept \textsl{set} (or cognates such as \textsl{aggregate} or \textsl{class}) and of the distinction between \textsl{standard} and \textsl{nonstandard} models which Kreisel explicitly proposed as additional case studies in \citeyearpar{Kreisel1967a,Kreisel1967b} but did not develop as systematically as the examples we have considered in \S\ref{casestudies}.\footnote{This category also includes the notion \textsl{mechanical procedure} (and its relationship to Church's Thesis) whose role within Kreisel's overall framework we have discussed in \S\ref{words}.}  Second, there are the notions of \textsl{finitist proof}, \textsl{predicative definability}, and \textsl{intuitionistic validity} with which Kreisel engaged substantially in his mathematical work in the 1950s and 1960s but mentions only briefly in central sources considered here.    In each case we have limited ourselves to a sketch of the basic ideas required to understand how Kreisel's various proposals and results can be understood as examples of informal rigour as it has been characterized here and a brief evaluation aimed at making his proposals more accessible to subsequent investigation.

\subsection{The discovery of set theoretic axioms}
\label{set}

As we discussed at the beginning of \S\ref{ch}, Kreisel's presentation of set theory in his survey of mathematical logic \citeyearpar{Kreisel1965a} contains not only a technical summary which was state-of-the-art for its time, but also a historical discussion of how Zermelo originally discovered his now-familiar axioms.   He begins by remarking  `Zermelo's informal derivation of his axioms can be analyzed by formulating explicitly properties of the cumulative-type structure' (p. 101).   In the introduction to  \citeyearpar[p. 202]{Kreisel1967a}, he also suggests that `the notion of aggregate which is analysed by means of the hierarchy (theory) of types'  should be understood as one of the `successes of mathematical logic'.  One of his clearest articulations of this point then appears in \citeyearpar{Kreisel1967b}:
\begin{quote}
{\footnotesize [The notion] set \textsl{of} something, first described clearly by Russell, and especially, Zermelo, has proved to be marvelously clear and comprehensive.  $\ldots$ [A]xioms which are evidently valid for the particular notion isolated by Zermelo (cumulative type structure) give a formal foundation $\ldots$ for a great deal of present day mathematical practice.  Zermelo's analysis furnishes an instance of a rigorous \textsl{discovery of axioms} (for the notion of set) $\ldots$ What one means here is that the intuitive notion of the cumulative type structure provides a coherent source of axioms $\ldots$ \hfill \citeyearpar[pp. 143-144]{Kreisel1967b} }
\end{quote}

In \citeyearpar{Kreisel1967b} Kreisel largely presupposes that Zermelo's \citeyearpar{Zermelo1930a}  original axiomatization of set theory should be understood as a successful application of informal rigour.\footnote{See \citep{Isaacson2011a} for a more extensive reconstruction of what Kreisel presupposed about Zermelo's account.}   A more extensive discussion of this point in given \S1 of \citeyearpar{Kreisel1967a} which is entitled `Basic objects, mathematical realism'.   Kreisel here suggests that the process by which Zermelo arrived at his axioms can be understood as continuing the method by which Peano and Dedekind arrived at their axioms by reflecting on the practice of arithmetic and analysis.  He then likens this process to that by which fundamental physical laws are discovered by an investigation of macroscopic objects and the laws which they obey before arguing on the basis of this analogy for the independent existence of mathematical objects and the objectivity of truth value.   A central step in his presentation is the following claim: `The obvious source of the mathematical properties used are insights which we interpret as being about external objects. The reliability of these insights is quite overwhelming $\ldots$' (p. 220).

Kreisel approvingly cites \citet{Godel1964} in this context and there is indeed an affinity between his case for realism and G\"odel's famous remark that `we do have something like a perception $\ldots$ of the objects of set theory, as is seen from the fact that the axioms force themselves upon us as being true'.\footnote{See \citep{Moss1971} for a more extensive discussion of Kreisel's mathematical realism.}  But unlike G\"odel, Kreisel does not speak further about a \textsl{general} faculty of mathematical intuition as a source of evidence.\footnote{In both \citeyearpar[\S 3.14]{Kreisel1965a} and \citeyearpar{Kreisel1967a} Kreisel speaks of how an intuitive faculty he refers to as \textsl{visualization} can serve as a form of evidence for finitist mathematics.   However his goal in this context was to account for how the visualizability of the iteration of combinatorial operations can provide finitary justification for transfinite induction up to (but perhaps not exceeding) the ordinal $\varepsilon_0$.   On the other hand he explicitly \textsl{disclaims} \citeyearpar[pp. 239-240]{Kreisel1967a} an analogy between this sort of intuitiability and what would be required to `visualize' all denumerable ordinals of the sort which is sometimes cited in order to justify set theoretic principles.}  If we wish to assimilate Kreisel's remarks about the discovery of set theoretic axioms to our prior treatment of informal rigour it thus seems that we must ultimately confront the following question: Should his use of the expression  `the intuitive notion of the cumulative type structure' be understood as describing an \textsl{abstract object} whose properties we grasp by something akin to G\"odelian intuition or as a \textsl{concept} whose properties we grasp via its constitutive principles akin to how Kreisel suggests we understand notions such as validity?

It is presumably only the second conceptual alternative which allows Kreisel's remarks to be reconstructed on the model of informal rigour we have proposed in \S\ref{oninfrig}.  For if set theoretic axioms are arrived at by the sort of intuitive faculty which G\"odel describes, it seems that their discovery is mediated by a process quite unlike that which Kreisel described as `philosophical proof' as exemplified by the examples of \S\ref{casestudies}.   But despite the apparent parallelism between Kreisel's description of mathematical realism and G\"odel's views, some evidence that he in fact preferred the conceptual alternative is provided by the more detail discussions of set theory in \citeyearpar[\S 1]{Kreisel1965a}, \citeyearpar[\S3b]{Kreisel1967a}, and \citeyearpar[\S A.2]{Kreisel1967c}.\footnote{\citet{Martin2005}  suggests that G\"odel should in fact also be read as a `conceptual realist' in something like the manner of the second alternative.}   In each of these sources Kreisel characterizes the cumulative in terms of the iterated sequences of ranks  $R(0) = \emptyset, R(\alpha + 1) = \mathcal{P}(R(\alpha)), R(\lambda) = \bigcup_{\alpha < \lambda} R(\alpha)$.  He then claims that reflection on this characterization  leads to the axioms of Zermelo (or Zermelo-Fraenkel) set theory.   

In \citeyearpar[\S1.1]{Kreisel1965a} Kreisel took the additional step of proposing `basic laws' which he describes as properties of this structure but are `not intended to ``define'' it axiomatically!'  (p. 101). These are stated in a two-sorted first-order language with variables for \textsl{sets} $x,y,z, \ldots$ and  $\xi,\eta,\zeta,\ldots$ for \textsl{types} (or `iterations')  together with the relations $x \in y$ and $x:\xi$ ($x$ \textsl{is of type} $\xi$).   Kreisel then states without proof that these principles have the same consequences as $\mathsf{ZF}$ in the language with just $\in$ (i.e. $\mathcal{L}^1_{\mathsf{Z}}$).\footnote{Kreisel's axioms include extensionality and comprehension for sets, the properness of the type structure (i.e. every set has a type, all members of a set of type $\xi$ have type $\eta < \xi$, the types form an infinite linear order).   It also includes a \textsl{supremum principle} which (taking into account a correction given on \citeyear[p. 272]{Kreisel1967a}) reads as follows: $\A x:\xi \E y \Phi(x,y) \rightarrow \E \eta \A x:\xi \E y:\eta \Phi(x,y)$.  Kreisel glosses this in terms of the reflection principle `what is true in the whole universe under consideration is true at some $\eta$'.   But although Kreisel cites \citep{Montague1965a} for the result that his system has the same set theoretic theorems as $\mathsf{ZF}$, the axiomatization in this paper -- i.e. what Montague calls \textsl{rank-free set theory} -- is actually a one-sorted second-order theory. \label{klt}}  We thus arrive at the further question whether the derivations in question can themselves be regarded as the components of an informally rigorous \textsl{argument} establishing that the $\mathsf{ZF}$ axioms indeed follow from  `the intuitive notion of the cumulative type structure'?

We will not attempt to provide a thorough answer here.  But it may also be noted that 
both Kreisel's axioms and the motivation he provides for them are similar to those proposed by \citet{Scott1974a} for what has come to be called his \textsl{theory of levels}.   Scott's derivation of the $\mathsf{ZF}$ axioms suggests that several of Kreisel's principles -- e.g. the linear ordering of types and his `least element' (or Regularity) axioms -- can be derived from more elementary assumptions about the type structure.  On the other hand, not only is Kreisel's theory first-order (meaning that it cannot be used to carry out his argument for the definiteness of CH), but the supremum principle described in note \ref{klt} is evidently included to allow for the derivation of the Replacement schema (certain instances of which we have seen in \S\ref{ch} are \textsl{non-definite} in Kreisel's sense).   The  question thus remains whether the sort of argument Kreisel indicates in \citeyearpar[\S 1]{Kreisel1965a} can be understood as providing a \textsl{precise} characterization of a specific $\mathcal{L}^1_{\mathsf{Z}}$- or $\mathcal{L}^2_{\mathsf{Z}}$-theory -- in the manner presumably called for by informal rigour -- rather than a piecemeal justification of a minimal set of axioms.\footnote{A similar issue would appear to arise for subsequent attempts to `justifiy' various set theoretic axioms on the basis of a given informal description of the iterative hierarchy. (See, e.g., \citep{Boolos1971},  \citep{Parsons1977}, \citep{Wang1977} and also the review \citep{Burgess1985a} for an assessment of these sources as applications of informal rigour similar to that suggested here.)   For although such accounts can be understood to provide justification for \textsl{accepting} certain basic axioms (e.g. Emptyset, Pairing, Union), their assimilation to the model of informal rigour we proposed in \S\ref{oninfrig} would also appear to require that the concept in question is articulated in a manner such that it also settles the status of Choice and Replacement (and perhaps other `axioms of infinity') as well as providing justification for \textsl{rejecting} other principles (e.g.  Constructibility or Determinacy) on the basis of similarly `intuitive' considerations.  But on Kreisel's view it seems likely that the set of `intuitively valid' $\mathcal{L}^1_{\mathsf{Z}}$-sentences will be highly complex -- e.g. since he desired it to decide formally independent statements like CH, it is unlikely to be recursively enumerable (see \S\ref{chsol}).   As we discuss further in \S\ref{intval}, Kreisel appears to have drawn a similar conclusion about the prospects for providing an informally rigorous analysis of the concept \textsl{intuitionistic validity}.   In the case of such (provably) complex concepts, one thus might reasonably begin to wonder at what point the method starts to lose its grip.}

\subsection{Standard versus nonstandard models of arithmetic}
\label{nonstandard}

Appendix B of \citeyearpar{Kreisel1967b} contains an argument about the relationship which nonstandard models of arithmetic models bear to the standard model which can be more straightforwardly assimmilated to the model of informal rigour presented in \S\ref{oninfrig}.\footnote{Kreisel's decision to include this section appears to have arisen from discussions he had with Abraham Robinson which -- like Kreisel's interaction with Myhill -- caused him to extend \citeyearpar{Kreisel1967b} after it was originally presented.   At this time Robinson was already widely known as an exponent of nonstandard analysis and in fact had delivered his address  `The metaphysics of the calculus' at the same 1965 conference.   It is likely, however, that Kreisel's interactions with Robison were also shaped by a paper \citeyearpar{Robinson1965a} which Robison has delivered a year earlier in which he explicitly presented a formalist reaction to Cohen's independence results in set theory.   As emerges more clearly in \citeyearpar{Kreisel1969a,Kreisel1971d} Kreisel thus regarded Robinson (along with Cohen) as a potential critic of informal rigour.}  Kreisel's argument is designed to answer the question posed at the beginning of the section:

{\footnotesize
\begin{quote} \textsl{Standard and nonstandard models}: The question we want to formulate is this: which comes first?   So to speak: which is more fundamental? \hfill \citeyearpar[p. 165]{Kreisel1967b}
\end{quote}
}
\noindent As Kreisel notes, this question is `technical' in the sense that the concept of \textsl{model} at issue is a mathematical notion rather than a truly `common' one.   For this reason he proposes to investigate it within the `basic conceptual framework of classical mathematics' together with its conventional existential assumptions about the existence of models (or sets).  He then notes

{\footnotesize
\begin{quote}
[C]urrent mathematics is full of related questions; when we look at an axiomatic system we do not merely look for properties common to all models of it (i.e. consequences), but for principal structures among them, e.g. minimal ones (the group defined by certain relations, and not only: all groups satisfying them), universal ones, etc. $\P$ In foundations too we find good answers to related questions: Which come first: (finite) ordinals or (finite) sets? $\ldots$ In short there is nothing outlandish in our question.
\end{quote}
}

Having framed the matter in this way, Kreisel then proceeds to propose a precise analysis of the relation \textsl{model} $\mathfrak{A}$ \textsl{comes before} (or \textsl{is more fundamental than}) \textsl{model} $\mathfrak{B}$.  Writing $\mathit{Fund}(\mathfrak{A},\mathfrak{B})$ for this relation Kreisel's definition can be stated in modern terms as follows:\footnote{\citet[p. 168]{Kreisel1967b} states this by simply saying that `$\mathfrak{B}$ is definable in $\mathfrak{A}$ by means of the language considered'.   But what is intended here seems to be precisely the definition of interpretability of structures employed in contemporary model theory -- i.e.  there exists a map $(\cdot)^*: \mathfrak{L}_{\mathfrak{B}} \rightarrow \mathfrak{L}_{\mathfrak{A}}$ associating the primitive expressions of $\mathfrak{L}_{\mathfrak{B}}$ with those of $\mathfrak{L}_{\mathfrak{A}}$ and a domain predicate $\delta(x)$ such that the $\mathfrak{L}_{\mathfrak{B}}$-structure $\mathfrak{B}^*$ with domain $B^* = \{a \in A : \mathfrak{A} \models \delta(a)\} \subseteq A$ and with non-logical symbols similarly interpreted in $\mathfrak{A}$ by their images under $(\cdot)^*$ -- e.g. $P^{\mathfrak{B}^*} = \{\vec{a} \in A^k : \mathfrak{A} \models P^*(\vec{a})\} \subseteq A^k$ for a $k$-ary predicate $P$ --  is such that $\mathfrak{B}^*$ is isomorphic to $\mathfrak{B}$.  See, e.g., \citep[\S 1.3]{Marker2002}.}
\begin{definition} $\mathit{Fund}(\mathfrak{A},\mathfrak{B})$ if and only if $\mathfrak{B}$ is model theoretically interpretable in $\mathfrak{A}$ but $\mathfrak{A}$ is not model theoretically interpretable in $\mathfrak{B}$.
\label{fund1}
 \end{definition}
 
Kreisel wished to employ this definition in order to argue for the conclusion that the standard model of arithmetic $\mathfrak{N} = \langle \mathbb{N},+,\times,0 \rangle$  `comes before' (or `is more fundamental than') nonstandard models of theories such as first-order Peano arithmetic [$\mathsf{PA}$].  In support of this contention, he cites the following two results:\footnote{Kreisel states these facts without citation.   It is reasonable to regard (\ref{fund2}ii) as folklore as it follows from the elementary observation that if an $\mathcal{L}^1_a$-model $\mathfrak{M}$ is nonstandard and $\delta(x)$ defined $\mathbb{N}$ in $\mathfrak{M}$, then the instance of the induction scheme for $\delta(x)$ would fail in $\mathfrak{M}$ from which it follows that it cannot satisfy $\mathsf{PA}$.   On the other hand, a version of (\ref{fund2}i) is stated by \citet[Proposition 4.1]{Scott1961} and can in fact be derived via an arithmetized completeness construction which \citet{Kreisel1950} had helped to popularize (see \citealp{Dean2020e}).}
\begin{theorem}
\begin{enumerate}[i)]
\item Let $\mathsf{T}$ be a consistent theory whose axioms $\mathrm{Ax}_{\mathsf{T}}$ are arithmetically definable -- i.e. $\mathrm{Ax}_{\mathsf{T}} = \{\psi : \mathfrak{N} \models \phi(\ulcorner \psi \urcorner)\}$ for some $\mathcal{L}^1_a$-formula $\phi(x)$.   Then there exists a nonstandard $\mathfrak{M} \models \mathsf{T}$ which is model theoretically interpretable in $\mathfrak{N}$.   
\item Let $\mathfrak{M} \models \mathsf{PA}$ be nonstandard.   Then $\mathfrak{N}$ is not model theoretically interpretable in $\mathfrak{M}$.   
\end{enumerate}
\label{fund2}
\end{theorem}

It is easy to see why Kreisel may have wished to exploit these results to counter the `formalistic' challenge that since standard and nonstandard models of arithmetic both satisfy some basic set of axioms (e.g. $\mathsf{PA}$) they are `equally fundamental' or -- even more radically -- `all models are just manners of speaking' (p. 168).  In particular, Definition \ref{fund1} might be taken to provide a precise definition of the notion `more fundamental than' which serves a role parallel to that played by Kreisel's proposed definition of mathematical definiteness (\ref{md}) in his CH argument.   Now let $\mathsf{T} \supseteq \mathsf{PA}$ and $\mathfrak{M} \models \mathsf{T}$ be a nonstandard model which is interpretable in $\mathfrak{N}$ as in Theorem \ref{fund2}i. This result together with Theorem \ref{fund2}ii) can then be invoked to conclude that $\mathit{Fund}(\mathfrak{N},\mathfrak{M})$ -- i.e.  that $\mathfrak{N}$ indeed `comes before' or `is more fundamental than' $\mathfrak{M}$.

This reconstruction suggests that Kreisel's argument in Appendix B of \citeyearpar{Kreisel1967b} can indeed be formulated in the manner illustrated by the examples we have considered in \S\ref{casestudies}.   But in order to evaluate such an argument as an instance of `philosophical proof', more would need to be said about how the relations `comes before' or is `more fundamental than' figure in our standing mathematical practices.  One might, for instance, attempt to explain the intended interpretation $\mathit{Fund}(\mathfrak{A},\mathfrak{B})$  `genetically' in the sense that the existence (or constructibility, etc.) of $\mathfrak{B}$ presupposes (or otherwise depends on) the existence  of $\mathfrak{A}$ (or its relative constructibility, etc.).    On the other hand, one might attempt to explain the relation `conceptually' in the sense that our ability to grasp (or intuit, etc.) the structure $\mathfrak{B}$ presupposes the ability to grasp the structure $\mathfrak{A}$.  

Kreisel says little to distinguish between these alternatives himself.  One can envision him thinking that the task of providing such an explanation fell to his formalist interlocutors who may have wished to maintain that nonstandard models of arithmetic are `equally fundamental' (or `come at the same time as') the standard model.  Nonetheless, it is also possible to cite  results which appear to tell \textsl{against} the plausibility of regarding Definition \ref{fund1} as demarcating a previously recognized notion on either of the proposed interpretations of $\mathit{Fund}$.  For instance consider the following facts, which are respectively consequences of results of Julia Robsinson \citeyearpar{Robinson1949} and Solomon Feferman \citeyearpar{Feferman1958}:
\begin{example}
\begin{enumerate}[i)]
\item Let $\mathfrak{Q} = \langle \mathbb{Q},+,\times,0\rangle$.   Then since $\mathfrak{Q}$ is model theoretically interpretable in $\mathfrak{N}$ and $\mathfrak{Q}$ is model theoretically interpretable in $\mathfrak{N}$,  $\neg \mathit{Fund}(\mathfrak{N},\mathfrak{Q})$ and $\neg \mathit{Fund}(\mathfrak{Q},\mathfrak{N})$ both hold.
\item Let $\mathfrak{M} \models \mathsf{TA}$ -- i.e. \textsl{true first-order arithmetic} -- be nonstandard.   Then $\mathfrak{M}$ is \textsl{not} model theoretically interpretable in $\mathfrak{N}$ and thus $\neg \mathit{Fund}(\mathfrak{N},\mathfrak{M})$.   
\end{enumerate}
\label{fund3}
\end{example}

Fact (\ref{fund3}i) shows that the structure of the natural numbers and that of the rational numbers are mutually interpretable within one another.   On the other hand, $\mathfrak{N}$ is typically understood as `coming before' $\mathfrak{Q}$ on the traditional `genetic' view about how the rational numbers are constructed as ratios of natural numbers.   But neither structure is `more fundamental than' the other according to Kreisel's criterion.  On the other hand, while Theorem \ref{fund2}i shows that 
every consistent $\mathcal{L}^1_a$-\textsl{theory} $\mathsf{T}$ with an arithmetically definable set of axioms has a nonstandard model interpretable within $\mathfrak{N}$,  (\ref{fund3}ii) shows that this result cannot be reformulated to speak of arbitrary \textsl{models} of $\mathsf{T}$ (as this class includes nonstandard models $\mathfrak{M} \models \mathsf{TA}$).    If we think of models as being given to us `conceptually' via their first-order theories, Kreisel's criterion thus only allows us to exclude as  `less fundamental than $\mathfrak{N}$' those nonstandard models which fail to satisfy some true arithmetical sentence. 

It seems likely that Kreisel was aware of results like these which tell against the descriptive adequacy of his definition of $\mathit{Fund}$.   In fact he remarked that `[I]t is not claimed that the definability criterion above formulates fully the question posed above. But it seems to be at least a sane step towards taking it seriously instead of leaving it at an empty level' \citeyearpar[p. 169]{Kreisel1967b}.  Thus despite its shortcomings the argument can still be seen as attempting to fulfill the goal of informal rigour  `not to leave undecided questions which can be decided by full use of evident properties of these intuitive notions'.

\subsection{Foundational standpoints}
\label{standpoints}

Much of Kreisel's work in the 1950s and early 1960s related to technical questions arising in light of the foundational programs we now call  finitism, predicativism, and intuitionism.   In contemporary philosophy of mathematics, these schools are often associated with specific historical figures, canonical texts, and (in certain cases) axiomatic systems.   On the other hand, Kreisel was not only working at a time when some of the programs were ongoing -- and thus to some extent still in flux -- but he was able to interact with several of the original authorities. 

We have already seen that Kreisel used the contrasts between foundational standpoints to structure his technical  survey of mathematical logic \citeyearpar{Kreisel1965a}.   But he also employs the same framework in \citeyearpar{Kreisel1967a} to systematically distinguish the forms of \textsl{evidence} which he took to characterize the standpoints.   
About this process he remarks

{\footnotesize
\begin{quote}
It is to be remarked that [certain kinds] of evidence have been suggested by the traditional literature on the Philosophy of Mathematics: what is done here is to sharpen the notions involved and to see whether they can be formalized, i.e. formulated with formal rigour $\ldots$ It is to be expected that mathematical logic will not for ever limp behind traditional notions, but that the technical results will suggest new significant kinds of preferred evidence.  \hfill \citeyearpar[p. 233]{Kreisel1967a}
\end{quote}}

As we have seen, Kreisel's creating subject argument is prototypical of how he envisioned such a process in the sense that it attempts to employ a mathematical result  to transform Brouwer's original weak counterexample into a formal refutation of Generalized Markov's Principle.   But in fact much of his earlier work also involved suggesting similar analyses of proposals in the works of figures such as Hilbert, Bernays, Russell, Poincar\'e, Brouwer, Heyting, and G\"odel.    Kreisel only started using the expression `informal rigour' to describe the method of reflecting on `traditional results' to obtain `sharpened notions' in the mid-1960s.  But it is evident from  allusions such as that given on \citeyearpar[p. 157]{Kreisel1967b} that he understood much of his earlier work to exemplify  such a process.

\subsubsection{Finitist proof}
\label{finproof}

Kreisel's earliest work in mathematical logic was informed by a careful reading of the two volumes of Hilbert and Bernays's \textsl{Grundlagen der Mathematik} \citeyearpar{Hilbert1934,Hilbert1939}.  The first two chapters of \citeyearpar{Hilbert1934} contain a detailed exposition of what Hilbert and Bernays themselves called `the finitary perspective'  for which Kreisel (e.g. \citeyear{Kreisel1958a,Kreisel1965a}) himself later popularized the terms `finitist' and `finitism'.   This includes a description of an arithmetical system containing terms introduced for recursively defined functions and quantifier free induction similar to a free variable formulation of what is now called \textsl{Primitive Recursive Arithmetic} $[\mathsf{PRA}]$.  \citet[p. 242]{Kreisel1951} originally adopted this characterization as a component of what he referred to as a \textsl{finitist proof} -- i.e.  one which involves only `free variable formalism with various forms of induction'.   

Kreisel retained this analysis until the late 1950s by which time he had started to study \textsl{proof theoretic reflection principles} -- i.e. schema typified by $\mathsf{Prov}_{\Sigma}(\ulcorner \phi \urcorner) \rightarrow \phi$ asserting that if an arithmetical formula $\phi$ is provable in a given axiomatic theory $\Sigma$, then $\phi$ is true in the standard model.   In the case that $\Sigma$ is a recursively axiomatizable theory certain instances of this schema will not be provable in $\Sigma$ in virtue of G\"odel's incompleteness theorems -- e.g., the case with $\phi = \bot$ is equivalent to $\mathrm{Con}(\Sigma)$.   This means the result of adding the prior scheme to $\Sigma$ will result in a stronger theory.  

In his address `Ordinal logics and the characterization of informal concepts of proof' \citeyearpar{Kreisel1960e} Kreisel described a procedure for iteratively adjoining statements of this form to obtain a sequence of theories $\Sigma_0 = \mathsf{PRA}, \Sigma_1, \Sigma_2, \ldots$ indexed by ordinals wherein $\Sigma_{\alpha +1}$ is obtained by closing $\Sigma_{\alpha}$ under the condition `if $\mathrm{Prov}_{\alpha}(\ulcorner \phi(\overline{n}) \urcorner)$ has been established by finitist means then, on the intended meaning of free variables, $\phi(\overline{n})$ is finitistically established'.\footnote{Kreisel's work on such \textsl{autonomous progressions} of theories was inspired by Turing's \citeyearpar{Turing1939} earlier \textsl{ordinal logics} and was extended further by \citet{Feferman1962}. \citet{Feferman1995a} and \citet{Franzen2004} provide popularized accounts in which Kreisel's role is minimized.}  About such a sequence Kreisel writes

{\footnotesize
\begin{quote}
[W]e propose to characterize finitist proofs by a precisely defined class of formal systems, namely the least \textsl{class} of systems $\Sigma_{\mu}$ containing a certain basic finitist apparatus and closed under the principle: if a proof predicate $\mathrm{Prov}_{\nu}(x)$ is recognized as such in a system $\Sigma_{\mu}$ of the class then the corresponding system $\Sigma_{\nu}$ also belongs to the class.  \P\  Our main result is that, in a precise sense, the theorems of this class are co-extensive with those of classical number theory when the latter is suitably interpreted. $\ldots$ Since each of our extensions is finitist this means at least that finitist results include essentially those of classical number theory.  \hfill \citeyearpar[p. 290]{Kreisel1960e}
\end{quote}
}

\noindent Kreisel would go on to describe similar progressions of theories in \citeyearpar[\S 3.4]{Kreisel1965a} and \citep{Kreisel1968a} wherein it is indeed shown that by starting with $\mathsf{PRA}$ and adjoining the so-called \textsl{uniform reflection principle} $\A x \mathrm{Prov}_{\Sigma_{\alpha}}(\ulcorner \phi(\dot{x}) \urcorner) \rightarrow \A x \phi(x)$ all of the theorems of $\mathsf{PA}$ are obtained (and potentially even stronger systems depending on the length of the iteration).   It is on the basis of this result which Kreisel appears to have updated his prior characterization of the extent of finitist proof to conclude that `finitist results include essentially those of classical number theory'.    

The fidelity of such a characterization of Hilbert and Bernays's original description of finitism remains controversial.\footnote{See, e.g., \citep[\S 1, \S 2]{Tait2005}, \citep{Dean2015a}, and \citep[\S 3]{Dean2017b}.} But such an identification is still prototypical of how Kreisel proposed to employ mathematical results together with reflection on the informal concepts to provide a precise characterization of the class of statements which are accessible from within a given foundational standpoint in terms of an axiomatic theory.  He would later describe this process as follows: 

{\footnotesize 
\begin{quote}
What principles of proof do we recognize as valid once we have understood (or, as one sometimes says, `accepted') certain given concepts? The process of recognizing the validity of such principles (including principles for defining new concepts, that is, formally, of extending a given language) is here conceived as a process of reflection; reflecting on the given concepts, reflecting on this process of reflection, and so forth. \hfill \citeyearpar[p. 489]{Kreisel1970a}
\end{quote}
}
\noindent  The centrality of reflection to Kreisel's conception of informal rigour is made clear in his exchange with Bar-Hillel which we have recounted in \ref{legbar}.  But as he does not attempt to provide a general description of how such a reflective process operates in the primary sources we have considered, we will not attempt a further reconstruction here. 

\subsubsection{Predicative definability}
\label{preddef}

Kreisel's interest in the notion of \textsl{predicativity} was continuous with his work on finitism.  At the end of \citeyearpar{Kreisel1960e} he suggested that an analysis of \textsl{predicative proof} can be obtained via a hierarchy of theories similar to that just described but formulated in the language of ramified second-order arithmetic and employing the extension principle  `if $<$ is proved to be a well-ordering in a system $\Sigma_{\mu}$, then the system with types indexed by $<$ is said to be proved in $\Sigma_{\mu}$ to be permissible (as a predicative proof predicate)' (p. 297).\footnote{Such an account was developed further by \citet{Feferman1964}.   After further summarizing Kreisel's analysis of predicative definability (Part I), Feferman then presents his own analysis of predicative provability  in terms of two transfinite progressions of theories which we calls $\mathsf{R}_{\alpha}$ and $\mathsf{H}_{\alpha}$ (Part II).  These may in turn be understood as proof-theoretic surrogates for the roles played by the hyperarithmetical sets and the $\Delta^1_1$-definable sets in Kreisel's proposal.  However one respect in which Kreisel and Feferman's account differ is with respect to how far along the ordinals predictive justification can be found for iteration of the set- or theory-forming operations. For reasons described below, Kreisel suggested taking $\omega^{ck}_1$ as a natural stopping point for predicative definability.  On the other hand, it was in this context in which Feferman originally proposed that the smaller ordinal $\Gamma_0$ serves as an upper bound on iterating predicative provability in virtue of the fact that it is the least non-autonomous ordinal of the progressions $\mathsf{R}_{\alpha}$ and $\mathsf{H}_{\alpha}$. Kreisel would later speaking approving of Feferman's proposal in \citeyearpar[pp. 174-177]{Kreisel1965a} and \citeyearpar[pp. 240-241]{Kreisel1967a}.}  But he also observed that such a proof-theoretic analysis does not touch the notion of \textsl{predicative definability} which figures centrally in the work of Russell and Poincar\'e.   

Kreisel had in fact already suggested such an analysis in his review of the paper in which \citet{Kleene1955a} had originally introduced the analytical and hyperarithmetical hierarchies: 

{\footnotesize
\begin{quote}
In the reviewer's opinion the class [of hyperarithmetical predicates] provides a precise and satisfactory definition of the notion of predicative sets (of integers), based on the concept of constructive ordinal.
\end{quote}
}

\noindent Kreisel  developed this point with considerable sophistication in several later papers of which his  address `La Pr\'edicativ\'e' \citeyearpar{Kreisel1960} is most central.   Therein he remarks

{\footnotesize
\begin{quote}
Without affirming the identification of predicative definitions (of sets of natural numbers) with the class [$\mathrm{HYP}$] of hyper-arithmetic definitions of Kleene, I will describe some results which bear on this identification: they demonstrate, for  [$\mathrm{HYP}$], properties which are evident for the intuitive notion of predicativity. \hfill (p. 373)
\end{quote}
}

Despite the tentativeness with which Kreisel expressed this proposal, the argument given in \citeyearpar{Kreisel1960} can in fact be understood as a paradigmatic instance of the method described in \S\ref{3ss} whereby the `intuitive notion of predicativity' is squeezed between two precise characteriztions of the hyperarithmetical sets.\footnote{Kreisel discusses the basis of this analysis further in \citeyearpar[\S 3, \S 8]{Kreisel1962c} and \citeyearpar[\S 1.5, \S 3.5]{Kreisel1965a}.   Reconstructions similar to that given here are also sketched in \citep[p. 604-605]{Feferman2005}, \citep[\S 4]{Walsh2016}, and \citep[\S 2.3]{Dean2017b}.}  The basis of Kreisel's analysis descends from two distinct conceptions of predicative definability which Kreisel locates in the work of Russell and Poincar\'e.   These can respectively be understood as giving rise to \textsl{narrow} and \textsl{wide} definitions which are then shown mathematically to both pick out the class of hyperarithmetical sets. 

The narrow conception of predicativity derives from the traditional negative characterization of an \textsl{impredicative definition} as one which violates Russell's well-known \textsl{vicious circle principle} -- i.e. `If, provided a certain collection had a total, it would have members only definable in terms of that total, then the said collection has no total' \citeyearpar[p. 225]{Russell1908}.   On this basis Kreisel proposed the following positive characterization of `\textsl{the fundamental idea of predicativity}' -- i.e. `in a predicative definition only quantifiers relating to already constructed sets are used' \citeyearpar[p. 377]{Kreisel1960}.

Kreisel suggested that this characterization could be sharpened in the manner of an informally rigorous argument by starting from another fundamental presupposition of Poincar\'e (e.g. \citeyear{Poincare1913b}) -- i.e. that since the natural numbers are given by intuition, they can be accepted as a completed infinite totality $\mathbb{N}$.    Kreisel also employed a related idea often credited to \citet{Weyl1918} -- i.e. since $\mathbb{N}$ is accepted as a completed totality, the \textsl{arithmetical sets} -- i.e. those sets $X \subseteq \mathbb{N}$ such that $X= \{n \in \mathbb{N} : \mathfrak{N} \models \phi^1(\overline{n})\}$ for some fomrula $\phi^1(x) \in \mathcal{L}^1_a$ whose definition involves only quantification of the natural numbers -- can be accepted as predicatively defined.   Building on this, he then described a structure  introduced by \citet{Kleene1955a} known as the \textsl{ramified analytic hierarchy}.  This hierarchy $\mathrm{RA}_{\xi}$ is defined up to an ordinal $\xi$ as follows: $\mathrm{RA}_0$ corresponds to the class of arithmetical sets; $\mathrm{RA}_{\alpha + 1}$ corresponds to the class of sets definable in the language $\mathcal{L}^2_a$ of second-order arithmetic with second-order quantifiers restricted to $\mathrm{RA}_{\alpha}$; $\textsc{RA}_{\alpha} = \bigcup_{\beta < \alpha}$ where $\alpha < \xi$ is a limit.

This definition leaves open the length of the iteration and thus does not itself provide a completely precise definition for this concept.   But in this regard Kreisel cites a theorem of \citet{Spector1955} who showed that the recursive well-orderings of the natural numbers correspond to the ordinals less than $\omega^{ck}_1$ -- i.e. the least ordinal which is not named within the system of ordinal notations known as \textsl{Kleene's} $O$.   And as he also observes, \citet{Kleene1959} had  shown that the class of hyperarithmetical sets $\mathrm{HYP}$ is in fact equal to $\textsc{RA}_{\omega^{ck}_1}$.\footnote{Although in \citeyearpar{Kreisel1960} Kreisel takes membership in $\textsc{RA}_{\omega^{ck}_1}$ as \textsl{defining}  `$X$ is hyperarithmetical', \citet[p. 210]{Kleene1955a}  had originally defined this notion in terms of the \textsl{hypearithemtical hieararchy} $\bigcup H_a$ formed by transfinite iterations of the Turing jump -- i.e. $X$ is hyperarithmetical just in case $X$ is recursive in $H_a$ for some $a \in O$.} On this basis Kreisel suggests that the membership of a set in this hierarchy `puts in a precise form the intuitive idea expressed by ``already'''  \citeyearpar[p. 377]{Kreisel1960} as it appears in the prior informal characterization of predicativity.   

Suppose we now introduce the primitive predicate $\mathit{Pred}(X)$ intended to express that $X$ is predicatively definable.   The first step in Kreisel's argument can be reconstructed as making a case for the adoption of the following principle as a narrow characterization of predicative definability: 
\begin{example}
$\A X \subseteq \mathbb{N}(X \in \textrm{HYP} \rightarrow \mathit{Pred}(X))$
\label{pred1}
\end{example} 
This expresses that the membership of a set of natural numbers in one of the classes $\mathrm{RA}_{\alpha}$ for $\alpha < \omega^{ck}_1$ is a \textsl{sufficient} condition for predicative definability.  

The wide characterization of predicativity which Kreisel employed was inspired by an alternative description  which he also locates in the writings of Poincar\'e:  

{\footnotesize
\begin{quote}
The following theorems concern an another idea of  \citet[p. 47]{Poincare1910} concerning the definition of predicativity: a definition of $D$ is said to be predicative if an enlargement of the class of the sets considered does not change the set defined by $D$. \hfill \citeyearpar[p. 378]{Kreisel1960} 
\end{quote}
}

\noindent Indeed a portion of the original passage which Kreisel cites reads as follows:\footnote{See \citep{Hallett2011} for further discussion of the relationship between Poincar\'e's characterization of predicativity and the contemporary notion of \textsl{absoluteness} which we are about to see that Kreisel took to be described here.}

{\footnotesize
\begin{quote}
I call a classification predicative if it is not changed by the introduction of new elements $\ldots$ What is here meant by the word `predicative' is best illustrated by an example. If I am to deposit a set of objects into a number of boxes two things can occur: either the objects already deposited are conclusively in their places, or, when I deposit a new object, I must always take the others out again (or at any rate some of them). In the first case I call the classification predicative, in the second not. \\ \hspace*{1ex} \hfill  \citep[p. 47]{Poincare1910}
\end{quote}
}

The informal notion described here is thus that of the rigidity of the class picked out by a definition relative to expansions of the domain in which it is evaluated.  But as Kreisel also notes

{\footnotesize
\begin{quote}
At first glance this notion is broader than the fundamental idea of  predicativity; because, 1) it allows the (of course, restricted) use of quantifiers relating to an indeterminate class of sets while the other idea does not use them; and 2) according to Kleene's lemma cited below, for each $[\alpha < \omega^{ck}_1]$, the set described by the variable $[X \in \mathrm{RA}_{\alpha}]$ can be enumerated by means of a function which is defined by a formula $D$ fulfilling the second condition of Poincar\'e. \hfill \citeyearpar[p. 378]{Kreisel1960} 
\end{quote}
}

Simplifying slightly, the wide definition of predicativity which Kreisel proposed on the basis of Poincar\'e's second characterization was that of $\Delta^1_1$-\textsl{definability}.  Recall  that a set $X \subseteq \mathbb{N}$ has this property -- conventionally abbreviated  $X \in \Delta^1_1$ --  just in case there exist formulas $\mathcal{L}^2_a$-formulas $\psi^2_1(x,X)$ and $\psi^2_2(x,X)$ not containing second-order quantifiers such that $X = \{n \in \mathbb{N} : \mathfrak{N} \models \E X \psi^2_1(\overline{n},X)\} =  \{n \in \mathbb{N} : \mathfrak{N} \models \A X \psi^2_1(\overline{n},X)\}$.  The Kreisel's wide definition can now be symbolized as follows: \begin{example}
$\A X \subseteq \mathbb{N}(Pred(X) \rightarrow X \in \Delta^1_1)$
\label{pred2}
\end{example} 

In order to understand the relation of this principle to the passages reproduced above it useful to recall that $\Delta^1_1$-definitions of sets are \textsl{absolute} between $\omega$-models of second order arithmetic.  For suppose that $\mathfrak{N}$ denotes the standard model of the second-order Peano axioms (wherein the second-order quantifiers range over all of $\mathcal{P}(\mathbb{N})$) and $\mathfrak{M}$ is a Henkin $\omega$-model of theory $\Delta^1_1$-$\mathsf{CA}_0$ consisting of first-order Peano arithmetic and comprehension for $\Delta^1_1$-formulas (wherein the second-order quantifiers range over some prescribed subset of  $\mathcal{P}(\mathbb{N})$ containing the $\Delta^1_1$-definable sets).  Then if $X$ is $\Delta^1_1$-definable -- say with $\E X \psi^2_1(x,X)$ as its $\Sigma^1_1$-characterization -- then  $X = \{n \in \mathbb{N} : \mathfrak{N} \models \E X \psi^2_1(\overline{n},X)\} =  \{n \in \mathbb{N} : \mathfrak{M} \models \E X \psi^2_1(\overline{n},X)\}$.   

The result of Kleene to which Kreisel alludes is as follows:
\begin{theorem}[\citealp{Kleene1955b}] For all $X \subseteq \mathbb{N}$, $X \in \Delta^1_1$ if and only if $X \in \mathrm{HYP}$. \label{kleenethm} \end{theorem}
\noindent A consequence is that $\Delta^1_1$-definable sets also exhibit the same absoluteness property with respect to $\omega$-models satisfying axioms (say $\mathsf{A}$) which entail that their second-order domain contains at least the class $\mathrm{HYP}$ of hyperarithmetical sets.   Taking into account both the left-to-right and right-to-left directions of this result, Kreisel writes

{\footnotesize 
\begin{quote}
This theorem clearly expresses Poincar\'e's idea. Since he considers the set of natural numbers to be well determined, we only consider $\omega$-models. But since he does not consider the totality of the sets of natural numbers to be well determined, the axioms $\mathsf{A}$ cannot distinguish between the different $\omega$-models (of $\mathsf{A}$). And a formula $[\phi(n)]$ can only be considered as an unambiguous definition if $[\phi(n)]$ defines the same set in each $\omega$-model of $\mathsf{A}$. If the set $\mathsf{A}$ itself is already (predicatively) well-defined we recover the hyper-arithmetic sets. \hfill \citeyearpar[p. 380]{Kreisel1960}
\end{quote}
}

This passage can be understood as providing justification for imposing (\ref{pred2}) as a \textsl{necessary} condition on the common concept of predicative definability which Kreisel locates in the work of Poincar\'e -- i.e. a predicative definition of a set $X$ should (at minimum) define it in a manner such that its extension is not altered by the addition of new sets from the class $\mathcal{P}(\mathbb{N})$ into the domain of second-order quantification.   But now if (\ref{pred1}) and (\ref{pred2}) are regarded as constitutive principles for  $Pred(X)$, then Theorem \ref{kleenethm} -- which in the notation we have adopted is expressed as 
\begin{example}
$\A X \subseteq \mathbb{N}(X \in \Delta^1_1 \leftrightarrow X \in \textrm{HYP} )$
\label{kleenethm2}
\end{example}
-- indeed succeeds in squeezing the intuitive concept of predicative definability between two precisely defined notions.  For in particular (\ref{pred1}), (\ref{pred2}), and (\ref{kleenethm2}) immediately yield
\begin{example}
$\A X \subseteq \mathbb{N}( Pred(X) \leftrightarrow X \in \textrm{HYP} )$
\end{example}

While again not entirely uncontroversial, such an analysis of predicativity has proven to be more influential in the subsequent development of generalized recursion theory and effective descriptive set theory.\footnote{See, e.g., \citep[\S 16]{Rogers1987}, \citep[\S II]{Sacks1990}, \citep{Moschovakis2016a}.}  But more germane to the current paper are the structural affinities between the argument just rehearsed and Kreisel's validity argument as reconstructed in \S\ref{validity}.   For to state the obvious parallel, the sorts of reflection on `traditional' accounts and mathematical practice which Kreisel gives for adopting (\ref{pred1}) and (\ref{pred2}) as constitutive principles for predicative definability can be compared to those which he gives for adopting (\ref{valarg}i) and (\ref{valarg}ii) as constitutive principles for logical validity.   Relative to this system of analogies, Kleene's Theorem \ref{kleenethm} thus plays the same role relative to Kreisel's analysis of predicativity which G\"odel's Completeness Theorem (\ref{valarg}iii) plays for his analysis of logical validity.\footnote{Although Kreisel does not note this analogy explicitly, he does say the following about a result of a theorem he credits to \citet{Addison1958} and \citet{Grzegorczyk1958}: `The $\ldots$ theorem expresses a definitional completeness of the class $[\mathrm{HYP}]$: if $X$ is defined by a property expressible [by a formula in the language of ramified analysis], then $X$ itself is definable in $[\mathrm{RA}_{\omega^{ck}_1}]$' \citeyearpar[p. 378]{Kreisel1960}.   (The result referred to states that if $\E X! \Phi(X)$ is true in $\mathfrak{N}$ then the (unique) set $A \subseteq \mathbb{N}$ satisfying this formula is hyperarithmetical and thus definable by a formula in the language of ramified analysis all of whose quantifiers are bounded by terms defining the level $\mathrm{RA}_{\alpha}$ for $\alpha < \omega^{ck}_1$.  See \citep[pp. 421-424]{Rogers1987} and note \ref{forcingnote} for additional discussion.)}

\subsubsection{Intuitionistic validity}
\label{intval}

Kreisel worked extensively on intuitionistic mathematics and formal systems from 1955 into the 1970s.   His contributions during can be broadly divided into three categories: i) in \citeyearpar{Kreisel1958,Kreisel1958f,Dyson1961,Kreisel1970} he investigated the metamathematical properties of the system of intuitionistic first-order logic known as the \textsl{Heyting Predicate Calculus} [$\mathsf{HPC}$]; ii) in (\citeyear{Kreisel1962a}, \citeyear[\S 2.3]{Kreisel1965a}) he proposed a formalization of the so-called \textsl{proof} (or \textsl{Brouwer-Heyting-Kolomogrov} [BHK]) interpretation of the intuitionistic logical connectives in the form of his so-called \textsl{Theory of Contructions}; iii) in \citeyearpar{Kreisel1959c,Kreisel1963,Kreisel1965a,Kreisel1968b}  he investigated the formalization of intuitionistic analysis, inclusive of introducing a modified realizability interpretation, continuity principles, and axioms for so-called \textsl{absolutely free} (or \textsl{lawless}) sequence .   We have discussed the bearing of iii) on Kreisel's creating subject argument in \S\ref{gmp}.   But it is also evident that his contributions to i) and ii) -- which we will see are related also to iii) -- formed a prominent part of the background behind what he says about \textsl{both} informal rigour \textsl{and} completeness proofs in \citeyearpar{Kreisel1967a,Kreisel1967b}.

A central question in this regard is the following: 
\begin{example}[(Q)]
\textsl{Can an informally rigorous argument he given that the set of sentences in the language of first-order logic which are valid when their quantifiers and connectives are assigned their intended intuitionistic interpretations coincide with the theorems of a given formal system}?
\end{example}
Much of Kreisel's work on intuitionism engages with this question.  But before it can be addressed using the model of informal rigour we have proposed in \S\ref{oninfrig}, more needs to be said about the concept of intuitionistic validity itself.   A useful summary of some of the complexities which arise in this regard by the following comments of Troelstra on Kreisel's validity argument for classical first-order logic we have considered in \S\ref{validity}:

{\footnotesize
\begin{quote}
Obviously, the extension of the class of valid sentences depends on the mathematical assumptions about possible structures. Classically, $\mathit{Val}$ is determined by the extension of 
the class of domains and relations; intuitionistically $\mathit{Val}$ not only depends on this extension, but also on the class of proofs of logically compound statements concerning intuitionistic 
structures. This dependence on the class of proofs is implicit in the axioms we postulate for certain mathematical objects; $\dots$ for example, where relations with a lawless parameter are included 
in the range of the relation  quantifiers in the definition of  $\mathit{Val}$, the extension of  $\mathit{Val}$ is determined by the axioms for lawless sequences, and the axioms express something about the possible proofs for statements of a certain form $\dots$. As a result, the dependence of $\mathit{Val}$ on mathematical assumptions is much more striking and essential in the intuitionistic case.   For example, on the more or less plausible assumption of Church's thesis, [$\mathsf{HPC}$] is incomplete. \\ \hspace*{1ex} \hfill \citep[p. 103]{Troelstra1977a}
\end{quote}
}

As this passage suggests, the characterization of the concept of \textsl{intuitionistic validity} depends on the interpretation of the intuitionistic connectives which is in turn intended to be given by the proof interpretation.   This suggests that the class of intuitionistic validities will  depend on how the notion of \textsl{constructive proof} is characterized.  But as Troelstra points out, such a characterization may itself depend on the principles which are accepted within the practice of intuitionistic mathematics.  Such a choice of principles may thus in turn have implications for metamathematical results involving intuitionistic formal systems.    

The work which Kreisel conducted under the first rubric mentioned above may be understood as illustrating this point in a particularly paradigmatic manner.  But it should also be kept in mind that the contemporaneous work he conducted under the second rubric can also be understood as an attempt to provide a characterization of the class of constructive proofs.   In particular, a central part of Kreisel's intention in proposing his Theory of Constructions was to characterize the conditions under which a construction $p$ can serve as a proof of a formula $\phi$ (notation $p:\phi$) via an inductive definition based on the logical form of $\phi$ so that the `constructive truth' (or `validity') of a formula coincides with its `constructive provability' in conformity with canonical expositions of the proof interpretation such as that given by Heyting (e.g. \citeyear{Heyting1956}).\footnote{Two important features of the Theory of Constructions are that it is `type and logic free' and that it conforms to the principle that the relation between constructive proofs and the formulas they demonstrate is \textsl{decidable} (in the sense expressed by the first of Kreisel's creating subject axioms CS1 discussed in \S\ref{gmpis}).  About the latter feature,  \cite{Dummett2000} would later remark `[I]n the explicatory clauses,  the sentential operators are applied only to decidable statements, and the quantifiers only to decidable predicates$\dots$ Hence the intuitive explanation of the logical constants may be claimed as genuine \textsl{explanations}, since, in order to understand them, it is necessary to know already, not the full meanings of the logical constants to which they relate, but only their meanings in a very restricted type of context.  The standard explanations of the intuitionistic logical constants are thus free of the circular character of the intuitive explanations of the classical ones; and, indeed,  in Kreisel's mathematical theory of constructions, these explanations appear as actual definitions' (p. 281).  In regard to this remark, it should be recalled that Kreisel originally proposed that the clause in the proof interpretation for implication which is now taken to read `a proof of $\phi \rightarrow \psi$ is a construction $p_1$ which permits us to transform any proof of $\phi$ into a proof of $\psi$' (which is akin to a $\Pi^0_2$-statement) should be replaced with `a proof of $\phi \rightarrow \psi$ is a pair $\langle p_1,p_2 \rangle$ such that $p_1$ is a construction which permits us to transform any proof of $\phi$ into a proof of $\psi$ and $p_2$ is a proof of this fact' (which is akin to a $\Delta^0_1$ -- and hence decidable -- statement).   Although this suggestion was originally subsumed within the BHK interpretation by \citet{Troelstra1977} it was later removed by \citet{Troelstra1988} (at which time the referent of `K' in the acronym was also changed from `Kreisel' to `Kolmogorov').   See \citep{Scott1970a}, \citep{Sundholm1983}, and \citep{Dean2016c} for more on the bearing of these developments on the characterization of intuitionistic validity.}   In particular, Kreisel's original proposal in \citeyearpar{Kreisel1962a} can be understood as leading to the following characterization:
\begin{example} $\mathit{iVal}(\phi)$ just in case there exists a construction $p$ such that $p:\phi$.
\label{ival}
\end{example}
We will take this as the common concept of intuitionistic validity which an informally rigorous argument should seek to characterize.   For without going into details about the Theory of Constructions itself, (\ref{ival}) at least coincides with the \textsl{form} of the proof interpretation which is widely accepted by practicing intuitionists.    

A basic challenge which stands in the way of relating the completeness of $\mathsf{HPC}$ to (Q) derives from the fact that most of the mathematical results which have been obtained about intuitionistic validity pertain not to $\mathit{iVal}$ but rather to other characterizations of this notion.   Of these, the first we must consider is defined in terms of so-called \textsl{internal interpretations} as also introduced originally by Kreisel:\footnote{The term `internal interpretation' is due to \citet[p. 156 ff]{Dummett2000}.} 

{\footnotesize
\begin{quote}
Presumably, it would be possible to introduce a so-called semantic definition of intuitionistic truth along the line of Tarski's definition of truth, only now the logical constants in the definition would have to be interpreted intuitionistically. $\dots$ It is of course clear that such a truth definition does not explain the meaning of the logical constants, but presupposes it.  \hfill \citeyearpar[p. 318]{Kreisel1958b}
\end{quote}
}
On this basis, Kreisel  proposed to adapt Tarski's definition of truth in a model to the language of $\mathsf{HPC}$ by replacing the classical notion of \textsl{set} by the corresponding intuitionistic notion of \textsl{species}.  If the latter notion is itself accepted as precise, then such a characterization can be understood as leading to an intuitionistic counterpart $iV_1$ to the predicate $V$ -- i.e. classical model theoretic validity -- as it figures in the argument of \S\ref{validity}.   In particular, \citet[p. 139] {Kreisel1962} states the following definition in the context of posing the question of the completeness of $\mathsf{HPC}$:\footnote{In \citeyearpar[p. 317]{Kreisel1958b} Kreisel refers to a less detailed formulation of this definition as the `intended interpretation' of the language of $\mathsf{HPC}$.   A similar definition is also described by \citep[\S 7]{Troelstra1977} and \citep[\S 13.1.7]{Troelstra1988a} under the name of `intuitionistic validity'.  This terminological decisions aside, it still seems reasonable to distinguish $iV_1$ as a definition which can potentially be stated precisely in an appropriate intuitionistic metatheory from $\mathit{iVal}$ which is not only informal but also intended to express the (genuinely) intended interpretation of intuitionistic validity in terms of constructive provability.}
\begin{example} $iV_1(\phi)$ just in case $\forall D \forall P^*_1 \ldots \forall P^*_k \phi_D(P^*_1,\ldots,P^*_k)$ where $\phi$ contains the predicate letters $P_1,\ldots,P_k$ or arities $r_1,\ldots,r_k$, $D$ ranges over arbitrary species, $P^*_1$ ranges over subspecies of $D^{r_i}$, and $\phi_D(P^*_1,\ldots,P^*_k)$ denotes the result of restricting the quantifiers of $\phi$ to $D$ and interpreting $P_i$ as $P^*_i$.
\label{iv1}
\end{example}

\vspace{11pt}

We will now adopt the notation $\proves_{\mathsf{HPC}}$ for formal provability in $\mathsf{HPC}$ (see, e.g., \citealp[\S 6.2]{Dalen2008}).  In his work on $\mathsf{HPC}$, Kreisel was primarily concerned not with the classical notion of completeness with respect to this relation but rather with the constructive provability of a conditional of the form `if $\phi$ is intuitionistically valid, then $\proves_{\mathsf{HPC}} \phi$'.  The notions of formal completeness which figure in the first of Kreisel's results can now be introduced as follows:
\begin{example}
\begin{enumerate}[i)]
\item Strong completeness: For all first-order formulas $\phi$, $iV_1(\phi) \ \Rightarrow \ \proves_{\mathsf{HPC}} \phi$.
\item Weak completeness: For all first-order formulas $\phi$, $iV_1(\phi) \Rightarrow  \texttt{Not Not}  \vdash_{\mathsf{HPC}} \varphi$.
\end{enumerate}
\end{example}
Here $\Rightarrow$ and $\texttt{Not}$ are respectively intended to express intuitionistic implication and negation in the metalanguage -- e.g. the weak completeness of $\mathsf{HPC}$ states that if $\phi$ is valid in the sense of (\ref{iv1}), then the assumption that $\phi$ is not provable in $\mathsf{HPC}$ can be constructively shown to lead to a contradiction.  

Part of the interest of results involving these definitions derives from the fact that that $iV_1$ can be shown to be equivalent to an alternative characterization of intuitionistic validity in terms of \textsl{spreads} (i.e. the intuitionistic analogue of an infinite tree) originally due to \citet{Beth1955,Beth1956}.  Beth had claimed that  $\mathsf{HPC}$ can be shown to be complete relative to his proposed semantics by `intuitionistic arguments'.  This claim was examined critically by G\"odel and Kreisel in \citeyearpar{Kreisel1958} and then by Dyson and Kreisel in  \citeyearpar{Dyson1961} in relation to the topological semantics which had been proposed by \citet{Mostowski1948} and which \citet{Kreisel1958f} had also investigated.   Although we will not go further into details here, these definitions are in turn connected with the principles which are assumed to govern choice sequences, in particular in regard to the existence of so-called \textsl{lawless sequences} (as observed by Troelstra in the passage above).\footnote{The notion of lawless sequence was introduced by \cite{Kreisel1958d} (independently of Brouwer) under the name `absolutely free choice sequence'  -- i.e. one which which has no general restriction on future values except finite number of initial values.  (He also later noted in \citeyearpar[p. 224]{Kreisel1968b} that one could find  a related notion -- i.e. a `restriction of restriction' -- in a remark of Brouwer.)  While it is generally agreed that this notion has not been useful in intuitionistic analysis (e.g. \citealp[p. 283]{Myhill1967}), it has proven to be of interest in logic.  Indeed, \cite{Burgess1981} states `no \textsl{general} refutation of all intuitionistic [non-theorems] is known which does not involve lawlessness'. }

The Beth and topological semantics for $\mathsf{HPC}$ can themselves be understood as giving rise to distinct precise characterizations of intuitionistic validity -- say $iV_2$ and $iV_3$ -- which can be shown to coincide extensionally with $iV_1$.\footnote{See \citep[\S13]{Troelstra1988a} for the relevant definitions and results.}
On the other hand, a sequence of metamathematical results which Kreisel was involved in obtaining can be understood as illustrating that we should not expect there to be a constructive completeness proof for $\mathsf{HPC}$ which is capable of functioning in an informally rigorous analysis of $\mathit{iVal}$ in the same way he understood G\"odel's Completeness Theorem for classical first-order logic to function in his argument for coincidence of $\mathit{Val}$ and $V$.   The first two of these can now be stated in their original form as follows:   
\begin{theorem} [G\"{o}del-Kreisel]
$\mathsf{HPC}$ is strongly complete if and only if Markov's principle holds. 
\label{Godel-Kreisel}
\end{theorem}
\begin{theorem}[Dyson-Kreisel] 
$\mathsf{HPC}$ is weakly complete if and only if weak Markov's principle holds. 
\end{theorem}
\noindent Here Markov's principle and weak Markov's principle correspond to the statements
\vspace{5pt}
\begin{examples}
\item[(MP$_0$)]  $\forall \alpha_B \neg\neg \exists n A(n, \alpha) \to \forall \alpha_B \exists n A(n, \alpha)$  
\item[(wMP$_0$)]  $\forall \alpha_B \neg\neg \exists n A(n, \alpha) \to \neg\neg \forall \alpha_B \exists n A(n, \alpha)$ 
\end{examples}
where $\alpha$ is a variable over choice sequences, $A (n, \alpha)$ is a primitive recursive relation, and the subscript $B$ means that the $\alpha$ are chosen from the full binary spread $B$ with values $0$ and $1$.

As noted by \citet[p. 6]{Moschovakis2019}, MP$_0$ implies GMP over $\mathsf{FIM}_0$ as formulated in \S\ref{gmp}.  In light of this,  Kreisel's creating subject argument -- which we have reconstructed as an informally rigorous \textsl{refutation} of GMP -- is also potentially related to work on the completeness of $\mathsf{HPC}$.  One way of formulating this is to observe that if we were to work in an intuitionistic metatheory which is sufficient to establish both establish Theorem \ref{Godel-Kreisel} and carry out Kreisel's proof of Theorem \ref{neggmp} we could then obtain the following consequence:
\begin{corollary}\label{non-comp}
\textnormal{\texttt{Not for all}} $\phi , iV_1(\phi) \Rightarrow \ \vdash_{\mathsf{HPC}} \phi$.
\end{corollary}

Although it seems likely that Kreisel was aware of this connection at the time of \citeyearpar{Kreisel1967b}, he was also cautious in his use of the specific term `incompleteness' in regard to $\mathsf{HPC}$.    For instance in remarking on a refutation of GMP in \citeyearpar[p. 182]{Kreisel1967b} based on an alternative definition of choice sequences Kreisel notes that `this is not yet sufficient to establish incompleteness of Heyting's predicate calculus \textsl{without some assumption} \textsl{on constructive functions} (such as Church's thesis)'.   Here `Church's thesis' should be understood not in the more familiar sense discussed in \S\ref{words} but rather to to refer one of a class of formal principles about choice sequences whose consistency with intuitionistic analysis Kreisel was involved in investigating.  One example is the statement
\begin{example}[(CT)] $\forall \alpha \exists x \forall y (\E z(T(x,y,z) \mand U(z) = \alpha(y))$ \end{example}
which can be understood as stating that all functions of type $\mathbb{N} \rightarrow \mathbb{N}$ are recursive.\footnote{Here  $T(x,y,z)$ is Kleene's predicate formalizing `the Turing machine with index $y$ halts with computation sequence $z$' and $U(z)$ is a primitive recursive function whose definition formalizes `the output encoded by $z$ is $y$'.  See, e.g., \citep[\S 4.3]{Troelstra1988}.}

A characteristic results which \citet[pp. 133-134]{Kreisel1970} obtained about CT is the following:
 
\begin{theorem}\label{incomp}
\textnormal{CT} implies that the set $($or species$)$ of intuitionistic validities -- i.e. $I_ 1 =_{\mathrm{df}} \{\ulcorner \varphi \urcorner : iV_1(\phi) \}$ -- is not recursively enumerable.   
\end{theorem}

\noindent There is indeed a straightforward sense in which Theorem \ref{incomp} bears on the prospect of formulating an informally rigorous analysis of  intuitionistic validity on the model of his classical validity argument (as presented in \S \ref{validity}).   For on the other hand, it is generally accepted that the axioms of $\mathsf{HPC}$ are intuitionistically valid in the specific sense of $\mathit{iVal}$.   On the other hand, it is a commonplace of both classical and intuitionistic mathematics that the set of theorems of a recursively axiomatized system such as $\mathsf{HPC}$ is recursively enumerable.  It would thus seem that Theorem \ref{incomp} presents a substantial obstacle to providing an informally rigorous analysis of $\mathit{iVal}$ on the model of the classical validity argument.   About this situation Kreisel remarked as follows:

{\footnotesize
\begin{quote}
As observed in \citeyearpar[p. 140]{Kreisel1962} Church's thesis implies incompleteness of Heyting predicate calculus, and this is sharpened [by the proof of Theorem \ref{incomp}]  $\ldots$ Regarding Church's thesis as neither plausible nor refuted, we can say that the notion of \emph{constructive validity of first order formulae} depends on problematic properties of the basic notion of constructive function; like second order validity \citeyearpar[p. 157]{Kreisel1967b}, and unlike first order validity $\ldots$ in the classical case. \hfill \citeyearpar[p. 126]{Kreisel1970}
\end{quote}}

This passage brings together all three of the arguments we have considered in \S\ref{casestudies} in regard to Kreisel's engagement with informal rigour -- i.e. his validity argument (which affirms that the set of classical first-order validities coincides with the classically derivable ones and hence \textsl{is} r.e.), his CH argument (which suggests that the class of second-order classical validities \textsl{is not} r.e.) and his creating subject argument (which bears on the characterizations of constructive functions and thus, via Theorem \ref{incomp}, whether the of set intuitionistic of validities is r.e.).   Several of his other remarks in \citeyearpar{Kreisel1970} make clear that Kreisel was well aware of these connections.  But rather than pursuing his further remarks here, we will conclude with the following brief reflections on how subsequent developments bear on our framing question (Q).  

Theorem \ref{incomp} is illustrative of a class of results which show how several different precise analyses of intuitionistic validity lead to sets of sentences which are \textsl{complex} in the sense now studied within computability theory.  For instance a well-known theorem of \citet{McCarty1988} states that if in addition to CT, $\mathrm{MP}_0$ is also assumed then $I_1$ is not even \textsl{arithmetically definable}.  And if we similarly define $iV_4(\phi)$ to mean that $\phi$ is \textsl{realizable} in the original sense of \citet{Kleene1945}, then \citet{Plisko1984} showed that the corresponding set of G\"odel numbers  $I_4$ is $\Pi^1_1$-complete.   On the other hand, if intuitionistic validity is analyzed in terms of Kripke models ($iV_5$) and the latter are understood from the perspective of a classical metatheory then $\mathsf{HPC}$ \textsl{can} be shown to be complete in the conventional sense.  It is then not difficult to show that the corresponding set of G\"odel numbers $I_5$ is r.e. but not recursive.\footnote{See, respectively, \citep[2.6.6]{Troelstra1988} and \citep[8.8.2]{Sorensen2006}.} 

Such results may be understood as engaging with informal rigour in at least four ways.   First,  they testify to the diversity of precise analyses of the common concept of intuitionistic validity which we have suggested stands behind the fundamental characterization (\ref{ival}).   Second, they testify to the dependence of the results which can be obtained on whether a classical or constructive metatheory is adopted to reason about these notions.   Third, they illustrate the apparent willingness of theorists to employ definitions and results from computability theory to interpret the significance of the results in both the classical and constructive settings.\footnote{If we adopt Kreisel's original perspective on Church's Thesis described in \S\ref{words}, then the process of drawing conclusions from the results discussed in the prior paragraph should itself be understood as an application of informal rigour.  On the other hand, the situation is complicated by the fact that formulas in the language of arithmetic or analysis need not have intuitionistically equivalent prenex normal forms (see, e.g., \citealp{Burr2004}).  This in turn suggests that additional care is required when we attempt to read off conclusions from classifications of sets or predicates relative to the classical formulations of the arithmetical and analytical hierarchies.}  Fourth -- and perhaps most significantly -- they illustrate what Kreisel appears to have regarded as a \textsl{divergence} or \textsl{bifurcation} of different potential analyses of intuitionistic validity which stands in contrast to the \textsl{convergence} of precise analyses of classical validity which he understood to be effected by the argument of \S\ref{validity}.   

These considerations illustrate from yet another perspective how Kreisel's method of informal rigour  stands at the crossroads of many developments within mathematical logic and the foundations of mathematics.   In the particular case of intuitionistic validity, they also suggest that the application of the method may not yet have reached its full potential -- e.g. in regard to its interaction with computability theory with both formalized and informal metatheory.   But this is a possibility we happily defer to another occasion.


{\footnotesize
\begin{thebibliography}{211}
\providecommand{\natexlab}[1]{#1}
\providecommand{\url}[1]{\texttt{#1}}
\expandafter\ifx\csname urlstyle\endcsname\relax
  \providecommand{\doi}[1]{doi: #1}\else
  \providecommand{\doi}{doi: \begingroup \urlstyle{rm}\Url}\fi

\bibitem[Addison(1958)]{Addison1958}
J.~Addison.
\newblock Separation principles in the hierarchies of classical and effective
  descriptive set theory.
\newblock \emph{Fundamenta mathematicae}, 46\penalty0 (2):\penalty0 123--135,
  1958.

\bibitem[Bell(1977)]{Bell1977a}
J.~Bell.
\newblock \emph{Boolean-{V}alued {M}odels and {I}ndependence {P}roofs in {S}et
  {T}heory}.
\newblock Clarendon Press, Oxford, 1977.

\bibitem[Beth(1955)]{Beth1955}
E.~Beth.
\newblock Semantic entailment and formal derivability.
\newblock \emph{Mededelingen der Koninklijke Nederlandse Akademie van
  Wetenschappen, Afd. Letterkunde}, 18\penalty0 (13):\penalty0 309--342, 1955.

\bibitem[Beth(1956)]{Beth1956}
E.~Beth.
\newblock Semantic construction of intuitionistic logic.
\newblock \emph{{Mededelingen der Koninklijke Nederandse Akademie van
  Wetenschappen, Afd. Letterkunde}}, 19\penalty0 (11):\penalty0 357--388, 1956.

\bibitem[Blanchette(1996)]{Blanchette1996}
P.~Blanchette.
\newblock {Frege and Hilbert on consistency}.
\newblock \emph{The Journal of Philosophy}, 93\penalty0 (7):\penalty0 317--336,
  1996.

\bibitem[Boolos(1971)]{Boolos1971}
G.~Boolos.
\newblock The {I}terative {C}onception of a {S}et.
\newblock \emph{The Journal of Philosophy}, 68:\penalty0 215--232, 1971.

\bibitem[Boolos(1985)]{Boolos1985b}
G.~Boolos.
\newblock Nominalist platonism.
\newblock \emph{The Philosophical Review}, 94\penalty0 (3):\penalty0 327--344,
  1985.

\bibitem[Brouwer(1948)]{Brouwer1948a}
L.~E. Brouwer.
\newblock Essentieel-negatieve eigenschappen.
\newblock \emph{KNAW Proceedings}, pages 963--964, 1948.
\newblock {Reprinted as ``Essentially negative propreties'' in
  \cite{Brouwer1976}, pp. 478-479.}

\bibitem[Burgess(1981)]{Burgess1981}
J.~P. Burgess.
\newblock The completeness of intuitionistic propositional calculus for its
  intended interpretation.
\newblock \emph{Notre Dame Journal of Formal Logic}, 22\penalty0 (1):\penalty0
  17--28, 1981.

\bibitem[Burgess(1985)]{Burgess1985a}
J.~P. Burgess.
\newblock {R}eview of \cite{Boolos1971}, \cite{Parsons1977}, \cite{Wang1977}
  and others.
\newblock \emph{The Journal of Symbolic Logic}, 50\penalty0 (2):\penalty0
  544--547, 1985.

\bibitem[Burr(2004)]{Burr2004}
W.~Burr.
\newblock The intuitionistic arithmetical hierarchy.
\newblock In \emph{Logic Colloquium}, volume~99, pages 51--59, 2004.

\bibitem[Button and Walsh(2018)]{Button2018}
T.~Button and S.~Walsh.
\newblock \emph{Philosophy and Model Theory (with a historical appendix by
  Wilfrid Hodges)}.
\newblock Oxford Univeristy Press, Oxford, 2018.

\bibitem[Cappelen(2018)]{Cappelen2018}
H.~Cappelen.
\newblock \emph{Fixing language: An essay on conceptual engineering}.
\newblock Oxford University Press, Oxford, 2018.

\bibitem[Carnap(1950)]{Carnap1950}
R.~Carnap.
\newblock Empiricism, semantics, and ontology.
\newblock \emph{Revue Internationale de Philosophie}, 4\penalty0 (11), 1950.
\newblock Reprinted in the Supplement to \textsl{Meaning and Necessity: A Study
  in Semantics and Modal Logic}, enlarged edition (University of Chicago Press,
  1956).

\bibitem[Carnap(1937)]{Carnap1937}
R.~Carnap.
\newblock \emph{Logical syntax of language (Translated by Amethe Smeaton)}.
\newblock Routledge, London, 1937.

\bibitem[Carnap(1945)]{Carnap1945}
R.~Carnap.
\newblock Two {C}oncepts of {P}robability.
\newblock \emph{Philosophy and Phenomenological Research}, 5:\penalty0
  513--532, 1945.

\bibitem[Carnap(1956)]{Carnap1956}
R.~Carnap.
\newblock \emph{Meaning and necessity: a study in semantics and modal logic},
  volume~30.
\newblock University of Chicago Press, 1956.

\bibitem[Carnap(1962)]{Carnap1962}
R.~Carnap.
\newblock \emph{Logical Foundations of Probability}.
\newblock University of Chicago Press, Chicago, 1962.

\bibitem[Chuaqui(1972)]{Chuaqui1972}
R.~Chuaqui.
\newblock Forcing for the impredicative theory of classes.
\newblock \emph{The Journal of Symbolic Logic}, 37\penalty0 (1):\penalty0
  1--18, 1972.

\bibitem[Cohen(1963)]{Cohen1963a}
P.~Cohen.
\newblock The independence of the continuum hypothesis.
\newblock \emph{Proceedings of the National Academy of Sciences}, 50\penalty0
  (6):\penalty0 1143--1148, 1963.

\bibitem[Cohen(1971)]{Cohen1971}
P.~Cohen.
\newblock Comments on the {F}oundations of {S}et {T}heory.
\newblock In \emph{Axiomatic Set Theory}, volume 13, Part I of
  \emph{Proceedings of Symposia in Pure Mathematics}, pages 9--15. American
  Mathematical Society, Providence, R.I., 1971.

\bibitem[Cohen(2002)]{Cohen2002}
P.~Cohen.
\newblock The discovery of forcing.
\newblock \emph{Rocky Mountain Journal of Mathematics}, 32\penalty0 (4), 2002.

\bibitem[Cohen and Hersh(1967)]{Cohen1967}
P.~Cohen and R.~Hersh.
\newblock {Non-Cantorian set theory}.
\newblock \emph{Scientific American}, 217\penalty0 (104-116):\penalty0 105,
  1967.

\bibitem[Dean(2016)]{Dean2016b}
W.~Dean.
\newblock Squeezing feasibility.
\newblock In A.~Beckmann, L.~Bienvenu, and N.~Jonoska, editors, \emph{Pursuit
  of the Universal: 12th Conference on Computability in Europe, CiE 2016,
  Paris, France, June 27 - July 1, 2016, Proceedings}, pages 78--88, 2016.

\bibitem[Dean(2020)]{Dean2020e}
W.~Dean.
\newblock Incompleteness via paradox and completeness.
\newblock \emph{The Review of Symbolic Logic}, 13\penalty0 (2):\penalty0
  541--592, 2020.

\bibitem[Dean and Walsh(2017)]{Dean2017b}
W.~Dean and S.~Walsh.
\newblock The prehistory of the subsystems of second-order arithmetic.
\newblock \emph{The Review of Symbolic Logic}, 10\penalty0 (2):\penalty0
  357--396, 2017.

\bibitem[Dean(2015)]{Dean2015a}
W.~Dean.
\newblock Arithmetical reflection and the provability of soundness.
\newblock \emph{Philosophia Mathematica}, 23\penalty0 (1):\penalty0 31--64,
  2015.

\bibitem[Dean and Kurokawa(2016)]{Dean2016c}
W.~Dean and H.~Kurokawa.
\newblock {Kreisel's Theory of Constructions, the Kreisel-Goodman paradox, and
  the second clause}.
\newblock In T.~Piecha and P.~Schroeder-Heister, editors, \emph{Advances in
  Proof-Theoretic Semantics}, Trends in Logic, pages 27--63. Springer, Berlin,
  2016.

\bibitem[Dragalin(1988)]{Dragalin1988}
A.~G. Dragalin.
\newblock \emph{Mathematical Intuitionism}.
\newblock American Mathematical Society, Providence, R.I., 1988.
\newblock Translation of the {R}ussian original from 1979.

\bibitem[Drake(1974)]{Drake1974}
F.~Drake.
\newblock \emph{Set {T}heory: {A}n {I}ntroduction to {L}arge {C}ardinals},
  volume~76 of \emph{Studies in Logic and the Foundations of Mathematics}.
\newblock North-Holland, Amsterdam, 1974.

\bibitem[Dummett(2000)]{Dummett2000}
M.~Dummett.
\newblock \emph{{Elements of Intuitionism}}.
\newblock Oxford Logic Guides. Oxford University Press, Oxford, 2000.

\bibitem[Dyson and Kreisel(1961)]{Dyson1961}
V.~Dyson and G.~Kreisel.
\newblock {Analysis of Beth's semantic construction of intuitionistic logic}.
\newblock Technical Report DA-04-200-ORD-997, {Applied mathematics and
  statistical laboratories, Stanford University}, 1961.

\bibitem[Etchemendy(1990)]{Etchemendy1990}
J.~Etchemendy.
\newblock \emph{The {C}oncept of {L}ogical {C}onsequence}.
\newblock Harvard University Press, Cambridge, Mass., 1990.

\bibitem[Etchemendy(2008)]{Etchemendy2008}
J.~Etchemendy.
\newblock Reflections on consequence.
\newblock In D.~Patterson, editor, \emph{New essays on Tarski and philosophy},
  pages 263--299. Oxford: Oxford University Press, 2008.

\bibitem[Ewald and Sieg(2013)]{Hilbert2013}
W.~Ewald and W.~Sieg, editors.
\newblock \emph{{David Hilbert's Lectures on the Foundations of Logic and
  Arithmetic 1917 -- 1933}}.
\newblock Springer, Berlin, 2013.

\bibitem[Ewald(1996)]{Ewald1996}
W.~Ewald.
\newblock \emph{From {K}ant to {H}ilbert: {A} {S}ource {B}ook in the
  {F}oundations of {M}athematics.}
\newblock Oxford University Press, New York, 1996.

\bibitem[Feferman(1964a)]{Feferman1964a}
S.~Feferman.
\newblock Some applications of the notions of forcing and generic sets.
\newblock \emph{Fundamenta mathematicae}, 56\penalty0 (325):\penalty0 45, 1964.

\bibitem[Feferman et~al.(1990)]{Godel1990}
S.~Feferman et~al., editors.
\newblock \emph{{Kurt G\"odel Collected works. {V}ol. {II}. Publications
  1938--1974}}.
\newblock Oxford Univeristy Press, Oxford, 1990.

\bibitem[Feferman(1958)]{Feferman1958}
S.~Feferman.
\newblock Arithmetically definable models of formalized arithmetic.
\newblock \emph{Notices of the American Mathematical Society}, 5:\penalty0
  679--680, 1958.

\bibitem[Feferman(1962)]{Feferman1962}
S.~Feferman.
\newblock Transfinite recursive progressions of axiomatic theories.
\newblock \emph{The Journal of Symbolic Logic}, 27\penalty0 (3):\penalty0
  259--316, 1962.

\bibitem[Feferman(1964)]{Feferman1964}
S.~Feferman.
\newblock Systems of predicative analysis.
\newblock \emph{The Journal of Symbolic Logic}, 29:\penalty0 1--30, 1964.

\bibitem[Feferman(1995)]{Feferman1995a}
S.~Feferman.
\newblock {Turing in the land of $O(z)$}.
\newblock In R.~Herken, editor, \emph{The Universal Turing Machine A
  Half-Century Survey}, pages 103--134. Springer, Berlin, 1995.

\bibitem[Feferman(2005)]{Feferman2005}
S.~Feferman.
\newblock Predicativity.
\newblock In S.~Shapiro, editor, \emph{The Oxford Handbook of Philosophy of
  Mathematics and Logic}, pages 590--624. Oxford University Press, Oxford,
  2005.

\bibitem[Field(2015)]{Field2015a}
H.~Field.
\newblock What is logical validity?
\newblock In C.~R. Caret and O.~T. Hjortland, editors, \emph{Foundations of
  Logical Consequence}, pages 33--70. Oxford University Press, Oxford, 2015.

\bibitem[Field(1980)]{Field1980}
H.~Field.
\newblock \emph{Science without {N}umbers: {A} {D}efence of {N}ominalism}.
\newblock Princeton University Press, Princeton, 1980.

\bibitem[Field(1984)]{Field1984}
H.~Field.
\newblock Is mathematical knowledge just logical knowledge?
\newblock \emph{Philosophical Review}, 93:\penalty0 509--552, 1984.

\bibitem[Field(1989)]{Field1989}
H.~Field.
\newblock \emph{Realism, {M}athematics and {M}odality}.
\newblock Blackwell, New York, 1989.

\bibitem[Field(1991)]{Field1991}
H.~Field.
\newblock Metalogic and {M}odality.
\newblock \emph{Philosophical Studies}, 62\penalty0 (1):\penalty0 1--22, 1991.

\bibitem[Field(2008)]{Field2008}
H.~Field.
\newblock \emph{Saving Truth from Paradox}.
\newblock Oxford University Press, Oxford, 2008.

\bibitem[Franz\'en(2004)]{Franzen2004}
T.~Franz\'en.
\newblock \emph{Inexhaustibility: {A} Non-Exhaustive Treatment}.
\newblock Lecture Notes in Logic. A.K. Peters, Wellesley, Mass., 2004.

\bibitem[Frege(1879)]{Frege1879}
G.~Frege.
\newblock \emph{Begriffsschrift, eine der arithmetischen nachgebildete
  Formelsprache des reinen Denkens}.
\newblock Nebert, Halle, 1879.

\bibitem[Freudenthal(1976)]{Brouwer1976}
H.~Freudenthal, editor.
\newblock \emph{Brouwer, L.E.J. Collected Works: Geometry, Analysis, Topology
  and Mechanics}, volume~II.
\newblock Elsevier, 1976.

\bibitem[Friedman(1975)]{Friedman1975a}
H.~Friedman.
\newblock Some systems of second-order arithmetic and their use.
\newblock In \emph{Proceedings of the International Congress of Mathematicians,
  Vancouver 1974}, volume~1, pages 235--242. Canadian Mathematical Congress,
  1975.

\bibitem[Fujimoto(2012)]{Fujimoto2012}
K.~Fujimoto.
\newblock Classes and truths in set theory.
\newblock \emph{Annals of Pure and Applied Logic}, 163\penalty0 (11):\penalty0
  1484--1523, 2012.

\bibitem[G{\"o}del(1929)]{Godel1929}
K.~G{\"o}del.
\newblock On the completeness of the calculus of logic.
\newblock In \emph{{Collected Works Volume I}}, pages 44--123. 1929.

\bibitem[G{\"o}del(1964)]{Godel1964}
K.~G{\"o}del.
\newblock {What is Cantor's continuum problem?}
\newblock In \emph{Collected works. {V}ol. {II}.}, pages 254--270. 1964.

\bibitem[G{\"o}del(1990)]{Godel1938}
K.~G{\"o}del.
\newblock The consistency of the axiom of choice and of the generalized
  continuum-hypothesis.
\newblock In  \citet{Godel1990}, pages 1--101.

\bibitem[G{\"o}del(1946)]{Godel1946}
K.~G{\"o}del.
\newblock {Remarks before the Princeton bicentennial conference on problems in
  mathematics}.
\newblock Reprinted in \cite{Godel1990}, pp. 144-153, 1946.

\bibitem[Grzegorczyk et~al.(1958)Grzegorczyk, Mostowski, and
  Ryll-Nardzewski]{Grzegorczyk1958}
A.~Grzegorczyk, A.~Mostowski, and C.~Ryll-Nardzewski.
\newblock The classical and the $\omega$-complete arithmetic.
\newblock \emph{The Journal of Symbolic Logic}, 23\penalty0 (2):\penalty0 pp.
  188--206, 1958.

\bibitem[H{\'a}jek and Pudl{\'a}k(1998)]{Hajek1998}
P.~H{\'a}jek and P.~Pudl{\'a}k.
\newblock \emph{Metamathematics of First-Order Arithmetic}.
\newblock Springer, Berlin, 1998.
\newblock First edition 1993.

\bibitem[Halbach(2020{\natexlab{a}})]{Halbach2020a}
V.~Halbach.
\newblock The substitutional analysis of logical consequence.
\newblock \emph{No{\^u}s}, 54\penalty0 (2):\penalty0 431--450,
  2020{\natexlab{a}}.

\bibitem[Halbach(2020{\natexlab{b}})]{Halbach2020b}
V.~Halbach.
\newblock Formal notes on the substitutional analysis of logical consequence.
\newblock \emph{Notre Dame J. Formal Logic}, 61\penalty0 (2):\penalty0
  317--339, 2020{\natexlab{b}}.

\bibitem[Halimi(2017)]{Halimi2017}
B.~Halimi.
\newblock Models as universes.
\newblock \emph{Notre Dame Journal of Formal Logic}, 58\penalty0 (1):\penalty0
  47--78, 2017.

\bibitem[Hallett(2011)]{Hallett2011}
M.~Hallett.
\newblock {Absoluteness and the Skolem Paradox}.
\newblock In \emph{Logic, Mathematics, Philosophy, Vintage Enthusiasms}, pages
  189--218. Springer, Berlin, 2011.

\bibitem[Hanson(1997)]{Hanson1997}
W.~H. Hanson.
\newblock The concept of logical consequence.
\newblock \emph{Philosophical Review}, 106\penalty0 (3):\penalty0 365--409,
  1997.

\bibitem[Henkin(1950)]{Henkin1950}
L.~Henkin.
\newblock Completeness in the {T}heory of {T}ypes.
\newblock \emph{The Journal of Symbolic Logic}, 15:\penalty0 81--91, 1950.

\bibitem[Heyting(1956)]{Heyting1956}
A.~Heyting.
\newblock \emph{Intuitionism. {A}n introduction}.
\newblock North-Holland, Amsterdam, 1956.

\bibitem[Hilbert(1929)]{Hilbert1929}
D.~Hilbert.
\newblock {Problcme der Grundlagen der Mathematik}.
\newblock \emph{Mathematische Annalen}, 102:\penalty0 1--9, 1929.
\newblock English translation as ``Problems of the Grounding of Mathematics''
  in \cite{Mancosu1998}, pp. 223-233.

\bibitem[Hilbert and Ackermann(1928)]{Hilbert1928}
D.~Hilbert and W.~Ackermann.
\newblock \emph{Grundz{\"u}ge der theoretischen Logik}.
\newblock Springer, Berlin, first edition, 1928.
\newblock Reprinted in \cite{Hilbert2013}.

\bibitem[Hilbert and Bernays(1934)]{Hilbert1934}
D.~Hilbert and P.~Bernays.
\newblock \emph{{Grundlagen der Mathematik}}, volume~I.
\newblock Springer, Berlin, 1934.
\newblock Second edition 1968.

\bibitem[Hilbert and Bernays(1939)]{Hilbert1939}
D.~Hilbert and P.~Bernays.
\newblock \emph{{Grundlagen der Mathematik}}, volume~II.
\newblock Springer, Berlin, 1939.
\newblock Second edition 1970.

\bibitem[Hilbert(1899)]{Hilbert1899}
D.~Hilbert.
\newblock \emph{{Grundlagen der Geometrie}}.
\newblock Teubner, Leipzig, 1899.

\bibitem[Hilbert and Ackermann(1938)]{Hilbert1938}
D.~Hilbert and W.~Ackermann.
\newblock \emph{Grundz{\"u}ge der theoretischen Logik}.
\newblock Springer, second edition, 1938.
\newblock Translated as \cite{Hilbert1950}.

\bibitem[Hilbert and Ackermann(1950)]{Hilbert1950}
D.~Hilbert and W.~Ackermann.
\newblock \emph{Principles of Mathematical Logic}.
\newblock Chelsea Publishing Company, New York, 1950.

\bibitem[Isaacson(2011)]{Isaacson2011a}
D.~Isaacson.
\newblock The {R}eality of {M}athematics and the {C}ase of {S}et {T}heory.
\newblock In Z.~Novak and A.~Simonyi, editors, \emph{Truth, Reference, and
  Realism}, pages 1--75. Central European University Press, Budapest, 2011.

\bibitem[Jech(2003)]{Jech2003}
T.~Jech.
\newblock \emph{Set {T}heory}.
\newblock Springer Monographs in Mathematics. Springer, Berlin, 2003.

\bibitem[Kanamori(2008)]{Kanamori2008}
A.~Kanamori.
\newblock Cohen and {S}et {T}heory.
\newblock \emph{Bulletin of Symbolic Logic}, 14\penalty0 (3):\penalty0
  351--378, 2008.

\bibitem[Kaye(1991)]{Kaye1991}
R.~Kaye.
\newblock \emph{Models of {P}eano {A}rithmetic}, volume~15 of \emph{Oxford
  Logic Guides}.
\newblock Oxford University Press, Oxford, 1991.

\bibitem[Kennedy and V{\"a}{\"a}n{\"a}nen(2017)]{Kennedy2017}
J.~Kennedy and J.~V{\"a}{\"a}n{\"a}nen.
\newblock Squeezing arguments and strong logics.
\newblock In \emph{15th International Congress of Logic, Methodology and
  Philosophy of Science}. College Publications, 2017.

\bibitem[Kleene(1952)]{Kleene1952}
S.~Kleene.
\newblock \emph{Introduction to {M}etamathematics}.
\newblock North-Holland, Amsterdam, 1952.

\bibitem[Kleene(1945)]{Kleene1945}
S.~C. Kleene.
\newblock On the interpretation of intuitionistic number theory.
\newblock \emph{The Journal of Symbolic Logic}, 10\penalty0 (4):\penalty0 pp.
  109--124, 1945.

\bibitem[Kleene(1955{\natexlab{a}})]{Kleene1955a}
S.~C. Kleene.
\newblock Arithmetical predicates and function quantifiers.
\newblock \emph{Transactions of the American Mathematical Society}, 79\penalty0
  (2):\penalty0 312--340, 1955{\natexlab{a}}.

\bibitem[Kleene(1955{\natexlab{b}})]{Kleene1955b}
S.~C. Kleene.
\newblock Hierarchies of number-theoretic predicates.
\newblock \emph{Bulletin of the American Mathematical Society}, 61\penalty0
  (3):\penalty0 193--213, 1955{\natexlab{b}}.

\bibitem[Kleene and Vesley(1965)]{Kleene1965}
S.~C. Kleene and R.~E. Vesley.
\newblock \emph{The Foundations of Intuitionistic Mathematics}.
\newblock North-Holland, Amsterdam, 1965.

\bibitem[Kleene(1959)]{Kleene1959}
S.~C. Kleene.
\newblock Quantification of number-theoretic functions.
\newblock \emph{Compositio Mathematica}, 14:\penalty0 23--40, 1959.

\bibitem[Kreisel(1950)]{Kreisel1950}
G.~Kreisel.
\newblock {Note on arithmetic models for consistent formulae of the predicate
  calculus}.
\newblock \emph{Fundamenta mathematicae}, 37:\penalty0 265--285, 1950.

\bibitem[Kreisel(1951)]{Kreisel1951}
G.~Kreisel.
\newblock On the interpretation of non-finitist proofs--{Part I}.
\newblock \emph{Journal of Symbolic Logic}, 16\penalty0 (4):\penalty0 241--267,
  1951.

\bibitem[Kreisel(1952)]{Kreisel1952c}
G.~Kreisel.
\newblock On the concepts of completeness and interpretation of formal systems.
\newblock \emph{Fundamenta mathematicae}, 39:\penalty0 103--127, 1952.

\bibitem[Kreisel(1953)]{Kreisel1953}
G.~Kreisel.
\newblock {Note on arithmetic models for consistent formulae of the predicate
  calculus. II}.
\newblock In \emph{Actes du XIeme Congres International de Philosophie}, volume
  XIV, pages 39--49, Amsterdam, 1953. North-Holland.

\bibitem[Kreisel(1954)]{Kreisel1954}
G.~Kreisel.
\newblock Remark on complete interpretations by models.
\newblock \emph{Archive for Mathematical Logic}, 2\penalty0 (1):\penalty0 4--9,
  1954.

\bibitem[Kreisel(1955)]{Kreisel1955a}
G.~Kreisel.
\newblock Models, translations and interpretations.
\newblock In T.~Skolem, editor, \emph{{Mathematical interpretation of formal
  systems}}, pages 26--50. North Holland, Amsterdam, 1955.

\bibitem[Kreisel(1956)]{Kreisel1956}
G.~Kreisel.
\newblock Some uses of metamathematics.
\newblock \emph{The British Journal for the Philosophy of Science}, 7\penalty0
  (26):\penalty0 161--173, 1956.

\bibitem[Kreisel(1958{\natexlab{a}})]{Kreisel1958}
G.~Kreisel.
\newblock Mathematical significance of consistency proofs.
\newblock \emph{The Journal of Symbolic Logic}, 23\penalty0 (2):\penalty0
  155--182, 1958{\natexlab{a}}.

\bibitem[Kreisel(1958{\natexlab{b}})]{Kreisel1958a}
G.~Kreisel.
\newblock {Hilbert's programme}.
\newblock \emph{Dialectica}, 12\penalty0 (3-4):\penalty0 346--372,
  1958{\natexlab{b}}.

\bibitem[Kreisel(1958{\natexlab{c}})]{Kreisel1958b}
G.~Kreisel.
\newblock Elementary completeness properties of intuitionistic logic with a
  note on negations of prenex formulae.
\newblock \emph{The Journal of symbolic logic}, 23\penalty0 (3):\penalty0
  317--330, 1958{\natexlab{c}}.

\bibitem[Kreisel(1958{\natexlab{d}})]{Kreisel1958c}
G.~Kreisel.
\newblock {The non-derivability of $\neg (x){A}(x) \rightarrow (\E x) \neg
  {A}(x)$, ${A}$ primitive recursive, in intuitionistic formal systems}.
\newblock \emph{Journal of Symbolic Logic}, 23:\penalty0 567--457,
  1958{\natexlab{d}}.

\bibitem[Kreisel(1958{\natexlab{e}})]{Kreisel1958d}
G.~Kreisel.
\newblock Wittgenstein's remarks on the foundations of mathematics.
\newblock \emph{British Journal for the Philosophy of Science}, 9\penalty0
  (34):\penalty0 135--158, 1958{\natexlab{e}}.

\bibitem[Kreisel(1958{\natexlab{f}})]{Kreisel1958f}
G.~Kreisel.
\newblock A remark on free choice sequences and the topological completeness
  proofs.
\newblock \emph{Journal of Symbolic Logic}, pages 369--388, 1958{\natexlab{f}}.

\bibitem[Kreisel(1958{\natexlab{g}})]{Kreisel1958k}
G.~Kreisel.
\newblock A remark on free choice sequences and the topological completeness
  proofs.
\newblock \emph{The Journal of Symbolic Logic}, 23\penalty0 (4):\penalty0
  369--388, 1958{\natexlab{g}}.

\bibitem[Kreisel(1959{\natexlab{a}})]{Kreisel1959}
G.~Kreisel.
\newblock {Analysis of Cantor-Bendixson theorem by means of the analytic
  hierarchy}.
\newblock \emph{Bulletin de lAcademie Polonaise des Sciences. Serie des
  Sciences Mathematiques, Astronomiques et Physiques}, 7:\penalty0 621--626,
  1959{\natexlab{a}}.

\bibitem[Kreisel(1959{\natexlab{b}})]{Kreisel1959c}
G.~Kreisel.
\newblock Interpretation of analysis by means of constructive functionals of
  finite types.
\newblock In A.~Heyting, editor, \emph{Constructivity in Mathematic}, pages
  101--128. North Holland, Amsterdam, 1959{\natexlab{b}}.

\bibitem[Kreisel(1960{\natexlab{a}})]{Kreisel1960}
G.~Kreisel.
\newblock La pr\'edicativit\'e.
\newblock \emph{Bulletin de la Soci\'et\'e Math\'ematique de France},
  88:\penalty0 371--391, 1960{\natexlab{a}}.

\bibitem[Kreisel(1960{\natexlab{b}})]{Kreisel1960e}
G.~Kreisel.
\newblock Ordinal logics and the characterization of informal concepts of
  proof.
\newblock In J.~A. Todd, editor, \emph{Proceedings International Congress of
  Mathematicians (1958, Edinburgh)}, pages 14--21. Cambridge University Press,
  1960{\natexlab{b}}.

\bibitem[Kreisel(1961)]{Kreisel1961a}
G.~Kreisel.
\newblock Set theoretic problems suggested by the notion of potential infinity.
\newblock In \emph{{Infinitistic Methods: Proceedings of the Symposium on
  Foundations of Mathematics in Warsaw}}, pages 103--140. Pergamon Press, 1961.

\bibitem[Kreisel(1962{\natexlab{a}})]{Kreisel1962}
G.~Kreisel.
\newblock On weak completeness of intuitionistic predicate logic.
\newblock \emph{Journal of Symbolic Logic}, 27:\penalty0 139--158,
  1962{\natexlab{a}}.

\bibitem[Kreisel(1962{\natexlab{b}})]{Kreisel1962a}
G.~Kreisel.
\newblock Foundations of intuitionistic logic.
\newblock In \emph{Studies in Logic and the Foundations of Mathematics},
  volume~44, pages 198--210. North Holland, Amsterdam, 1962{\natexlab{b}}.

\bibitem[Kreisel(1962{\natexlab{c}})]{Kreisel1962c}
G.~Kreisel.
\newblock The axiom of choice and the class of hyperarithmetic functions.
\newblock \emph{Indagationes Mathematicae}, 24:\penalty0 307--319,
  1962{\natexlab{c}}.

\bibitem[Kreisel(1963)]{Kreisel1963}
G.~Kreisel.
\newblock {Stanford Report on the foundations of analysis (with contributions
  by W.A. Howard, W. W. Tait, and R.J. Parikh)}.
\newblock Technical report, Stanford, 1963.

\bibitem[Kreisel(1965)]{Kreisel1965a}
G.~Kreisel.
\newblock Mathematical logic.
\newblock In T.~Saaty, editor, \emph{Lectures on Modern Mathematics, Vol. III},
  pages 95--195. Wiley, New York, 1965.

\bibitem[Kreisel(1967{\natexlab{a}})]{Kreisel1967a}
G.~Kreisel.
\newblock Mathematical logic: What has it done for the philosophy of
  mathematics?
\newblock In R.~Schoenman, editor, \emph{Bertrand Russell, Philosopher of the
  Century}, pages 201--272. Allen and Unwin, London, 1967{\natexlab{a}}.

\bibitem[Kreisel(1967{\natexlab{b}})]{Kreisel1967b}
G.~Kreisel.
\newblock Informal {R}igour and {C}ompleteness {P}roofs.
\newblock In I.~Lakatos, editor, \emph{{Problems in the philosophy of
  mathematics: Proceedings of the International Colloquium in the Philosophy of
  Science, London, 1965}}, pages 138--186. North-Holland, Amsterdam,
  1967{\natexlab{b}}.

\bibitem[Kreisel(1967{\natexlab{c}})]{Kreisel1967c}
G.~Kreisel.
\newblock \emph{Elements of {M}athematical {L}ogic. {M}odel theory
  \textnormal{(with Jean-Louis Krivine)}}.
\newblock Studies in Logic and the Foundations of Mathematics. North-Holland,
  Amsterdam, 1967{\natexlab{c}}.

\bibitem[Kreisel(1968{\natexlab{a}})]{Kreisel1968}
G.~Kreisel.
\newblock A survey of proof theory.
\newblock \emph{The Journal of Symbolic Logic}, 33\penalty0 (3):\penalty0
  321--388, 1968{\natexlab{a}}.

\bibitem[Kreisel(1968{\natexlab{b}})]{Kreisel1968b}
G.~Kreisel.
\newblock Lawless sequences of natural numbers.
\newblock \emph{Compositio Mathematica}, 20:\penalty0 222--248,
  1968{\natexlab{b}}.

\bibitem[Kreisel(1969{\natexlab{a}})]{Kreisel1969a}
G.~Kreisel.
\newblock Two {N}otes on the {F}oundations of {S}et-{T}heory.
\newblock \emph{Dialectica}, 23\penalty0 (2):\penalty0 93--114,
  1969{\natexlab{a}}.

\bibitem[Kreisel(1969{\natexlab{b}})]{Kreisel1969c}
G.~Kreisel.
\newblock Informal {R}igour and {C}ompleteness {P}roofs (excerpt).
\newblock In J.~Hintikka, editor, \emph{The Philosophy of Mathematics}, pages
  78--94. Oxford Univeristy Press, Oxford, 1969{\natexlab{b}}.

\bibitem[Kreisel(1970{\natexlab{a}})]{Kreisel1970}
G.~Kreisel.
\newblock Church's {T}hesis: A kind of reducibility axiom for constructive
  mathematics.
\newblock In A.~Kino, J.~Myhill, and R.~Vesley, editors, \emph{Intuitionism and
  {P}roof {T}heory}, pages 121--150. North-Holland, Amsterdam,
  1970{\natexlab{a}}.

\bibitem[Kreisel(1970{\natexlab{b}})]{Kreisel1970a}
G.~Kreisel.
\newblock Principles of proof and ordinals implicit in given concepts.
\newblock \emph{Studies in Logic and the Foundations of Mathematics},
  60:\penalty0 489--516, 1970{\natexlab{b}}.

\bibitem[Kreisel(1971)]{Kreisel1971d}
G.~Kreisel.
\newblock Observations on popular discussions of foundations.
\newblock In D.~Scott, editor, \emph{Axiomatic Set Theory}, pages 189--198.
  American Mathematical Society, Providence, 1971.

\bibitem[Kreisel(1976)]{Kreisel1976}
G.~Kreisel.
\newblock {What have we learnt from Hilbert's second problem?}
\newblock In \emph{Mathematical Developments Arising from Hilbert Problems,
  Part 2}, volume~67, pages 93--130. American mathematical society, Providence,
  R. I., 1976.

\bibitem[Kreisel(1977)]{Kreisel1977b}
G.~Kreisel.
\newblock Review of \cite{Weston1976}.
\newblock \emph{Zentralblatt f{\"u}r Mathematik}, Zbl 0341.02053:\penalty0
  31--32, 1977.

\bibitem[Kreisel(1980)]{Kreisel1980}
G.~Kreisel.
\newblock {Kurt {G}\"odel: 28 April 1906-14 January 1978}.
\newblock \emph{Biographical Memoirs of Fellows of the Royal Society},
  26:\penalty0 149--224, 1980.

\bibitem[Kreisel(1987{\natexlab{a}})]{Kreisel1987a}
G.~Kreisel.
\newblock {Church's Thesis and the Ideal of Informal Rigour}.
\newblock \emph{Notre Dame journal of formal logic}, 28\penalty0 (4):\penalty0
  499--519, 1987{\natexlab{a}}.

\bibitem[Kreisel(1987{\natexlab{b}})]{Kreisel1987b}
G.~Kreisel.
\newblock G{\"o}del's excursions into intuitionistic logic.
\newblock In P.~Weingartner and L.~Schmetterer, editors, \emph{G{\"o}del
  Remembered}, pages 67--179. Bibliopolis, Napoli, 1987{\natexlab{b}}.

\bibitem[Kreisel(1989)]{Kreisel1989}
G.~Kreisel.
\newblock Logical aspects of the axiomatic method: on their significance in
  (traditional) foundations and in some (now) common or garden varieties of
  mathematics.
\newblock In H.~Ebbinghaus, J.~Fern{\'a}ndez-Prida, M.~Garrido, D.~Lascar, and
  M.~Artalejo, editors, \emph{Logic Colloquium'87: Proceedings of the
  Colloquium Held in Granada, Spain, July 20-25, 1987}, pages 183--217. North
  Holland, 1989.

\bibitem[Kreisel and L{\'e}vy(1968)]{Kreisel1968a}
G.~Kreisel and A.~L{\'e}vy.
\newblock Reflection principles and their use for establishing the complexity
  of axiomatic systems.
\newblock \emph{Zeitschrift f\"ur Mathematische Logik und Grundlagen der
  Mathematik}, 14:\penalty0 97--142, 1968.

\bibitem[Kreisel et~al.(1975)Kreisel, Mints, and Simpson]{Kreisel1975}
G.~Kreisel, G.~Mints, and S.~Simpson.
\newblock The use of abstract language in elementary metamathematics: Some
  pedagogic examples.
\newblock In \emph{Logic Colloquium}, pages 38--131. Springer, 1975.

\bibitem[Kripke(1965)]{Kripke1965}
S.~Kripke.
\newblock {Semantical analysis of intuitionistic logic. I, Formal systems and
  recursive functions}.
\newblock In J.~N. Crossley and M.~A.~E. Dummett, editors, \emph{Formal Systems
  and Recursive Functions}. North-Holland, 1965.

\bibitem[Kripke(2019)]{Kripke2019}
S.~A. Kripke.
\newblock Free choice sequences: A temporal interpretation compatible with
  acceptance of classical mathematics.
\newblock \emph{Indagationes Mathematicae}, 30\penalty0 (3):\penalty0 492--499,
  2019.

\bibitem[Lakatos(1967)]{Lakatos1967}
I.~Lakatos.
\newblock \emph{Problems in the philosophy of mathematics: proceedings of the
  International Colloquium in the Philosophy of Science, London, 1965},
  volume~1.
\newblock North Holland, Amsterdam, 1967.

\bibitem[Langford(1942)]{Langford1942}
C.~Langford.
\newblock The notion of analysis in moore's philosophy in the philosophy of
  g.e. moore.
\newblock In P.~Schlipp, editor, \emph{The Philosophy of G.E. Moore}, pages
  321--342. Open Court, La Salle, 1942.

\bibitem[L{\'e}vy(1960)]{Levy1960}
A.~L{\'e}vy.
\newblock Axiom schemata of strong infinity in axiomatic set theory.
\newblock \emph{Pacific Journal of Mathematics}, 10\penalty0 (1):\penalty0
  223--238, 1960.

\bibitem[Mancosu(1998)]{Mancosu1998}
P.~Mancosu, editor.
\newblock \emph{From {B}rouwer to {H}ilbert: The Debate on the Foundations of
  Mathematics in the 1920s}.
\newblock Oxford University Press, Oxford, 1998.

\bibitem[Marker(2002)]{Marker2002}
D.~Marker.
\newblock \emph{Model {T}heory}, volume 217 of \emph{Graduate Texts in
  Mathematics}.
\newblock Springer-Verlag, New York, 2002.

\bibitem[Markov(1956)]{Markov1956}
A.~Markov.
\newblock On a principle of constructive mathematical logic.
\newblock In \emph{Proceedings of the Third All-Union Mathematical Congress,},
  pages 146--147, Moscow, 1956.

\bibitem[Martin(2005)]{Martin2005}
D.~A. Martin.
\newblock {G{\"o}del's conceptual realism}.
\newblock \emph{Bulletin of Symbolic Logic}, 11\penalty0 (2):\penalty0
  207--224, 2005.
\newblock ISSN 1079-8986.

\bibitem[McCarty(2008)]{McCarty2008}
C.~McCarty.
\newblock {Completeness and incompleteness for intuitionistic logic}.
\newblock \emph{Journal of Symbolic Logic}, 73\penalty0 (4):\penalty0
  1315--1327, 2008.

\bibitem[McCarty(1988)]{McCarty1988}
C.~McCarty.
\newblock Constructive validity is nonarithmetic.
\newblock \emph{The Journal of Symbolic Logic}, 53\penalty0 (4):\penalty0
  1036--1041, 1988.

\bibitem[McGee(1997)]{McGee1997}
V.~McGee.
\newblock How we learn mathematical language.
\newblock \emph{The Philosophical Review}, 106\penalty0 (1):\penalty0 35--68,
  1997.

\bibitem[Montague(1963)]{Montague1963}
R.~Montague.
\newblock Syntactical treatments of modality, with corollaries on reflexion
  principles and finite axiomatizability.
\newblock \emph{{Acta Philosophica Fennica}}, 16:\penalty0 153--167, 1963.

\bibitem[Montague(1965)]{Montague1965a}
R.~Montague.
\newblock Set theory and higher-order logic.
\newblock \emph{Studies in Logic and the Foundations of Mathematics},
  40:\penalty0 131--148, 1965.

\bibitem[Moore(1903)]{Moore1903}
G.~E. Moore.
\newblock \emph{Principia Ethica}.
\newblock Cambridge University Press, 1903.

\bibitem[Moore(1987)]{Moore1987}
G.~H. Moore.
\newblock The origins of forcing.
\newblock In \emph{Studies in Logic and the Foundations of Mathematics Logic
  Colloquium '86}. North Holland, 1987.

\bibitem[Moschovakis(2019)]{Moschovakis2019}
J.~R. Moschovakis.
\newblock {Markov's Principle, Markov's Rule} and the notion of constructive
  proof.
\newblock Forthcoming in \textsl{Intuitionism, Computation, and Proof: Selected
  themes from the research of G. Kreisel}, 2019.

\bibitem[Moschovakis(2016)]{Moschovakis2016a}
Y.~N. Moschovakis.
\newblock Hyperarithmetical sets.
\newblock In \emph{Martin Davis on Computability, Computational Logic, and
  Mathematical Foundations}, pages 107--149. Springer, 2016.

\bibitem[Moss(1971)]{Moss1971}
J.~Moss.
\newblock Kreisel's work on the philosophy of mathematics -- {I. Realism}.
\newblock In \emph{Studies in Logic and the Foundations of Mathematics},
  volume~61, pages 411--438. Elsevier, 1971.

\bibitem[Mostowski(1950)]{Mostowski1950}
A.~Mostowski.
\newblock Some impredicative definitions in the axiomatic set-theory.
\newblock \emph{Fundamenta mathematicae}, 38:\penalty0 110--124, 1950.

\bibitem[Mostowski(1948)]{Mostowski1948}
A.~Mostowski.
\newblock Proofs of non-deducibility in intuitionistic functional calculus.
\newblock \emph{The Journal of Symbolic Logic}, 13\penalty0 (4):\penalty0
  204--207, 1948.

\bibitem[Myhill(1967)]{Myhill1967}
J.~Myhill.
\newblock Notes towards an axiomatization of intuitionistic analysis.
\newblock \emph{Logique et analyse}, 9\penalty0 (35/36):\penalty0 280--297,
  1967.

\bibitem[Odifreddi(1999)]{Odifreddi1999a}
P.~Odifreddi.
\newblock Reducibilities.
\newblock In E.~Griffor, editor, \emph{Handbook of {C}omputability {T}heory},
  volume 140, pages 89--119. Elsevier, 1999.

\bibitem[Odifreddi(1983{\natexlab{a}})]{Odifreddi1983}
P.~Odifreddi.
\newblock Forcing and reducibilities.
\newblock \emph{The Journal of Symbolic Logic}, 48\penalty0 (2):\penalty0
  288--310, 1983{\natexlab{a}}.

\bibitem[Odifreddi(1983{\natexlab{b}})]{Odifreddi1983b}
P.~Odifreddi.
\newblock {Forcing and Reducibilities II}.
\newblock \emph{The Journal of Symbolic Logic}, 48\penalty0 (3):\penalty0
  724--743, 1983{\natexlab{b}}.

\bibitem[Odifreddi(1996)]{Odifreddi1996}
P.~Odifreddi, editor.
\newblock \emph{Kreiseliana. About and Around Georg Kreisel}.
\newblock A.K. Peters, Wellesley, Mass., 1996.

\bibitem[Parsons(1977)]{Parsons1977}
C.~Parsons.
\newblock What is the iterative conception of set?
\newblock In R.~E. Butts and J.~Hintikka, editors, \emph{Logic, Foundations of
  Mathematics and Computability Theory}, pages 335--367. Reidel, Dordrecht,
  1977.

\bibitem[Parsons(1990)]{Parsons1990a}
C.~Parsons.
\newblock The uniqueness of the natural numbers.
\newblock \emph{Iyyun}, 39\penalty0 (1):\penalty0 13--44, 1990.

\bibitem[Plisko(1984)]{Plisko1984}
V.~E. Plisko.
\newblock Absolute realizability of predicate formulas.
\newblock \emph{Mathematics of the {USSR}-Izvestiya}, 22\penalty0 (2):\penalty0
  291--308, 1984.

\bibitem[Poincar{\'e}(1910)]{Poincare1910}
H.~Poincar{\'e}.
\newblock {\"Uber transfinite Zahlen}.
\newblock In \emph{{Sechs Vortr{\"a}ge {\"u}ber ausgew{\"a}hlte Gegenst{\"a}nde
  aus der reinen Mathematik und mathematischen Physik: auf Einladung der
  Wolfskehl-Kommission der K{\"o}niglichen Gesellschaft der Wissenschaften
  gehalten zu G{\"o}ttingen vom 22.-28. April 1909}}, volume~4. Teubner,
  Leipzig, 1910.
\newblock Reprinted in \cite{Ewald1996}, vol. pp. 1071-1074.

\bibitem[Poincar\'e(1913)]{Poincare1913b}
H.~Poincar\'e.
\newblock Intuition and logic in mathematics.
\newblock In \emph{The foundations of science}, pages pp. 210--223. The Science
  Press, New York, 1913.
\newblock {Translated by G. Halstad}.

\bibitem[Potter(2004)]{Potter2004}
M.~Potter.
\newblock \emph{Set {T}heory and {I}ts {P}hilosophy}.
\newblock Oxford University Press, New York, 2004.

\bibitem[Quine(1970)]{Quine1970a}
W.~Quine.
\newblock \emph{Philosophy of Logic}.
\newblock Harvard University Press, Cambridge, Mass., 1970.

\bibitem[Read(1994)]{Read1994}
S.~Read.
\newblock Formal and material consequence.
\newblock \emph{Journal of Philosophical Logic}, 23\penalty0 (3):\penalty0
  247--265, 1994.

\bibitem[Reck(2012)]{Reck2012}
E.~Reck.
\newblock Carnapian explication: A case study and critique.
\newblock In \emph{Carnap's ideal of explication and naturalism}, pages
  96--116. Springer, Berlin, 2012.

\bibitem[Robinson(1965)]{Robinson1965a}
A.~Robinson.
\newblock Formalism 64.
\newblock In \emph{Logic, methodology and philosophy of science}, pages
  228--246. North-Holland, Amsterdam, 1965.

\bibitem[Robinson(1968)]{Robinson1968}
A.~Robinson.
\newblock Some thoughts on the history of mathematics.
\newblock \emph{Compositio Mathematica}, 20:\penalty0 188--193, 1968.

\bibitem[Robinson(1949)]{Robinson1949}
J.~Robinson.
\newblock Definability and decision problems in arithmetic.
\newblock \emph{The Journal of Symbolic Logic}, 14:\penalty0 98--114, 1949.

\bibitem[Rogers(1987)]{Rogers1987}
H.~Rogers.
\newblock \emph{Theory of recursive functions and effective computability}.
\newblock MIT Press, Cambridge, MA, 1987.
\newblock First edition 1967.

\bibitem[Russell(1908)]{Russell1908}
B.~Russell.
\newblock Mathematical logic as based on the theory of types.
\newblock \emph{American journal of mathematics}, 30\penalty0 (3):\penalty0
  222--262, 1908.

\bibitem[Sacks(1990)]{Sacks1990}
G.~Sacks.
\newblock \emph{Higher {R}ecursion {T}heory}.
\newblock Springer-Verlag, 1990.

\bibitem[Scott(1970)]{Scott1970a}
D.~Scott.
\newblock Constructive validity.
\newblock In \emph{{Symposium on Automatic Demonstration}}, pages 237--275,
  Berlin, 1970. Springer.

\bibitem[Scott(1961)]{Scott1961}
D.~Scott.
\newblock On {C}onstructing {M}odels for {A}rithmetic.
\newblock In \emph{Infinitistic {M}ethods}, pages 235--255. Pergamon, Oxford,
  1961.

\bibitem[Scott(1974)]{Scott1974a}
D.~Scott.
\newblock Axiomatizing {S}et {T}heory.
\newblock In T.~J. Jech, editor, \emph{Axiomatic Set Theory}, pages 207--214.
  American Mathematical Society, Providence, 1974.

\bibitem[Shapiro(1985)]{Shapiro1985a}
S.~Shapiro.
\newblock Second-order languages and mathematical practice.
\newblock \emph{The Journal of Symbolic Logic}, 50\penalty0 (3):\penalty0
  714--742, 1985.

\bibitem[Shapiro(1987)]{Shapiro1987}
S.~Shapiro.
\newblock Principles of reflection and second-order logic.
\newblock \emph{Journal of Philosophical Logic}, 16\penalty0 (3):\penalty0
  309--333, 1987.

\bibitem[Shapiro(1991)]{Shapiro1991}
S.~Shapiro.
\newblock \emph{Foundations without {F}oundationalism: {A} {C}ase for
  {S}econd-{O}rder {L}ogic}, volume~17 of \emph{Oxford Logic Guides}.
\newblock Oxford University Press, 1991.

\bibitem[Shepherdson(1952)]{Shepherdson1952}
J.~Shepherdson.
\newblock Inner {M}odels for {S}et {T}heory. {I}{I}.
\newblock \emph{Journal of Symbolic Logic}, 17:\penalty0 225--237, 1952.

\bibitem[Simpson(2009)]{Simpson2009}
S.~Simpson.
\newblock \emph{Subsystems of second order arithmetic}.
\newblock Cambridge University Press, Cambridge, second edition, 2009.

\bibitem[Smith(2011)]{Smith2011a}
P.~Smith.
\newblock Squeezing arguments.
\newblock \emph{Analysis}, 71\penalty0 (1):\penalty0 22--30, 2011.

\bibitem[Smorynski(1984)]{Smorynski1984}
C.~Smorynski.
\newblock {Lectures on nonstandard models of arithmetic}.
\newblock In G.~Lolli, G.~Longo, and A.~Marqa, editors, \emph{Logic
  {C}olloquium '82}, pages 1--70. North-Holland, Amsterdam, 1984.

\bibitem[Smullyan(1961)]{Smullyan1961}
R.~Smullyan.
\newblock \emph{{Theory of formal systems}}.
\newblock Princeton University Press, Princeton, 1961.

\bibitem[S{\o}rensen and Urzyczyn(2006)]{Sorensen2006}
M.~H. S{\o}rensen and P.~Urzyczyn.
\newblock \emph{Lectures on the Curry-Howard isomorphism}.
\newblock Elsevier, Amsterdam, 2006.

\bibitem[Spector(1955)]{Spector1955}
C.~Spector.
\newblock Recursive well-orderings.
\newblock \emph{The Journal of Symbolic Logic}, 20\penalty0 (2):\penalty0
  151--163, 1955.

\bibitem[Steel(2014)]{Steel2014}
J.~Steel.
\newblock G{\"o}del's program.
\newblock In J.~Kennedy, editor, \emph{Interpreting G{\"o}del}, pages 153--179.
  {C}ambridge {U}niversity {P}ress, Cambridge, 2014.

\bibitem[Sundholm(1983)]{Sundholm1983}
G.~Sundholm.
\newblock Constructions, proofs and the meaning of logical constants.
\newblock \emph{Journal of Philosophical Logic}, 12\penalty0 (2):\penalty0
  151--172, 1983.

\bibitem[Tait(2005)]{Tait2005}
W.~Tait.
\newblock \emph{The {P}rovenance of {P}ure {R}eason}.
\newblock Logic and Computation in Philosophy. Oxford University Press, New
  York, 2005.

\bibitem[Tarski(1936)]{Tarski1936}
A.~Tarski.
\newblock Der {W}ahrheitsbegriff in den formalisierten {S}prache.
\newblock \emph{Studia Philosophica}, 1:\penalty0 261--405, 1936.

\bibitem[Tarski and Vaught(1956)]{Tarski1956c}
A.~Tarski and R.~L. Vaught.
\newblock Arithmetical extensions of relational systems.
\newblock \emph{Compositio mathematica}, 13:\penalty0 81--102, 1956.

\bibitem[Troelstra(1969)]{Troelstra1969}
A.~S. Troelstra.
\newblock \emph{Principles of Intuitionism}.
\newblock Number~95 in Lecture Notes in Mathematics. Springer, Amsterdam, 1969.

\bibitem[Troelstra(1977{\natexlab{a}})]{Troelstra1977}
A.~S. Troelstra.
\newblock Aspects of constructive mathematics.
\newblock In J.~Barwise, editor, \emph{{Handbook of Mathematical Logic}},
  volume~90, pages 973--1052. Elsevier, 1977{\natexlab{a}}.

\bibitem[Troelstra and van Dalen(1988{\natexlab{a}})]{Troelstra1988a}
A.~S. Troelstra and D.~van Dalen.
\newblock \emph{{Constructivism in Mathematics, An Introduction, Volume 2}}.
\newblock North-Holland, Amsterdam, 1988{\natexlab{a}}.

\bibitem[Troelstra(1977{\natexlab{b}})]{Troelstra1977a}
A.~S. Troelstra.
\newblock \emph{Choice sequences: A chapter of intuitionistic mathematics}.
\newblock Oxford Logic Guides. Oxford University Press, Oxford,
  1977{\natexlab{b}}.

\bibitem[Troelstra and van Dalen(1988{\natexlab{b}})]{Troelstra1988}
A.~Troelstra and D.~van Dalen.
\newblock \emph{{Constructivism in Mathematics, An Introduction, Volume 1}}.
\newblock North-Holland, Amsterdam, 1988{\natexlab{b}}.

\bibitem[Turing(1936)]{Turing1936}
A.~Turing.
\newblock {On computable numbers, with an application to the
  \textsl{Entscheidungsproblem}}.
\newblock \emph{Proceedings of the London mathematical society}, 42\penalty0
  (2):\penalty0 230--265, 1936.

\bibitem[Turing(1939)]{Turing1939}
A.~Turing.
\newblock Systems of logic based on ordinals.
\newblock \emph{Proceedings of the London Mathematical Society}, 2\penalty0
  (1):\penalty0 161--228, 1939.

\bibitem[V{\"a}{\"a}n{\"a}nen(2012)]{Vaananen2012}
J.~V{\"a}{\"a}n{\"a}nen.
\newblock Second order logic or set theory?
\newblock \emph{Bulletin of Symbolic Logic}, 18\penalty0 (1):\penalty0 91--121,
  2012.

\bibitem[V{\"a}{\"a}n{\"a}nen and Wang(2015)]{Vaananen2015}
J.~V{\"a}{\"a}n{\"a}nen and T.~Wang.
\newblock Internal categoricity in arithmetic and set theory.
\newblock \emph{Notre Dame Journal of Formal Logic}, 56\penalty0 (1):\penalty0
  121--134, 2015.

\bibitem[V\"a\"an\"anen(2019)]{Vaananen2019}
J.~V\"a\"an\"anen.
\newblock {A extension of a theorem of Zermelo}.
\newblock \emph{The Bulletin of Symbolic Logic}, 25\penalty0 (2):\penalty0
  208--212, 2019.

\bibitem[V{\"a}{\"a}n{\"a}nen(2020)]{Vaananen2020}
J.~V{\"a}{\"a}n{\"a}nen.
\newblock Tracing internal categoricity.
\newblock \emph{Theoria}, 2020.
\newblock Online first.

\bibitem[Van~Atten(2016)]{Van-Atten2016}
M.~Van~Atten.
\newblock {Troelstra's Paradox and Markov's Principle}.
\newblock Technical Report hal-01415603, Hal Archives, 2016.

\bibitem[Van~Atten(2018)]{Van-Atten2018}
M.~Van~Atten.
\newblock {The Creating Subject, the Brouwer-Kripke Schema, and infinite
  proofs}.
\newblock \emph{Indagationes Mathematicae}, 29\penalty0 (6):\penalty0
  1565--1636, 2018.

\bibitem[van Dalen(2008)]{Dalen2008}
D.~van Dalen.
\newblock \emph{Logic and Structure}.
\newblock Springer, 4th edition, 2008.

\bibitem[Van~Dalen(1982)]{Van-Dalen1982a}
D.~Van~Dalen.
\newblock The creative subject and heyting's arithmetic.
\newblock \emph{Banach Center Publications}, 9\penalty0 (1):\penalty0 379--382,
  0 1982.

\bibitem[Vesley(1972)]{Vesley1972}
R.~Vesley.
\newblock {Choice sequences and Markov's principle}.
\newblock \emph{Compositio Mathematica}, 24\penalty0 (1):\penalty0 33--53,
  1972.

\bibitem[Walmsley(2002)]{Walmsley2002}
J.~Walmsley.
\newblock Categoricity and indefinite extensibility.
\newblock \emph{Proceedings of the Aristotelian Society}, 102:\penalty0
  239--257, 2002.

\bibitem[Walsh(2016)]{Walsh2016}
S.~Walsh.
\newblock Predicativity, the russell-myhill paradox, and church's intensional
  logic.
\newblock \emph{Journal of Philosophical Logic}, 45\penalty0 (3):\penalty0
  277--326, 2016.

\bibitem[Wang(1977)]{Wang1977}
H.~Wang.
\newblock Large sets.
\newblock In \emph{Logic, Foundations of Mathematics, and Computability
  Theory}, pages 309--333. Springer, 1977.

\bibitem[Weston(1977)]{Weston1977}
T.~Weston.
\newblock The continuum hypothesis is independent of second-order
  $\mathsf{ZF}$.
\newblock \emph{Notre Dame Journal of Formal Logic}, 18\penalty0 (3):\penalty0
  499--503, 1977.

\bibitem[Weston(1974)]{Weston1974}
T.~Weston.
\newblock \emph{The {C}ontinuum {H}ypothesis: {I}ndependence and
  {T}ruth-{V}alue}.
\newblock PhD thesis, Massachusetts Institute of Technology, 1974.

\bibitem[Weston(1976)]{Weston1976}
T.~Weston.
\newblock Kreisel, the {C}ontinuum {H}ypothesis and {S}econd {O}rder {S}et
  {T}heory.
\newblock \emph{Journal of Philosophical Logic}, 5\penalty0 (2):\penalty0
  281--298, 1976.

\bibitem[Weyl(1918)]{Weyl1918}
H.~Weyl.
\newblock \emph{Das Kontinuum. Kritische Untersuchungen {\"u}ber die Grundlagen
  der Analysis.}
\newblock Verlag von Veit \& Comp., Leipzig, 1918.

\bibitem[Zermelo(1930)]{Zermelo1930a}
E.~Zermelo.
\newblock \"{U}ber {G}renzzahlen und {M}engenbereiche: {N}eue {U}ntersuchungen
  \"uber die {G}rundlagen der {M}engenlehre.
\newblock \emph{Fundamenta Mathematicae}, 16:\penalty0 29--47, 1930.
\newblock Translated in \cite{Ewald1996} vol.~2 pp.~1219-1233.

\bibitem[Zermelo(2010)]{Zermelo1930d}
E.~Zermelo.
\newblock {\textsl{Bericht} and die Notgemeinschaft der Deutschen Wissenschaft
  \"uber meine Forschungen betreffend die \textsl{Grundlagen der Mathematik}}.
\newblock In H.~Ebbinghaus and A.~Kanamori, editors, \emph{{Ernst Zermelo:
  Collected Works}}, pages 432--443. Springer, 2010.

\end{thebibliography}
}

\end{document}